\documentclass[11pt]{amsart}
\usepackage{color}
\usepackage{psfrag}
\usepackage{fullpage}
\usepackage{graphicx}
\usepackage{hyperref}
\usepackage{subfigure}
\usepackage{enumerate}
\newtheorem{theorem}{Theorem}[section]
\newtheorem{remark}[theorem]{Remark}

\title{Enforcing the non-negativity constraint and maximum principles 
for diffusion with decay on general computational grids}
\author{H.~Nagarajan} \author{K.~B.~Nakshatrala}
\address{Harsha Nagarajan, Graduate student, Department of Mechanical Engineering, Texas A\&M University, College Station, Texas 77843.}
\address{Correspondence to: Kalyana Babu Nakshatrala, Department of Mechanical Engineering, 216 Engineering/Physics Building, Texas A\&M University, College Station, Texas 77843. TEL:+1-979-845-1292}
\email{knakshatrala@tamu.edu}

\date{\today}

\begin{document}

\begin{abstract}
  In this paper, we consider anisotropic diffusion with decay, which takes the form $\alpha(\mathbf{x}) 
  c(\mathbf{x}) - \mathrm{div}[\mathbb{D}(\mathbf{x})\mathrm{grad}[c(\mathbf{x})]] = f(\mathbf{x})$ with 
  decay coefficient $\alpha(\mathbf{x}) \geq 0$, and diffusivity coefficient $\mathbb{D}(\mathbf{x})$ to 
  be a second-order symmetric and positive definite tensor. It is well-known that this particular equation 
  is a second-order elliptic equation, and satisfies a maximum principle under certain regularity 
  assumptions. However, the finite element implementation of the classical Galerkin formulation for 
  both anisotropic and isotropic diffusion with decay does not respect the maximum principle. Put 
  differently, the classical Galerkin formulation violates the discrete maximum principle for diffusion 
  with decay even on structured computational meshes. 
  
  We first show that the numerical accuracy of the classical Galerkin formulation deteriorates 
  dramatically with an increase in $\alpha$ for isotropic media and violates the discrete maximum 
  principle. However, in the case of isotropic media, the extent of violation decreases with the mesh 
  refinement. We then show that, in the case of anisotropic media, the classical Galerkin formulation 
  for anisotropic diffusion with decay violates the discrete maximum principle even at lower values of 
  decay coefficient and does not vanish with mesh refinement. We then present a methodology for 
  enforcing maximum principles under the classical Galerkin formulation for anisotropic diffusion 
  with decay on general computational grids using optimization techniques. Representative 
  numerical results (which take into account anisotropy and heterogeneity) are presented to 
  illustrate the performance of the proposed formulation. 
\end{abstract}
\keywords{maximum principle; discrete maximum principle; non-negative solutions; 
  convex quadratic programming; anisotropic diffusion with decay; general computational grids}

\maketitle 


\section{INTRODUCTION}
\label{Sec:Decay_Introduction}

In this paper we consider \emph{heterogeneous anisotropic diffusion with decay}, which 
takes the form: $\alpha(\mathbf{x}) c(\mathbf{x}) - \mathrm{div}[\mathbb{D}(\mathbf{x}) 
\mathrm{grad}[c(\mathbf{x})]] = f(\mathbf{x})$ with $\alpha(\mathbf{x}) \geq 0$ and 
$\mathbb{D}(\mathbf{x})$ is a symmetric and positive definite second-order tensor. 
This equation is a linear second-order elliptic partial differential equation 
\cite{Gilbarg_Trudinger}. There are many important problems in Mathematical Physics 
which give rise to this equation \cite{Tikhonov_Samarskii}. Also, this equation 
arises in numerical methods and mathematical analysis of transient problems 
\cite{Ladyzhenskaya_Mathematical_Physics}. Some of these cases include:
\begin{enumerate}[(a)]
\item For certain gases, the diffusion process is accompanied by a decay of the molecules 
  of the diffusing gas, and the decay is proportional to the concentration of the gas. 
  Such a phenomenon can be modeled as a diffusion equation with decay. 
\item For certain problems, the governing equation of diffusion in a moving domain can be 
  transformed into a diffusion equation with decay. 
\item Application of the method of horizontal lines to the transient diffusion equation 
  (which is a linear parabolic partial differential equation) gives rise to a diffusion 
  equation with decay. 
\end{enumerate}

\subsection{Maximum principles and discrete maximum principles}
From the theory of partial differential equations, it is well-known that the diffusion equation 
with decay satisfies a maximum principle under appropriate regularity assumptions. In some 
cases (but not always) the physically important condition that the concentration is non-negative 
is a direct consequence of a maximum principle. \emph{It is important to note that the classical 
maximum principle for diffusion with decay is different from the classical maximum principle for 
pure diffusion equation (see Theorem \ref{Theorem:Decay_maximum_principle} and Remark 
\ref{Remark:Decay_pure_diffusion} in this paper).} 

It is imperative that predictive numerical simulations employ accurate and reliable 
discretization methods. The resulting discrete systems must inherit or mimic fundamental 
properties of continuous systems. The non-negative constraint and maximum principles are 
some of the essential properties of diffusion-type equations. However, it is well-known (and 
also discussed below) that many numerical formulations (including the popular ones) may 
not give non-negative solutions or satisfy maximum principles for these types of equations 
on general computational grids. \emph{Another point to note is that the satisfaction of 
maximum principles and the non-negative constraint by a numerical formulation will be altered 
by the presence of the decay term.} (That is, the conditions under which a numerical 
formulation satisfies maximum principles and the non-negative constraint for pure diffusion 
can be different from those for diffusion with decay.) This leads us to discrete maximum 
principles.

The discrete analogy of a maximum principle is commonly referred to as a discrete 
maximum principle (DMP). Some factors that affect discrete maximum principles 
are: numerical formulation, mesh size, element type, nature of the computational 
domain (e.g., presence/absence of holes), and properties of the medium -- decay 
coefficient, diffusivity coefficient, anisotropy, and heterogeneity.  

\subsection{Prior numerical works}
Numerous numerical formulations have been developed for both isotropic and anisotropic 
diffusion equations. These formulations are based on finite difference methods 
\cite{Morton_Mayers,Strikwerda}, finite volume method \cite{Patankar,
Eigestad_Aavatsmark_Espedal_CompGeoSci_2002_v6_p381,Droniou_Eymard_NumerMath_2006_v105_p35}, 
finite element method \cite{Braess,Hughes}, mixed method \cite{Brezzi_Fortin,Raviart_Thomas_MAFEM_1977_p292,
Nedelec_NumerMath_1980_v35_p315,Nedelec_NumerMath_1986_v50_p57,
Brezzi_Douglas_Marini_NumerMath_1985_v47_p217,
Brezzi_Douglas_Durran_Marini_NumerMath_1987_v51_p237,
Brezzi_Douglas_Fortin_Marini_MMNA_1987_v21_p581,
Masud_Hughes_CMAME_2002_v191_p4341,Nakshatrala_Turner_Hjelmstad_Masud_CMAME_2006_v195_p4036}, 
discontinuous Galerkin method \cite{Arnold_Brezzi_Cockburn_Marini_SIAMJNA_2002_v39_p1749,
Hughes_Masud_Wan_CMAME_2006_v195_p3347,Brezzi_Hughes_Marini_Masud_SIAMJSC_2005_v22_p119}, 
spectral element method \cite{Karniadakis_Sherwin}, and mimetic method 
\cite{Hyman_Morel_Shashkov_Steinberg_CompGeoSci_2002_v6_p333,
Kuznetsov_Lipnikov_Shashkov_CompGeoSci_2004_v8_p301,Brezzi_Lipnikov_Shashkov_SIAMJNA_2005_v43_p1872,
Lipnikov_Shashkov_Svyatskiy_JCP_2006_v211_p473}. Most of these methods can be extended to diffusion 
with decay. However, none of the aforementioned \emph{specific} formulations satisfy maximum principles 
(for both pure diffusion equation, and diffusion with decay). 

Lately, there is a surge in research activity on enforcing maximum principles, especially 
for diffusion-type equations. However, these earlier works differ from the proposed 
formulation as they have one or more following limitations: 
\begin{enumerate}[(a)]
\item The studies did not consider anisotropy and heterogeneity. It should be noted that 
  developing numerical formulations that satisfy for isotropic diffusion is much easier 
  than anisotropic diffusion, and there are practical solutions to satisfy maximum principles 
  under the classical Galerkin formulation for homogeneous isotropic medium. These include: 
  \begin{enumerate}[(i)]
  \item Any one-dimensional mesh with linear elements satisfies maximum principles under 
    the classical Galerkin formulation. 
  \item Any mesh with acute-angled triangles or (even right-angled triangles) will satisfy 
    maximum principles. Under certain milder restriction, a Delaunay mesh will also satisfy 
    maximum principles. Now, with advances in computational geometry, software packages are 
    available which can produce Delaunay meshes for reasonably complex geometries. For example, 
    \textsf{CGAL} \cite{cgal}, \textsf{Qhull} \cite{Barber96thequickhull,qhull}, \textsf{Triangle} 
    \cite{shewchuk96b}.
  \item A mesh with rectangular elements with some restrictions on the aspect ratio satisfies the 
  discrete maximum principle \cite{Christie_Hall_IJNME_1984_v20_p549}.  In particular, a mesh 
  with square elements satisfy the discrete maximum principle. 
  \end{enumerate}
  None of the above conditions (except Condition (i)) ensure the satisfaction of maximum 
  principles and the non-negative constraint if the diffusivity tensor is anisotropic. 
  \item The studies did not consider the effect of decay. The decay term affects the 
    classical maximum principle of second-order elliptic partial differential equation 
    (see Theorem \ref{Theorem:Decay_maximum_principle}). Moreover, the decay terms 
    alters the conditions under which a formulation satisfies maximum principles. 
  \item The studies did not consider general computational grids, but instead derived 
    conditions on the mesh and on the properties of the medium. They limited their 
    studies to structured grids (rectangular elements, acute-angled triangles).
\end{enumerate}

We now briefly discuss some of the important works on discrete maximum principles. The earlier 
works on discrete maximum principles are from the finite difference literature. Some of these 
notable works are \cite{Varga_SIAMJNA_1966_v3_p355,Ciarlet_AequationesMath_1970_v4_p338}. 
It is important to note that these studies did not consider anisotropy, and general computational grids. 
In References \cite{Hohn_Mittelmann_Computing_1981_v27_p145,Yanik_ComputMathAppl_1987_v14_p459,
Yanik_ComputMathAppl_1989_v17_p1431,Vejchodsky_Solin_MathComput_2007_v76_p1833}, sufficient 
conditions are derived for higher-order elements to satisfy discrete maximum principles, but 
the studies are restricted to one-dimensional problems or isotropic diffusion. 
Ciarlet and Raviart \cite{Ciarlet_Raviart_CMAME_1973_v2_p17} considered \emph{isotropic} 
diffusion with decay under the classical Galerkin formulation. The main goal of Reference 
\cite{Ciarlet_Raviart_CMAME_1973_v2_p17} is to get restrictions on the mesh to satisfy 
maximum principles, and not a methodology that works on general computational grids. 
Herrera and Valocchi \cite{Herrera_Valocchi_GW_2006_v44_p803} have employed \emph{flow-oriented 
  derivatives} to enforce the non-negative constraint. However, the methodology is limited in 
scope as it is restricted to rectangular grids, and a special form of the diffusivity 
tensor.  
References \cite{LePotier_CRM_2005_v341_p787,Lipnikov_Shashkov_Svyatskiy_Vassilevski_JCP_2007_v227_p492,
Yuan_Sheng_JCP_2008_v227_p6288,Lipnikov_Svyatskiy_Vassilevski_JCP_2009_v228_p703} addressed (pure) 
anisotropic diffusion using \emph{finite volume techniques}. All these papers are some variants of 
the idea proposed by LePotier \cite{LePotier_CRM_2005_v341_p787}, which is to choose the location 
of sampling points for the concentration in each cell in such a way to meet the non-negative 
constraint. The methodology (which is proposed under the finite volume method) cannot be easily 
modified to fit into the framework offered by the finite element method (at least, not in the 
present form presented in these references), and till date, there is no extension of this idea 
to the finite element method. Some notable works on discrete maximum principles and monotonicity 
are in the \emph{Multi-Point Flux Approximation} (MPFA) literature 
\cite{Mlacnik_Durlofsky_JCP_2006_v216_p337,Nordbotten_Aavatsmark_Eigestad_NumerMath_2007_v106_p255,
Keilegavlen_Nordbotten_Aavatsmark_AML_2009_v22_p1178}, and these works considered logically 
rectangular grids, or derived restrictions on the mesh and medium properties. 

Liska and Shashkov \cite{Liska_Shashkov_CiCP_2008_v3_p852} proposed a non-negative formulation 
for pure anisotropic diffusion equation based on \emph{conservative finite difference methods} \cite{Shashkov}. 
Nakshatrala and Valocchi \cite{Nakshatrala_Valocchi_JCP_2009_v228_p6726} have extended the variational 
multiscale and lowest order Raviart-Thomas \emph{mixed} formulations to produce non-negative solutions 
based on optimization techniques. Also, Reference \cite[Appendix]{Nakshatrala_Valocchi_JCP_2009_v228_p6726} 
discusses various conditions to satisfy the non-negative constraint. Another interesting work is by Burman and 
Ern \cite{Burman_Ern_CRM_2004_v338_p641} who have derived a \emph{nonlinear stabilized} Galerkin 
formulation that satisfies a discrete maximum principle on general grids but they considered isotropic diffusion. 
Other recent works on discrete maximum principle include \cite{Korotov_Krizek_Neittaanmaki_MathComp_2000_v70_p107,
Vanselow_AppMath_2001_v46_p13,Karatson_Korotov_NumerMath_2005_v99_p669,
Karatson_Korotov_JCAM_2006_v192_p75}, and all these works 
focused on getting restrictions on computational meshes to satisfy maximum principles. As discussed 
in Reference \cite{Nakshatrala_Valocchi_JCP_2009_v228_p6726}, the idea of getting restrictions on 
the mesh and medium properties in the case of anisotropic medium as the conditions are stringent, 
and in some cases a mesh may not even exist. 
This paper is an extension of the ideas presented in References 
\cite{Nakshatrala_Valocchi_JCP_2009_v228_p6726,Liska_Shashkov_CiCP_2008_v3_p852}.

\subsection{Main contributions of this paper}
The main contribution of the paper is to present a robust methodology for enforcing 
maximum principles and the non-negative constraint for \emph{anisotropic diffusion 
with decay}. The methodology is applicable for general computational grids with low-order 
finite elements. We also derive a (theoretical) sufficient condition for uniform computational 
meshes under which the classical Galerkin formulation for diffusion with decay satisfies the 
maximum principle for one-dimensional problems. 

\subsection{An outline and symbolic notation used in this paper}
The remainder of this paper is organized as follows. In Section \ref{Sec:Decay_Governing_equations}, 
we present governing equations for anisotropic diffusion with decay, and clearly outline the problem 
statement. In Section \ref{Sec:Decay_Nonnegative}, we present a methodology for enforcing the 
non-negativity constraint and maximum principles for anisotropic diffusion with decay on general 
computational grids. In Sections \ref{Sec:Decay_NR}, we illustrate the performance of the proposed 
formulation using representative numerical examples. Finally, conclusions are drawn in Section 
\ref{Sec:Decay_Conclusions}.

The symbolic notation adopted in this paper is as follows. Throughout this paper repeated 
indices do not imply summation. (That is, we do not employ Einstien's summation convention.) 
We shall make a distinction between vectors in the continuum and finite element settings. 
Similarly, we make a distinction between second-order tensors in the continuum setting 
versus matrices in the context of the finite element method. The continuum vectors are 
denoted by lower case boldface normal letters, and the second-order tensors will be denoted 
using \LaTeX{} blackboard font (for example, vector $\mathbf{x}$ and second-order tensor 
$\mathbb{D}$). In the finite element context, we shall denote the vectors using lower case 
boldface italic letters, and the matrices are denoted using upper case boldface italic 
letters. For example, vector $\boldsymbol{v}$ and matrix $\boldsymbol{K}$. Other notational 
conventions adopted in this paper are introduced as needed.

\section{GOVERNING EQUATIONS AND PROBLEM STATEMENT}
\label{Sec:Decay_Governing_equations}
Let $\Omega \subset \mathbb{R}^{nd}$ be a bounded open set, where ``$nd$'' denotes the number 
of spatial dimensions.  The boundary is denoted by $\partial \Omega$, which is assumed to be 
piecewise smooth. A spatial point is denoted by $\mathbf{x} \in \Omega$. The gradient and 
divergence with respect to $\mathbf{x}$ are denoted by $\mathrm{grad}[\cdot]$ and $\mathrm{div}
[\cdot]$, respectively. The concentration of a chemical species is denoted by $c(\mathbf{x})$. 
The boundary is divided into two parts: $\Gamma^{\mathrm{D}}$ and $\Gamma^{\mathrm{N}}$ such that 
$\Gamma^{\mathrm{D}} \cup \Gamma^{\mathrm{N}} = \partial \Omega$ and $\Gamma^{\mathrm{D}} \cap 
\Gamma^{\mathrm{N}} = \emptyset$. $\Gamma^{\mathrm{D}}$ is that part of the boundary on which 
Dirichlet boundary condition is prescribed, and $\Gamma^{\mathrm{N}}$ is the part of the 
boundary on which Neumann boundary condition is prescribed. The unit outward normal to 
the boundary is denoted by $\mathbf{n}(\mathbf{x})$. 
The governing equations for anisotropic diffusion with decay can be written as follows:
\begin{subequations}
  \label{Eqn:Decay_Governing_equations}
  \begin{align}
    \label{Eqn:Decay_diffusion}
    &\alpha(\mathbf{x}) c(\mathbf{x})-\mathrm{div}[\mathbb{D}(\mathbf{x}) \mathrm{grad}
    [c(\mathbf{x})]] = f(\mathbf{x}) \quad \mathrm{in} \; \Omega  \\
    & c(\mathbf{x}) = c^{\mathrm{p}}(\mathbf{x}) \quad \mathrm{on} \; 
    \Gamma^{\mathrm{D}} \\
    & \mathbf{n}(\mathbf{x}) \cdot \mathbb{D}(\mathbf{x}) \mathrm{grad}[c(\mathbf{x})] 
    = t^{\mathrm{p}}(\mathbf{x}) \quad \mathrm{on} \; \Gamma^{\mathrm{N}} 
  \end{align}
\end{subequations}
where $\alpha(\mathbf{x}) \geq 0$ is the decay coefficient, $\mathbb{D}(\mathbf{x})$ is the diffusivity 
tensor, $f(\mathbf{x})$ is the volumetric source/sink, $c^{\mathrm{p}}(\mathbf{x})$ is the prescribed 
concentration on the boundary, and $t^{\mathrm{p}}(\mathbf{x})$ is the prescribed flux on the boundary. 
The diffusivity tensor is symmetric, and assumed to be bounded and uniformly elliptic. That is, there 
exists two constants $\mathrm{0 < \xi_1 \leq \xi_2 < +\infty}$ such that 
\begin{align}
\label{Eqn:Decay_positive_definiteness_D}
  \xi_1 \mathbf{y}^{\mathrm{T}} \mathbf{y}\leq \mathbf{y}^{\mathrm{T}} \mathbb{D}
  (\mathbf{x}) \mathbf{y} \leq \xi_2 \mathbf{y}^{\mathrm{T}} \mathbf{y} \quad 
  \forall \mathbf{x} \in \Omega \; \mathrm{and} \; \forall \mathbf{y} \in 
  \mathbb{R}^{nd}
\end{align}
Equation \eqref{Eqn:Decay_Governing_equations} is a second-order elliptic partial differential 
equation, and from the theory of partial differential equations, it is known to satisfy the 
following maximum principle: 
\begin{theorem}[maximum principle]
  \label{Theorem:Decay_maximum_principle}
  Let $c(\mathbf{x})$ $\in C^{2}(\Omega) \cap C(\bar{\Omega})$, and $\alpha(\mathbf{x}) \in 
  C^{0}(\bar{\Omega})$ with $\alpha(\mathbf{x}) \geq 0$. In addition, $\mathrm{div}[\mathbb{D}
  (\mathbf{x})]$ exists and is bounded in $\Omega$. If $\alpha(\mathbf{x}) c(\mathbf{x}) - 
  \mathrm{div}[\mathbb{D}(\mathbf{x}) \mathrm{grad}[c(\mathbf{x})]] \geq 0$ in $\Omega$ 
  then $c(\mathbf{x})$ satisfies the following equation: 
  \begin{align}
    \label{Eqn:Decay_Maximum_principle}
    \mathop{\mathrm{min}}_{\mathbf{x} \in \bar{\Omega}} \; c(\mathbf{x}) \geq 
    \mathop{\mathrm{min}}_{\mathbf{x} \in \partial \Omega} \; c^{-}(\mathbf{x})
  \end{align}
  where 
  \begin{align}
    c^{-}(\mathbf{x}) := \min(c(\mathbf{x}),0)
  \end{align}
\end{theorem}
\noindent A proof to the above theorem can be found in any standard books on partial differential 
equations (e.g., References \cite{Fraenkel_maximum_principles,Gilbarg_Trudinger,Protter_Weinberger}). 
Few remarks about the above theorem and its implications are in order.
\begin{remark}
\label{Remark:Decay_max_min}
  If $c(\mathbf{x})$ satisfies $\alpha(\mathbf{x}) c(\mathbf{x}) - \mathrm{div}[\mathbb{D}
  (\mathbf{x})\mathrm{grad}[c(\mathbf{x})]] \leq 0$ in $\Omega$ (and the remaining conditions 
  in Theorem \ref{Theorem:Decay_maximum_principle} hold) then $c(\mathbf{x})$ satisfies the 
  following equation: 
  \begin{align}
    \mathop{\mathrm{max}}_{\mathbf{x} \in \bar{\Omega}} \; c(\mathbf{x}) 
    \leq \mathop{\mathrm{max}}_{\mathbf{x} \in \partial \Omega} \; c^{+}(\mathbf{x})
  \end{align}
  where 
  \begin{align}
    c^{+}(\mathbf{x}) := \max(c(\mathbf{x}),0)
  \end{align}
\end{remark}
\begin{remark}
  If $\alpha(\mathbf{x}) < 0$ then the equation \eqref{Eqn:Decay_diffusion} is called 
  Helmholtz equation. It should be noted that Helmholtz equation does not satisfy a 
  maximum principle similar to Theorem \ref{Theorem:Decay_maximum_principle}, and a 
  counterexample can be found in Reference \cite{Protter_Weinberger}. This implies 
  that the condition $\alpha(\mathbf{x}) \geq 0$ in Theorem \ref{Theorem:Decay_maximum_principle} 
  cannot be relaxed. 
\end{remark}
\begin{remark}
  It should be noted that one can find in the literature maximum principles even when 
  $c(\mathbf{x})$ does not belong to $C^{2}(\Omega)$ (and even when $c(\mathbf{x})$ 
  is only measurable, for example see Reference \cite{Trudinger_MathZ_1977_v156_p291}). 
  A detailed discussion of such results is beyond the scope of this paper, and is not central 
  to the development of the proposed numerical formulation. An interested reader on maximum 
  principles under weaker conditions can refer to \cite{Pucci_Serrin,Protter_Weinberger,
  Gilbarg_Trudinger,Fraenkel_maximum_principles} and references therein.
\end{remark}
\begin{remark}
  \label{Remark:Decay_pure_diffusion}
  For the case of pure diffusion (i.e., $\alpha(\mathbf{x}) = 0$) we have the following 
  maximum principle. Let $c(\mathbf{x})$ $\in C^{2}(\Omega) \cap C(\bar{\Omega})$, and 
  $\mathrm{div}[\mathbb{D}(\mathbf{x})]$ exists and bounded in $\Omega$. If $-\mathrm{div}
  [\mathbb{D}(\mathbf{x}) \mathrm{grad}[c(\mathbf{x})]] \geq 0$ in $\Omega$ then $c(\mathbf{x})$ 
  satisfies 
  \begin{align}
    \mathop{\mathrm{min}}_{\mathbf{x} \in \bar{\Omega}} \; c(\mathbf{x}) =
    \mathop{\mathrm{min}}_{\mathbf{x} \in \partial \Omega} \; c(\mathbf{x})
  \end{align}
\end{remark}

\begin{remark}
  It is important to note the difference in the maximum principles for pure diffusion 
  (which is given in Remark \ref{Remark:Decay_pure_diffusion}) and diffusion with decay 
  (which is given by Theorem \ref{Theorem:Decay_maximum_principle}). In the case of 
  $\alpha(\mathbf{x}) \geq 0$ (that is, diffusion with decay), the ``\emph{non-negative 
    minimum}'' occurs on the boundary, whereas in the case of $\alpha(\mathbf{x}) = 0$ 
  (that is, pure diffusion) the maximum principle says that the minimum occurs on the 
  boundary. 
\end{remark}

\subsection{Consequence of maximum principles}
Maximum principles have important mathematical consequences in the study of partial differential 
equations and physical implications in modeling. Maximum principles are often employed in proving 
well-posedness (in particular, uniqueness of solution), and obtaining point-wise estimates. For 
example, for Poisson's equation (which is a second-order elliptic partial differential equation) 
the uniqueness of solution is a direct consequence of the maximum principle \cite{McOwen}. 
To illustrate an important physical implication, let us apply the maximum principle outlined above 
to the transient diffusive system given by equation \eqref{Eqn:Decay_Governing_equations}. We shall 
assume that $\Gamma^{\mathrm{D}} = \partial \Omega$ (that is, we prescribe Dirichlet boundary conditions 
on the whole boundary). If $f(\mathbf{x}) \geq 0$ (i.e., we have volumetric source), and $c^{\mathrm{p}}
(\mathbf{x}) \geq 0$ (i.e., we have non-negative prescribed Dirichlet boundary conditions on the whole 
boundary); then from the maximum principle it can be inferred that the quantity $c(\mathbf{x})$ is 
non-negative in the whole domain. That is, 
\begin{align}
  c(\mathbf{x}) \geq 0 \quad \forall \mathbf{x} \in \bar{\Omega}
\end{align}
Now, the question is whether a given numerical formulation gives non-negative solutions if 
the prescribed data on the boundary is non-negative and the prescribed forcing function is 
a source. Also, whether a chosen numerical formulation gives solutions that are in accordance 
with maximum principles. This leads us to the problem statement and the approach taken in this paper. 

\begin{remark}
  Under certain conditions (on the forcing function and boundary conditions), the non-negative 
  constraint can be a special case of a maximum principle as shown above. However, it should 
  be noted that, in general, the non-negative constraint can be an independent result, and need 
  not be a consequence of any known maximum principle. For example, one can construct a 
  simple problem in which the non-negative constraint is not a consequence of the maximum 
  principle given in Theorem \ref{Theorem:Decay_maximum_principle}. To wit, one can have 
  a forcing function that is a source in some region and a sink in some other region of 
  the domain. For this case, the conditions given in Theorem \ref{Theorem:Decay_maximum_principle} 
  are not met, but still one may have the non-negative constraint on the concentration of the 
  diffusant. 
\end{remark}

\subsection{Problem statement and our approach}
The problem statement can be written as follows: develop a finite element formulation 
for anisotropic diffusion with decay that satisfies the non-negative constraint and 
maximum principles on general computational grids for low-order finite elements. 

\begin{figure}[h]
  \centering	
  \includegraphics[scale=0.85]{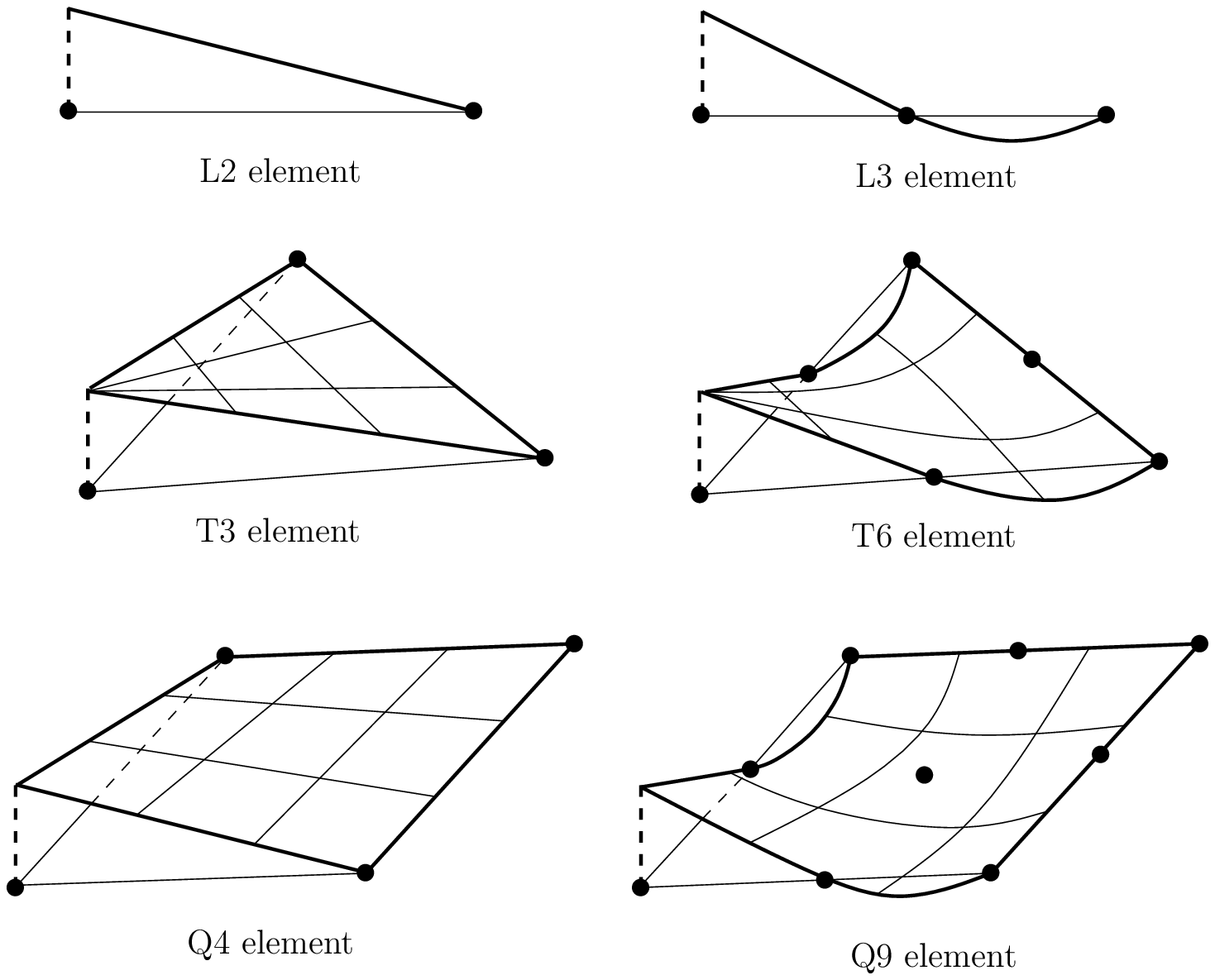}
  \caption{This figure illustrates how the proposed methodology of 
    enforcing the non-negative constraint and maximum principles 
    works for low-order finite elements like two-node linear element 
    (L2), three-node triangular element (T3), four-node quadrilateral 
    element (Q4). The proposed methodology does not  work for high-order 
    elements like three-node quadratic element (L3), six-node triangular 
    element (T6), nine-node quadrilateral element (Q9). In all the cases, 
    the nodal concentrations are non-negative. For low-order elements, 
    non-negative nodal concentrations ensures that the solution is 
    non-negative within the whole finite element. In the case of 
    high-order finite elements, enforcing non-negative nodal 
    concentrations does not imply non-negative concentration 
    throughout the element domain.} \label{Fig:Decay_T3T6_plot}
\end{figure}

The proposed methodology is based on the following key idea. We start with the finite 
element formulation of the classical Galerkin formulation, which has a variational 
statement. To this variational statement, we augment the bounds on the nodal concentrations 
given by the maximum principle. The resulting problem belongs to convex quadratic programming, 
and is solved by the active-set strategy. The proposed methodology works for all low-order 
finite elements (e.g., two-node linear element, three-node triangular element, four-node 
quadrilateral element, four-node tetrahedron element, and eight-node brick element) as 
nodal concentrations satisfying the maximum principle ensure that the maximum principle 
is met throughout the computational domain. The proposed methodology, in general, does 
not work for high-order elements as illustrated in Figure \ref{Fig:Decay_T3T6_plot}.

\section{WEAK FORMULATION AND DISCRETE MAXIMUM PRINCIPLE}
\label{Sec:Decay_Nonnegative}
Herein, we employ the classical (single-field) Galerkin formulation. We shall 
define the following function spaces:
\begin{subequations}
  \begin{align}
    \mathcal{P} &:= \left\{c(\mathbf{x}) \in H^{1}(\Omega) \; \big| \; c(\mathbf{x}) = 
      c^{\mathrm{p}}(\mathbf{x}) \; \mathrm{on} \; \Gamma^{\mathrm{D}}\right\} \\
    \mathcal{Q} &:= \left\{w(\mathbf{x}) \in H^{1}(\Omega) \; \big| \; w(\mathbf{x}) = 0 
      \; \mathrm{on} \; \Gamma^{\mathrm{D}}\right\} 
  \end{align}
\end{subequations}
where $H^{1}(\Omega)$ is a standard Sobolev space \cite{Brezzi_Fortin}. For weak solutions, 
we can relax the regularity requirement on the diffusivity tensor $\mathbb{D}(\mathbf{x})$. 
We shall assume that  each component of $\mathbb{D}(\mathbf{x})$ is square integrable, which 
is equivalent to saying that 
\begin{align}
  \int_{\Omega} \mathrm{tr}[\mathbb{D}(\mathbf{x})^{T} \mathbb{D}(\mathbf{x})] \; 
  \mathrm{d} \Omega < +\infty
\end{align}
where $\mathrm{tr}[\cdot]$ is the standard trace operator \cite{Chadwick} used in 
continuum mechanics. The classical Galerkin formulation for anisotropic diffusion 
with decay \eqref{Eqn:Decay_Governing_equations} reads: Find $c(\mathbf{x}) \in 
\mathcal{P}$ such that 
\begin{align}
  \label{Eqn:Decay_single_field_formulation}
  \mathcal{B}(w;c) = L(w) \quad \forall w(\mathbf{x}) \in \mathcal{Q}
\end{align}
where the bilinear form and linear functional are, respectively, defined as 
\begin{subequations}
  \label{Eqn:Decay_functionals_B_L}
  \begin{align}
    \mathcal{B}(w;c) &:= \int_{\Omega} \mathrm{grad}[w(\mathbf{x})] \cdot \mathbb{D}(\mathbf{x}) 
    \mathrm{grad}[c(\mathbf{x})] \; \mathrm{d} \Omega + \int_{\Omega} w(\mathbf{x})  
    \mathrm{\alpha}(\mathbf{x}) c(\mathbf{x}) \; \mathrm{d} \Omega
    \\
    L(w) &:= \int_{\Omega} w(\mathbf{x}) \; f(\mathbf{x}) \; \mathrm{d} \Omega + 
    \int_{\Gamma^{\mathrm{N}}} w(\mathbf{x}) \; t^{\mathrm{p}}(\mathbf{x}) \; \mathrm{d} \Gamma 
  \end{align}
\end{subequations}
It is well-known that the above weak form \eqref{Eqn:Decay_single_field_formulation} 
is equivalent to the following variational statement
\begin{align}
  \mathop{\mathrm{minimize}}_{c(\mathbf{x}) \in \mathcal{P}} \quad 
  \frac{1}{2} \mathcal{B}(c;c) - L(c)
\end{align}   

It may not be possible, in general, to obtain analytical solutions for equations 
\eqref{Eqn:Decay_single_field_formulation}--\eqref{Eqn:Decay_functionals_B_L} 
especially for realistic problems with complex geometries. In such situations 
one may have to resort to numerical solutions. Herein we employ the Finite 
Element Method (FEM). Let the domain $\Omega$ be decomposed into ``$Nele$'' 
non-overlapping open element subdomains. That is, 
\begin{align}
  \bar{\Omega} = \bigcup_{e = 1}^{Nele} \bar{\Omega}^{e}
\end{align}
where a superposed bar denotes the set closure. The boundary of $\Omega^e$ is denoted as 
$\partial \Omega^{e} := \bar{\Omega}^{e} - \Omega^{e}$. For a non-negative integer $m$, 
$\mathbb{P}^{m}(\Omega^{e})$ denotes the linear vector space spanned by polynomials up to 
$m$-th order defined on the subdomain $\Omega^{e}$. We shall define the following finite 
dimensional vector spaces of $\mathcal{P}$ and $\mathcal{Q}$:
\begin{subequations}
  \begin{align}
    \mathcal{P}^{h} &:= \left\{c^{h}(\mathbf{x}) \in \mathcal{P} \; \big| \; c^{h}(\mathbf{x}) 
      \in C^{0}(\bar{\Omega}), c^{h}(\mathbf{x}) \big|_{\Omega^e} \in \mathbb{P}^{k}(\Omega^{e}), 
      e = 1, \cdots, Nele \right\} \\
    \mathcal{Q}^{h} &:= \left\{w^{h}(\mathbf{x}) \in \mathcal{Q} \; \big| \; w^{h}(\mathbf{x}) 
      \in C^{0}(\bar{\Omega}), w^{h}(\mathbf{x}) \big|_{\Omega^e} \in \mathbb{P}^{k}(\Omega^{e}), 
      e = 1, \cdots, Nele \right\}       
  \end{align}
\end{subequations}
where $k$ is a non-negative integer. A corresponding finite element formulation can be written as: 
Find $c^{h}(\mathbf{x}) \in \mathcal{P}^{h}$ such that
\begin{align}
\label{Eqn:Decay_single_field_FE_formulation}
\mathcal{B}(w^{h};c^{h}) = L(w^{h}) \quad \forall w^{h}(\mathbf{x}) \in \mathcal{Q}^{h}
\end{align}

\subsection{A methodology for enforcing the non-negative constraint and maximum principles}
Before we present a methodology for enforcing the non-negative constraint and (discrete) maximum 
principles under the classical Galerkin formulation, we present some definitions and relevant 
results from numerical optimization. We shall use the symbols $\preceq$ and $\succeq$ to denote 
component-wise inequalities for vectors. That is, for given any two (finite dimensional) vectors 
$\boldsymbol{a}$ and $\boldsymbol{b}$ 
\begin{align}
  \boldsymbol{a} \preceq \boldsymbol{b} \quad 
  \mbox{means that } \quad a_i \leq b_i \; \forall i
\end{align}
Similarly one can define the symbol $\succeq$. Let us denote the standard inner-product in 
Euclidean spaces by $<\cdot;\cdot>$. A problem in quadratic programming takes the form
\begin{subequations}
  \label{Eqn:Decay_QP}
  \begin{align}
    \mathop{\mbox{minimize}}_{\boldsymbol{x}}  & \quad f_0(\boldsymbol{x}) := \frac{1}{2} 
    <\boldsymbol{x}; \boldsymbol{Q} \boldsymbol{x}> - <\boldsymbol{x}; \boldsymbol{g}> \\
    \label{Eqn:Decay_QP_standard_inequality}
    \mbox{subject to} & \quad \boldsymbol{A} \boldsymbol{x} \preceq \boldsymbol{b} 
    \quad \mbox{(inequality constraints)} \\
    & \quad \boldsymbol{C} \boldsymbol{x} = \boldsymbol{d} \quad \mbox{(equality constraints)}
  \end{align}
\end{subequations}
The above problem belongs to \emph{convex quadratic programming} if $\boldsymbol{Q}$ 
is positive semidefinite (which makes the objective function $f_0(\boldsymbol{x})$ 
to be convex). As the name suggests, convex quadratic programming is a special case 
of convex optimization. For further details on convex optimization and associated 
numerical algorithms see References \cite{Boyd_convex_optimization,Nocedal_Wright,
Luenberger_Ye_Nonlinear_Programming}. 

\begin{figure}[!h]
  \centering
  \subfigure[$\alpha = 1$, $c_{\mathrm{min}} =  0.8863$]{
    \includegraphics[scale=0.35]{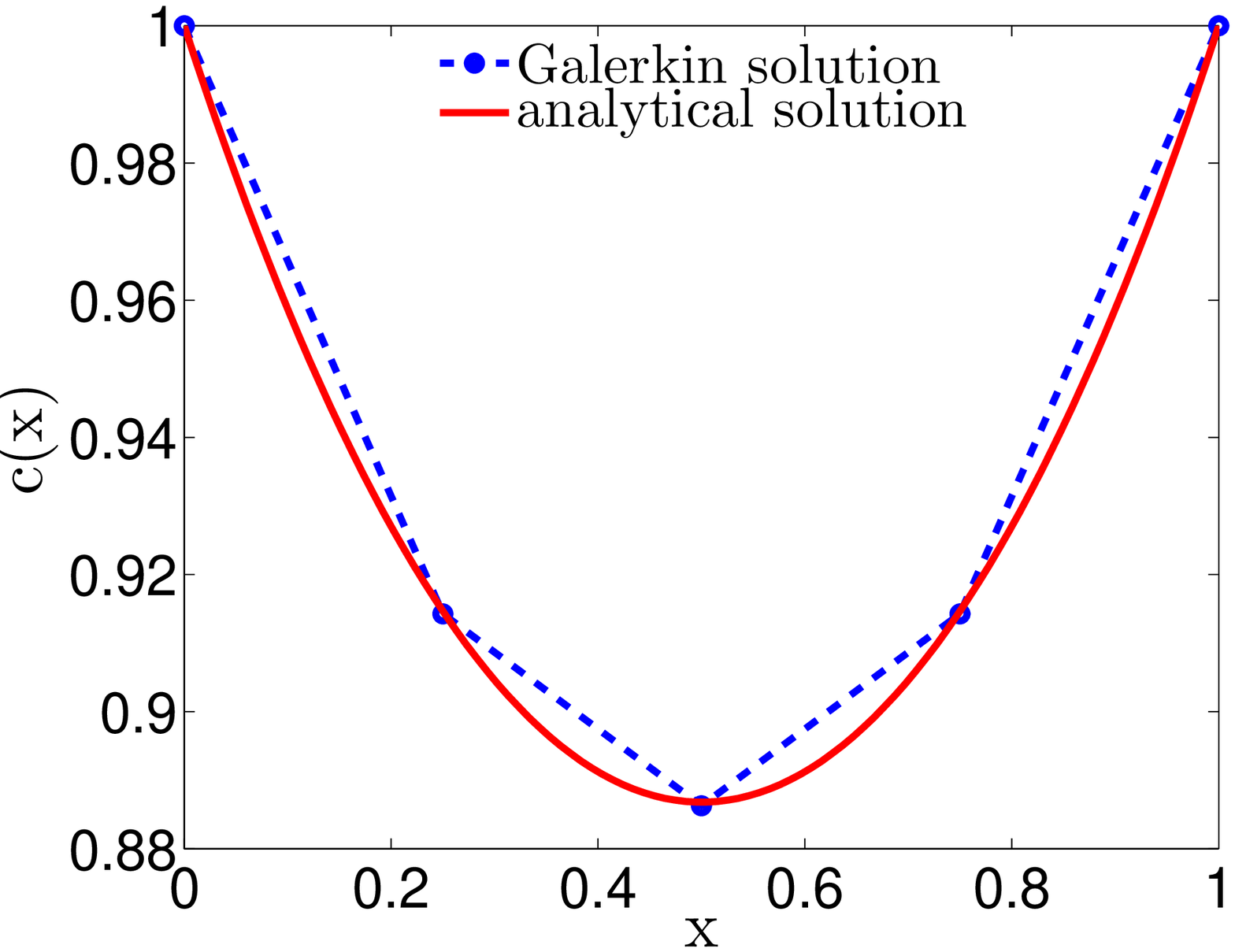}}        
  \subfigure[$\alpha = 100$, $c_{\mathrm{min}} =  -0.0068$]{
    \includegraphics[scale=0.35]{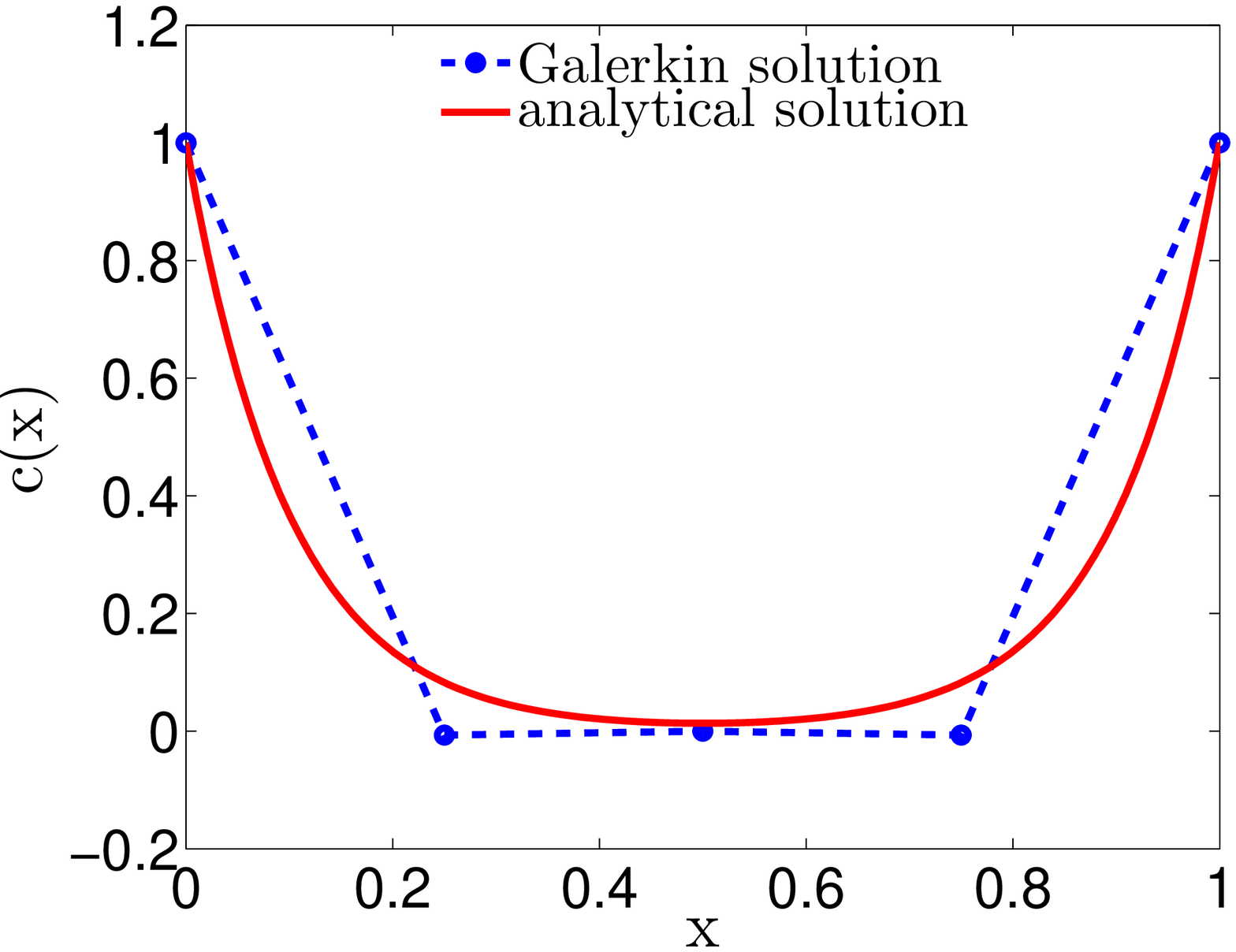}} 
  \subfigure[$\alpha = 500$, $c_{\mathrm{min}} = -0.1977$]{
    \includegraphics[scale=0.35]{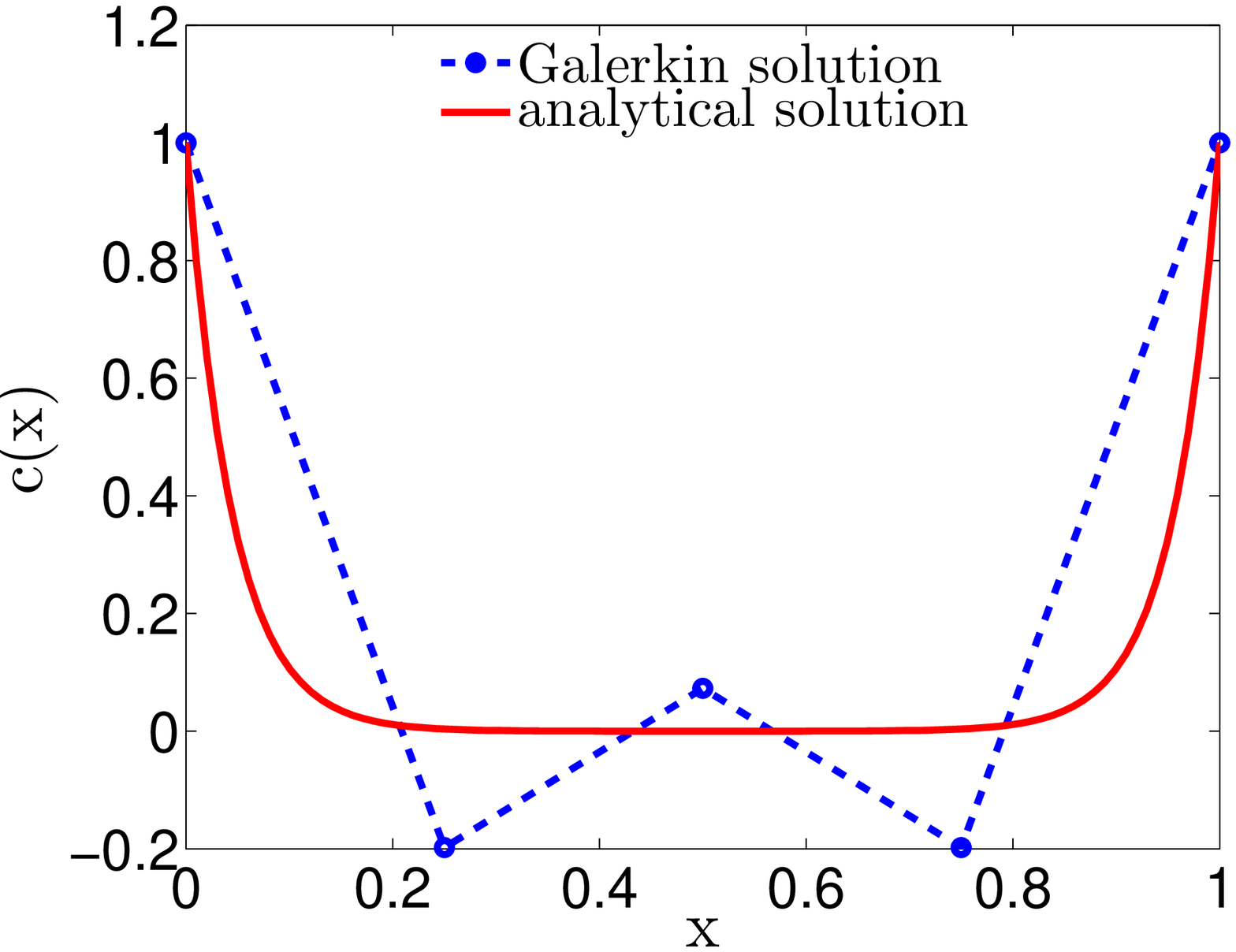}}        
  \subfigure[$\alpha = 1000$, $c_{\mathrm{min}} = -0.2378$]{
    \includegraphics[scale=0.35]{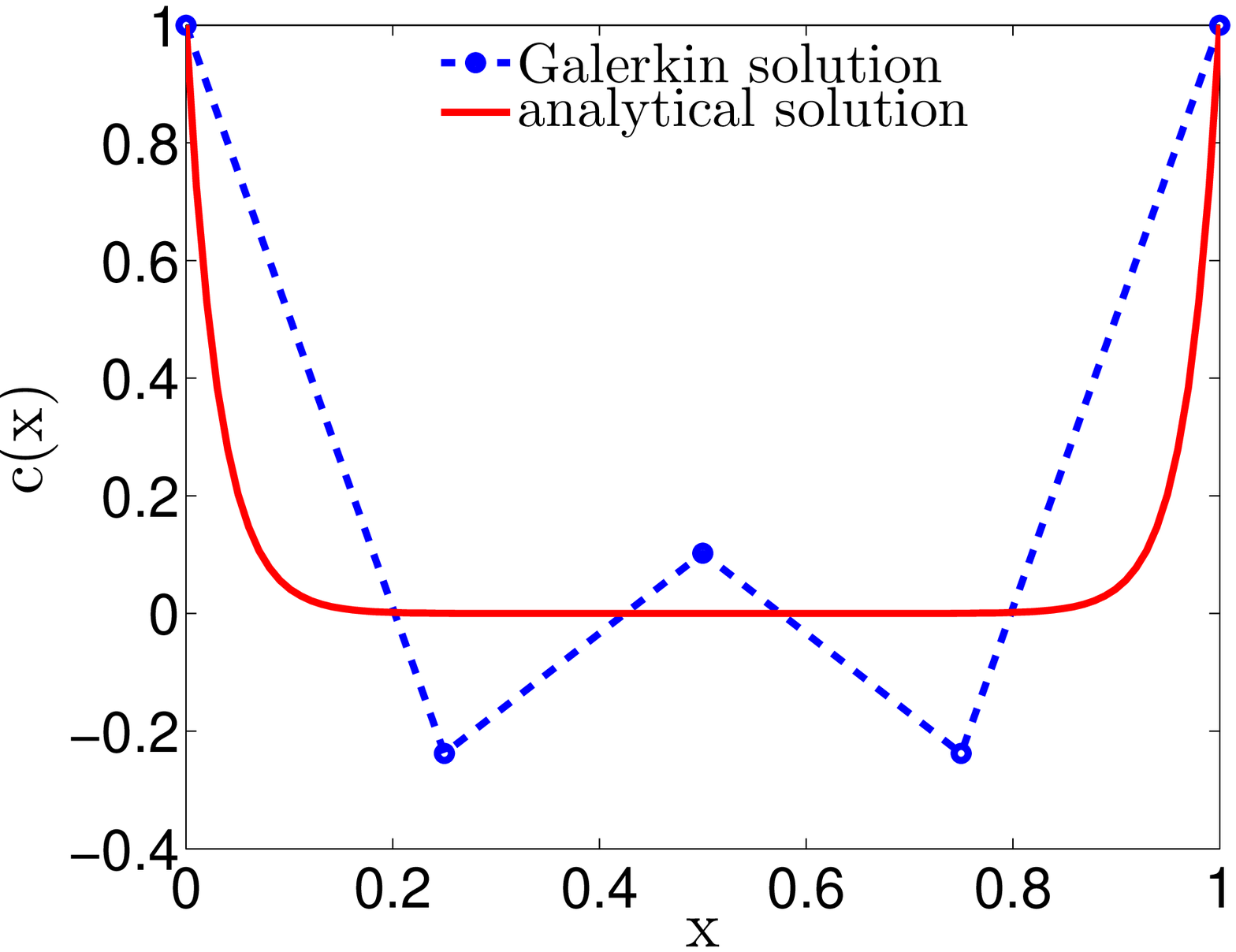}} 
  \caption{One-dimensional problem: Comparison of the numerical solution obtained 
  using the classical Galerkin formulation with the analytical solution for various 
  values of the decay coefficient $\alpha$. Note that the larger the value of $\alpha$, 
  the larger is the violation of the discrete maximum principle by the classical 
  Galerkin formulation.} 
  \label{Fig:Decay_1D_various_alpha}
\end{figure}

\begin{figure}[h]
  \centering
  \subfigure[Three-node triangular mesh with 32 elements]{
    \includegraphics[scale=0.335]{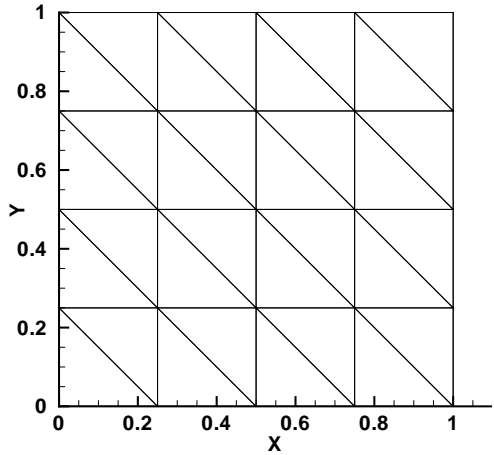}}        
  \subfigure[Classical Galerkin formulation]{
    \includegraphics[scale=0.335]{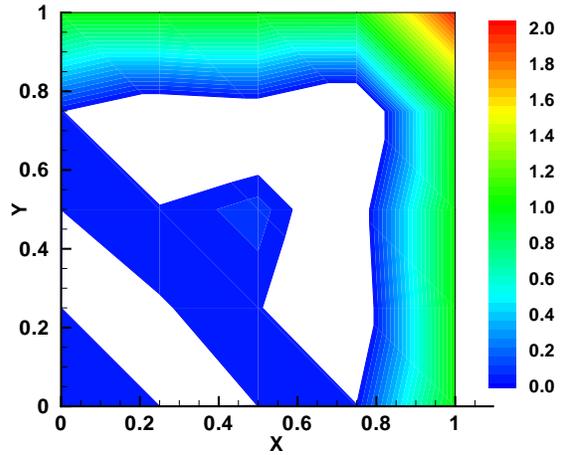}}        
  \subfigure[Proposed formulation]{
    \includegraphics[scale=0.335]{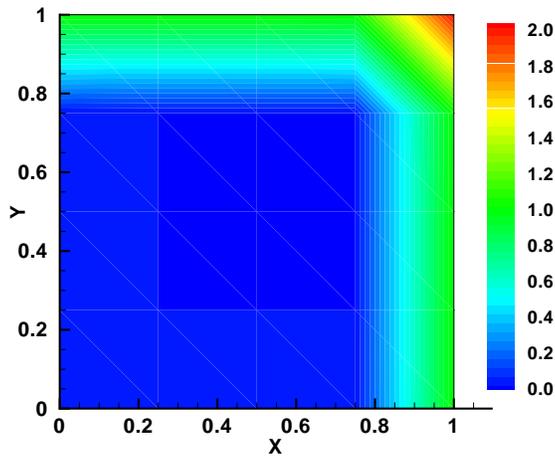}} 
  \caption{Two-dimensional problem with isotropic medium: This figure shows the 
    contours of concentration for $\alpha = 500$ on a \emph{coarse computational mesh} 
    under the Galerkin formulation and the proposed formulation. Regions that have 
    negative concentrations are indicated in white color. The proposed formulation 
    produced physically meaningful non-negative concentrations in the entire computational 
    domain, Under the classical Galerkin formulation, approximately $24\%$ of the total 
    number of nodes have negative concentration. Also, under the classical Galerkin 
    formulation, the  minimum concentration is -0.4049, which occurred inside the 
    domain thereby violating the maximum principle.}\label{Fig:Decay_2D_coarse_mesh}
\end{figure}

We now return to the finite element implementation of the classical Galerkin formulation 
of anisotropic diffusion with decay. After finite element discretization, the discrete 
equations take the form
\begin{align}
  \label{Eqn:Helmholtz_discrete}
  \boldsymbol{K} \boldsymbol{c} = \boldsymbol{f}
\end{align}
where $\boldsymbol{K}$ is a symmetric positive definite matrix, $\boldsymbol{c}$ is the 
vector containing nodal concentrations, and $\boldsymbol{f}$ is the load vector (arising 
from the forcing function). The corresponding minimization problem can be written as 
\begin{align}
  \label{Eqn:Helmholtz_minimization}
  \mathop{\mbox{minimize}}_{\boldsymbol{c} \in \mathbb{R}^{ndofs}} \quad  \frac{1}{2} 
  <\boldsymbol{c}; \boldsymbol{K}  \boldsymbol{c}> - <\boldsymbol{c}; \boldsymbol{f}>
\end{align}
where ``$ndofs$'' denotes the number of degrees of freedom in the finite element mesh 
(which is equal to the total number of nodes minus the number of nodes at which a Dirichlet 
boundary condition is enforced). As shown in Figures \ref{Fig:Decay_1D_various_alpha} 
and \ref{Fig:Decay_2D_coarse_mesh}, the finite element solution based on equation 
\eqref{Eqn:Helmholtz_discrete} produces unphysical negative concentrations even for 
simple problems. A formulation corresponding to \eqref{Eqn:Helmholtz_minimization} 
that satisfies the maximum principle (given by Theorem \ref{Theorem:Decay_maximum_principle} 
and Remark \ref{Remark:Decay_max_min})  can be written as 
\begin{subequations}
  \label{Eqn:Decay_DMP}
  \begin{align}
    &\mathop{\mbox{minimize}}_{\boldsymbol{c} \in \mathbb{R}^{ndofs}} \quad  
    \frac{1}{2} <\boldsymbol{c}; \boldsymbol{K}  \boldsymbol{c}>  - 
    <\boldsymbol{c}; \boldsymbol{f}> \\
    &\mbox{subject to} \quad c_{\mathrm{min}} \boldsymbol{1} \preceq \boldsymbol{c} 
    \preceq c_{\mathrm{max}} \boldsymbol{1} 
  \end{align}
\end{subequations}
where $\boldsymbol{1}$ is a vector of size $ndofs$ containing ones, and 
$c_{\mathrm{min}}$ and $c_{\mathrm{max}}$ are, respectively, given by
\begin{subequations}
  \begin{align}
    c_{\mathrm{min}} &:= \min_{\mathbf{x} \in \partial \Omega} c^{-}(\mathbf{x}) \; 
    \mathrm{where} \; c^{-}(\mathbf{x}) = \min\{c(\mathbf{x}),0\} \\
    c_{\mathrm{max}} &:= \max_{\mathbf{x} \in \partial \Omega} c^{+}(\mathbf{x}) \; 
    \mathrm{where} \; c^{+}(\mathbf{x}) = \max\{c(\mathbf{x}),0\}
  \end{align}
\end{subequations}
A corresponding formulation that satisfies the non-negative constraint can be obtained 
by setting $c_{\mathrm{min}} = 0$, and omitting the upper bound (which is equivalent to 
the condition $c_{\mathrm{max}} = +\infty$); and can be written as follows:
\begin{subequations}
  \label{Eqn:Decay_non-negative}
  \begin{align}
    &\mathop{\mbox{minimize}}_{\boldsymbol{c} \in \mathbb{R}^{ndofs}} \quad  
    \frac{1}{2} <\boldsymbol{c}; \boldsymbol{K}  \boldsymbol{c}>  - 
    <\boldsymbol{c}; \boldsymbol{f}> \\
    \label{Eqn:Decay_non-negative_constraint}
    &\mbox{subject to} \quad \boldsymbol{0} \preceq \boldsymbol{c}  
  \end{align}
\end{subequations}
where $\boldsymbol{0}$ is a vector of size $ndofs$ containing zeros. 
\begin{remark}
  The constraint $c_{\mathrm{min}} \boldsymbol{1} \preceq \boldsymbol{c} \preceq 
  c_{\mathrm{max}} \boldsymbol{1}$ can be rewritten in the standard form given by 
  equation \eqref{Eqn:Decay_QP_standard_inequality} as follows:
  \begin{subequations}
    \begin{align*}
      \boldsymbol{c} &\preceq c_{\mathrm{max}} \boldsymbol{1} \\
      - \boldsymbol{c} &\preceq -c_{\mathrm{min}} \boldsymbol{1} 
    \end{align*}
  \end{subequations}
  The constraint \eqref{Eqn:Decay_non-negative_constraint} can be put in the 
  standard form by rewriting it as: $-\boldsymbol{c} \preceq \boldsymbol{0}$.
\end{remark}

Comparing with equation \eqref{Eqn:Decay_QP}, it is evident that the above problems  
\eqref{Eqn:Decay_DMP} and \eqref{Eqn:Decay_non-negative} belong to convex quadratic 
programming. The first-order optimality conditions (which are given by Karush-Kuhn-Tucker 
conditions) corresponding to equation \eqref{Eqn:Decay_non-negative} take the following 
form:
\begin{subequations}
  \label{Eqn:Decay_first_order_optimality}
  \begin{align}
    &\boldsymbol{K} \boldsymbol{c} = \boldsymbol{f} + \boldsymbol{\lambda} \\
    \label{Eqn:Decay_QP_inequality}
    &\boldsymbol{c} \succeq \boldsymbol{0} \\
    \label{Eqn:Decay_QP_lambda_inequality}
    &\boldsymbol{\lambda} \succeq \boldsymbol{0} \\
    \label{Eqn:Decay_QP_complementary_conditions}
    &\lambda_{i} c_{i} = 0 \quad (i = 1, \cdots, ndofs)
  \end{align}
\end{subequations}
where $\boldsymbol{\lambda}$ is a vector of Lagrange multipliers enforcing 
the constraint \eqref{Eqn:Decay_QP_inequality}. Similarly, one can write 
first-order optimality conditions for the optimization problem given by 
equation \eqref{Eqn:Decay_DMP}.

\begin{remark}
  The above set of equations \eqref{Eqn:Decay_first_order_optimality} is not 
  linear because of the inequality constraints \eqref{Eqn:Decay_QP_inequality} 
  and \eqref{Eqn:Decay_QP_lambda_inequality} and complementary conditions 
  \eqref{Eqn:Decay_QP_complementary_conditions}.  
\end{remark}

\section{REPRESENTATIVE NUMERICAL RESULTS}
\label{Sec:Decay_NR}
In this section, we illustrate the performance of the proposed non-negative formulation 
for the anisotropic diffusion with decay using representative one- and two-dimensional 
problems. 
In all our numerical experiments we have employed the standard active-set strategy 
\cite{Luenberger_Ye_Nonlinear_Programming} to solve resulting convex quadratic 
programming problems. In all our numerical simulations we have taken the violated 
nodes under the classical Galerkin formulation as the initial active-set. This choice 
is motivated by the numerical studies reported by Nakshatrala and Valocchi 
\cite{Nakshatrala_Valocchi_JCP_2009_v228_p6726} in which it has been shown that the 
initial active-set based on the violated nodes from the underlying formulation, in most cases, 
takes fewer active-set strategy iterations (than, say, empty set as the initial guess). 

\subsection{One-dimensional problem}
\begin{figure}[htp]
  \centering
  \includegraphics[scale=0.45]{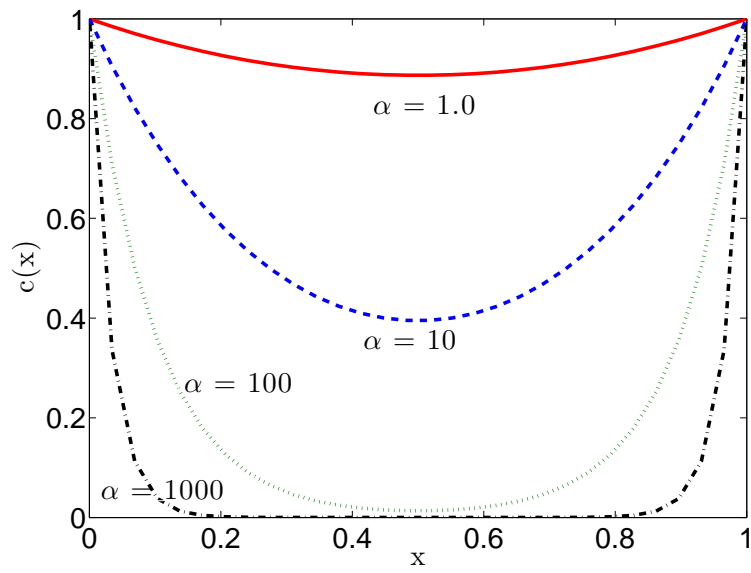}
  \caption{One-dimensional problem: Analytical solution for various 
  values of alpha.} \label{Fig:TransientDMP_1D_exactsol}
\end{figure}

Consider the following simple one-dimensional problem with homogeneous forcing function: 
\begin{subequations}
  \label{Eqn:Transient_1D_Helmholtz}
  \begin{align}
    &\alpha c(\mathrm{x}) - \frac{d^2 c}{d \mathrm{x}^2}  =  0 
    \quad \mathrm{in} \; \Omega := (0,1) \\
    & c(\mathrm{x}=0) = c(\mathrm{x}=1) = 1
  \end{align}
\end{subequations}
with $\alpha \geq 0$. The analytical solution to the above problem is given by 
\begin{align}
  c(\mathrm{x}) = \dfrac{1-\exp(-\sqrt{\alpha} )}{\exp(\sqrt{\alpha} )-\exp(-\sqrt{\alpha})}
  \exp(\sqrt{\alpha} \mathrm{x}) + \dfrac{\exp(\sqrt{\alpha} )-1}{\exp(\sqrt{\alpha} ) - 
    \exp(-\sqrt{\alpha} )}\exp(-\sqrt{\alpha} \mathrm{x}) 
\end{align}
In Figure \ref{Fig:TransientDMP_1D_exactsol}, the analytical solution is plotted for various 
values of $\alpha$. As one can see from the figure, sharp boundary layers exist for higher 
values of $\alpha$. For obtaining numerical results, the computational domain is divided 
into four equal-sized elements. The numerical results obtained using the classical Galerkin 
formulation are shown in Figure \ref{Fig:Decay_1D_various_alpha}, and the classical Galerkin 
formulation clearly violates the discrete maximum principle for higher value of $\alpha$. 
From this figure it is also observed that the larger the value of $\alpha$ the larger is 
the violation of the discrete maximum principle. However, for one-dimensional problems, 
the violation of discrete maximum principle decreases with mesh refinement (which is not 
true, in general, in higher spatial dimensions especially when anisotropy dominates).  

\begin{figure}[!h]
  \centering
  \subfigure{
    \includegraphics[scale=0.35]{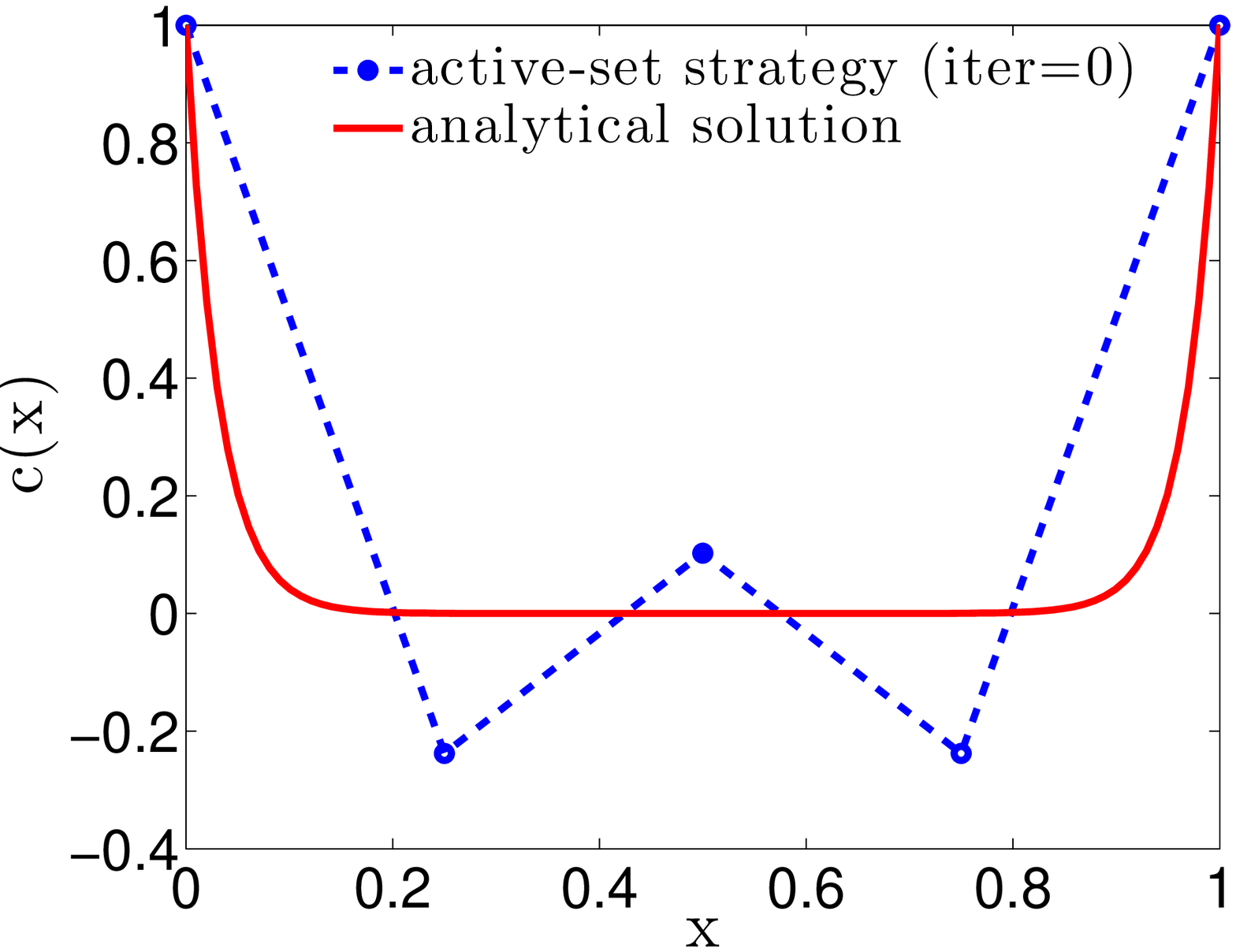}}        
  \subfigure{
    \includegraphics[scale=0.35]{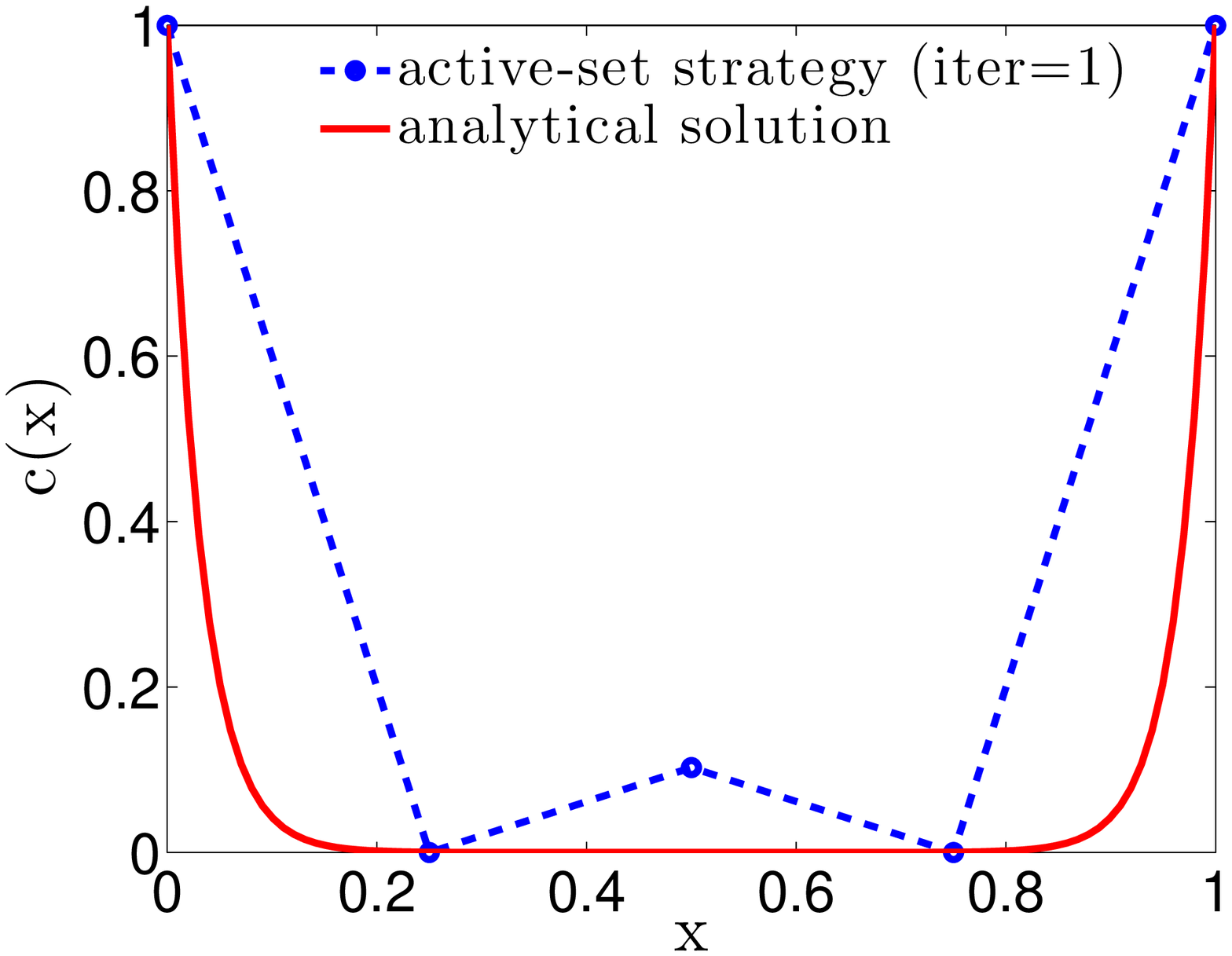}} 
  \subfigure{
    \includegraphics[scale=0.35]{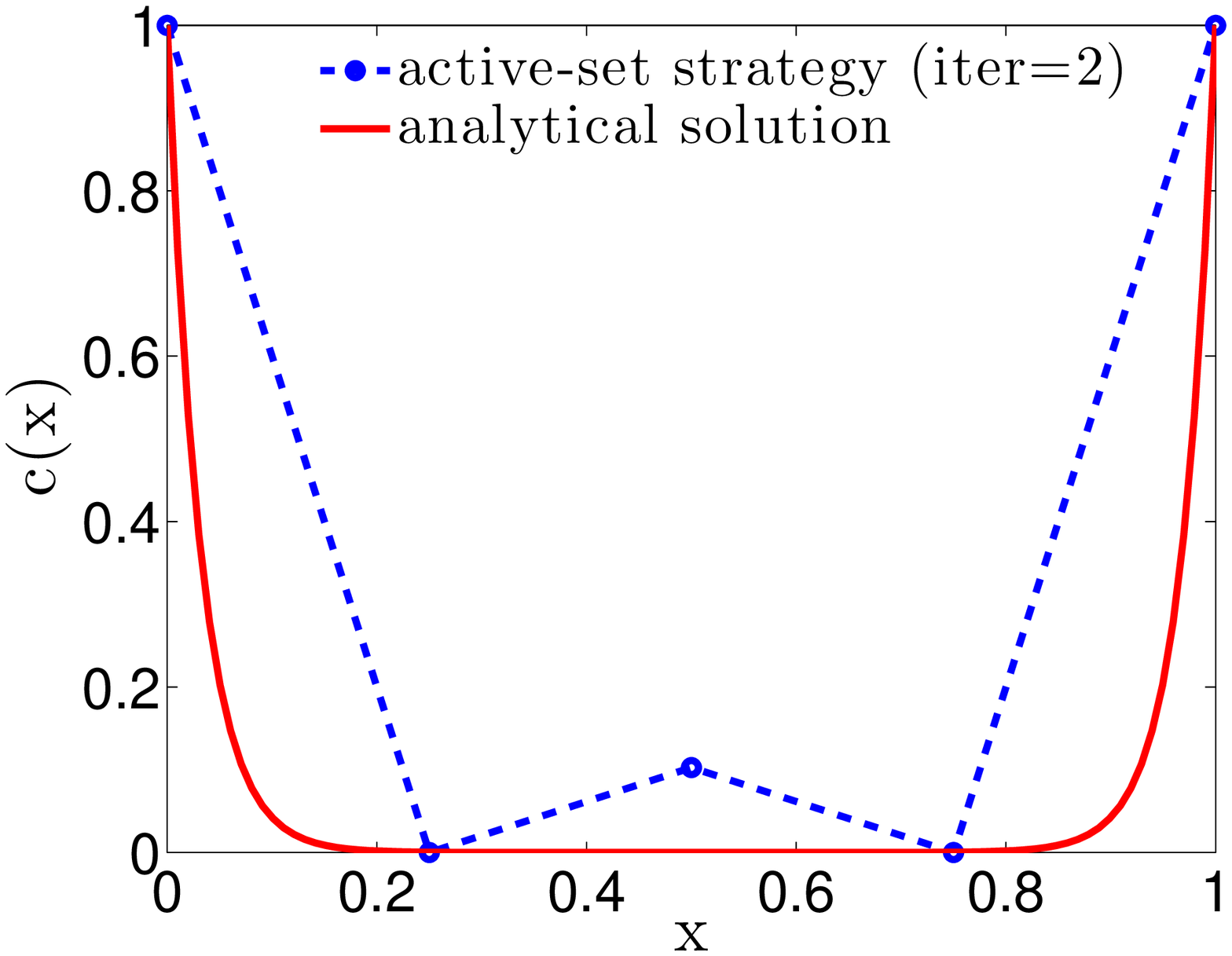}}        
  \subfigure{
    \label{Fig:Helmholtz_1D_converged}
    \includegraphics[scale=0.35]{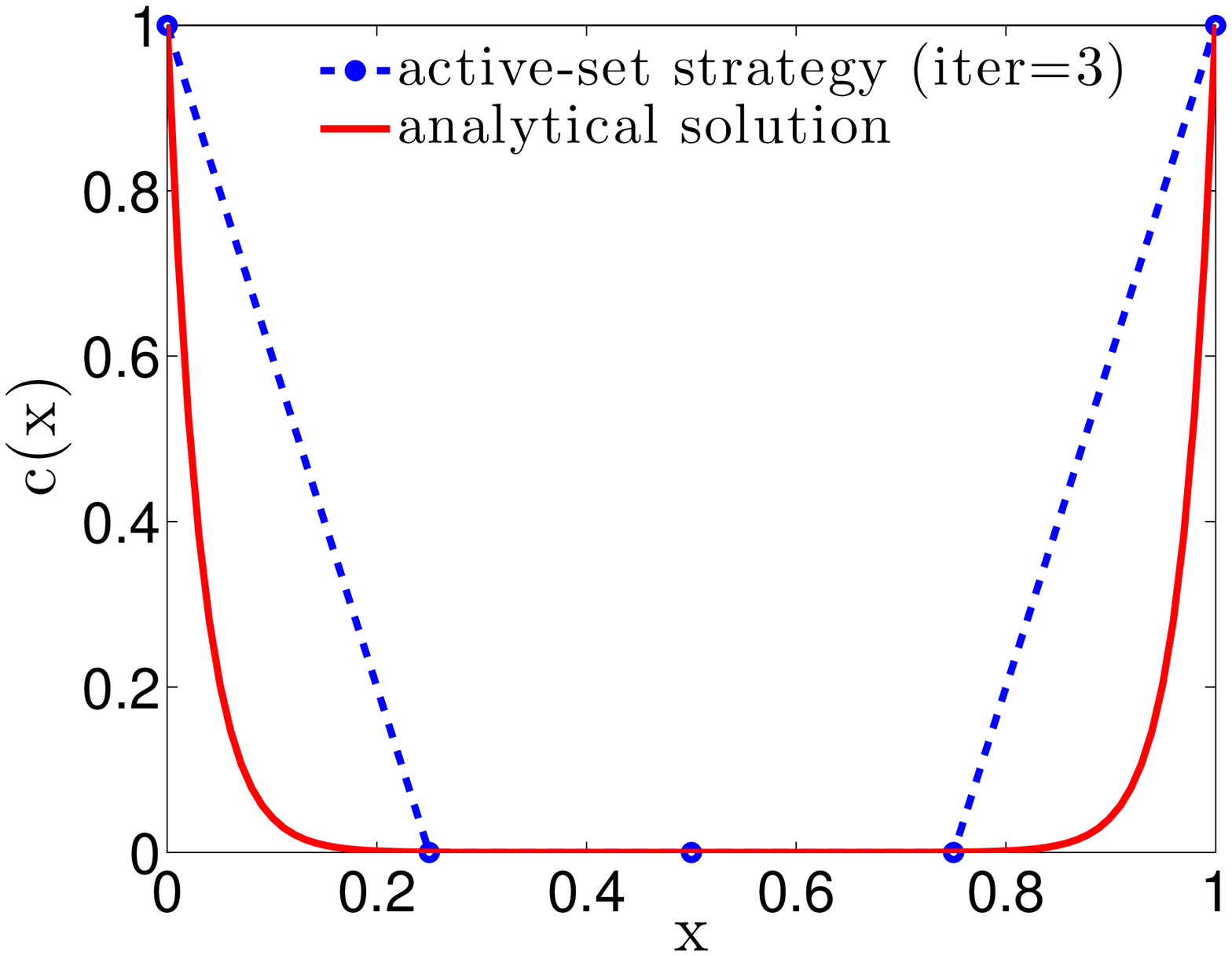}} 
  \caption{One-dimensional problem: This figure shows the variation of the numerical 
    solution under the proposed formulation at various active-set strategy iterations 
    for $\alpha = 1000$. The active-set strategy converged after three iterations. Note 
    that, for this problem, the converged numerical solution from the proposed formulation 
    matches exactly at nodes with the analytical solution.} 
    \label{fig:Helmholtz_1D_NonNeg_iters}
\end{figure}

\begin{figure}
  \includegraphics[scale=0.35]{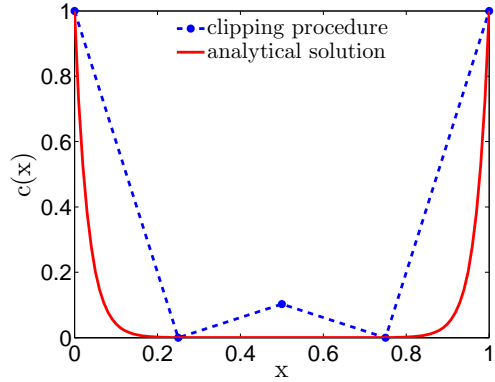}
  \caption{One-dimensional problem: This figure compares the analytical solution 
    with the numerical solution obtained using the ``clipping procedure,'' which 
    basically chops offs all the negative nodal concentrations obtained from the 
    Galerkin formulation by setting them to zero. We have taken $\alpha = 1000$ 
    in this numerical simulation. The corresponding numerical solution obtained 
    using the proposed methodolody is shown in Figure \ref{Fig:Helmholtz_1D_converged}, 
    and the proposed methodolody performs better than the clipping procedure.} 
  \label{Fig:Helmholtz_1D_clipping}
\end{figure}

We have solved the above one-dimensional problem using the proposed formulation, and have 
employed the active-set strategy to solve the resulting convex quadratic programming problem. 
We have taken $\alpha = 1000$, and have employed the same computational mesh as discussed 
above. Figure \ref{fig:Helmholtz_1D_NonNeg_iters} illustrates how the active-set strategy 
performed at various iterations, and the active-set strategy converged in three iterations. 
In Figure \ref{Fig:Helmholtz_1D_clipping}, we have shown the performance of the ``\emph{clipping 
procedure}" in which all the negative nodal concentrations from the Galerkin formulation 
are chopped off by setting them to zero. As discussed in the caption of the figure, the 
proposed methodology performs better than the clipping procedure. Also, the clipping 
procedure does not have a variational basis, and is considered as a ``variational crime."

\begin{figure}[h]
  \centering
      \subfigure[$\alpha = 1$]{
      \includegraphics[scale=0.45]{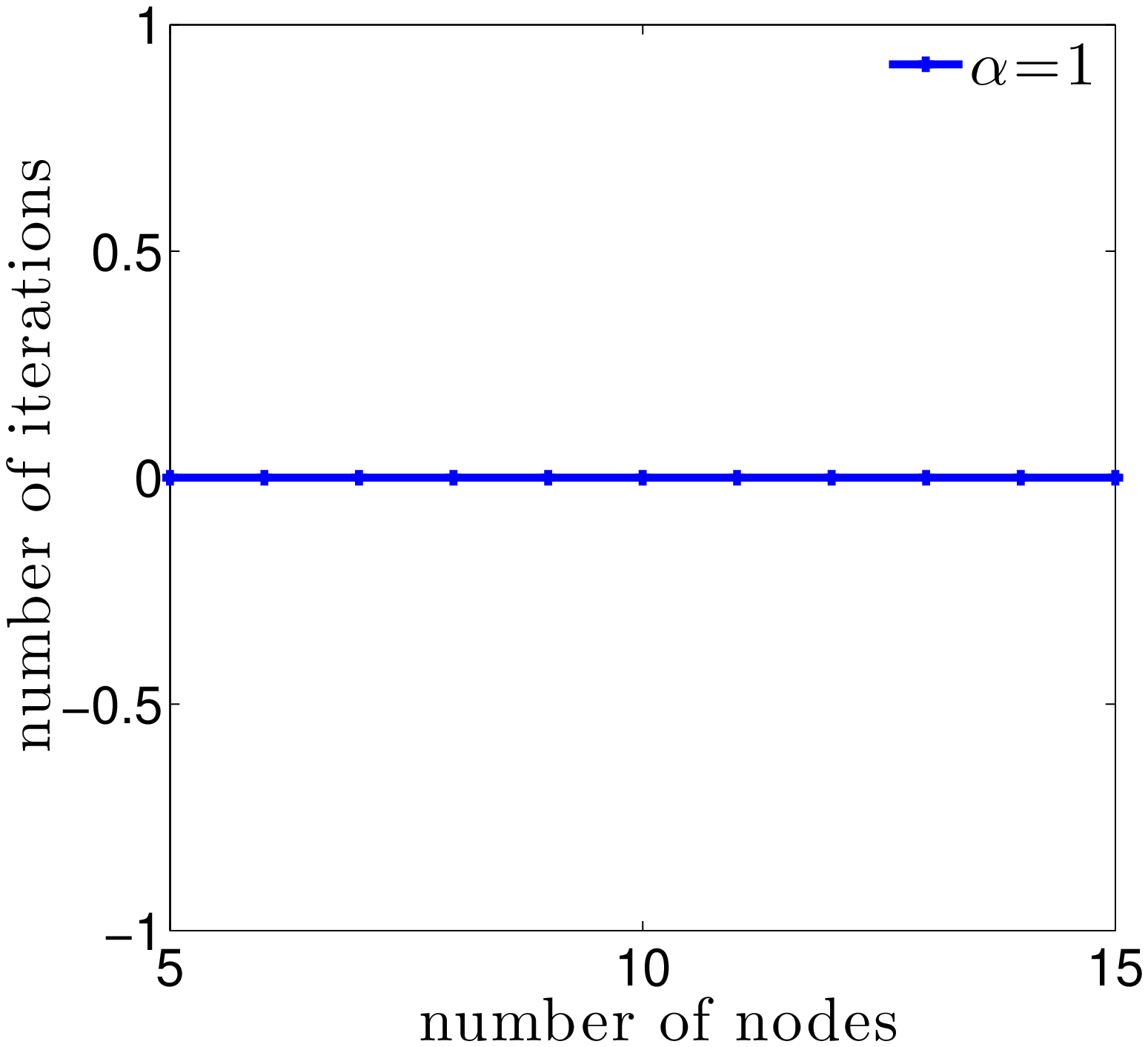}}        
      \subfigure[$\alpha = 100$]{
      \includegraphics[scale=0.45]{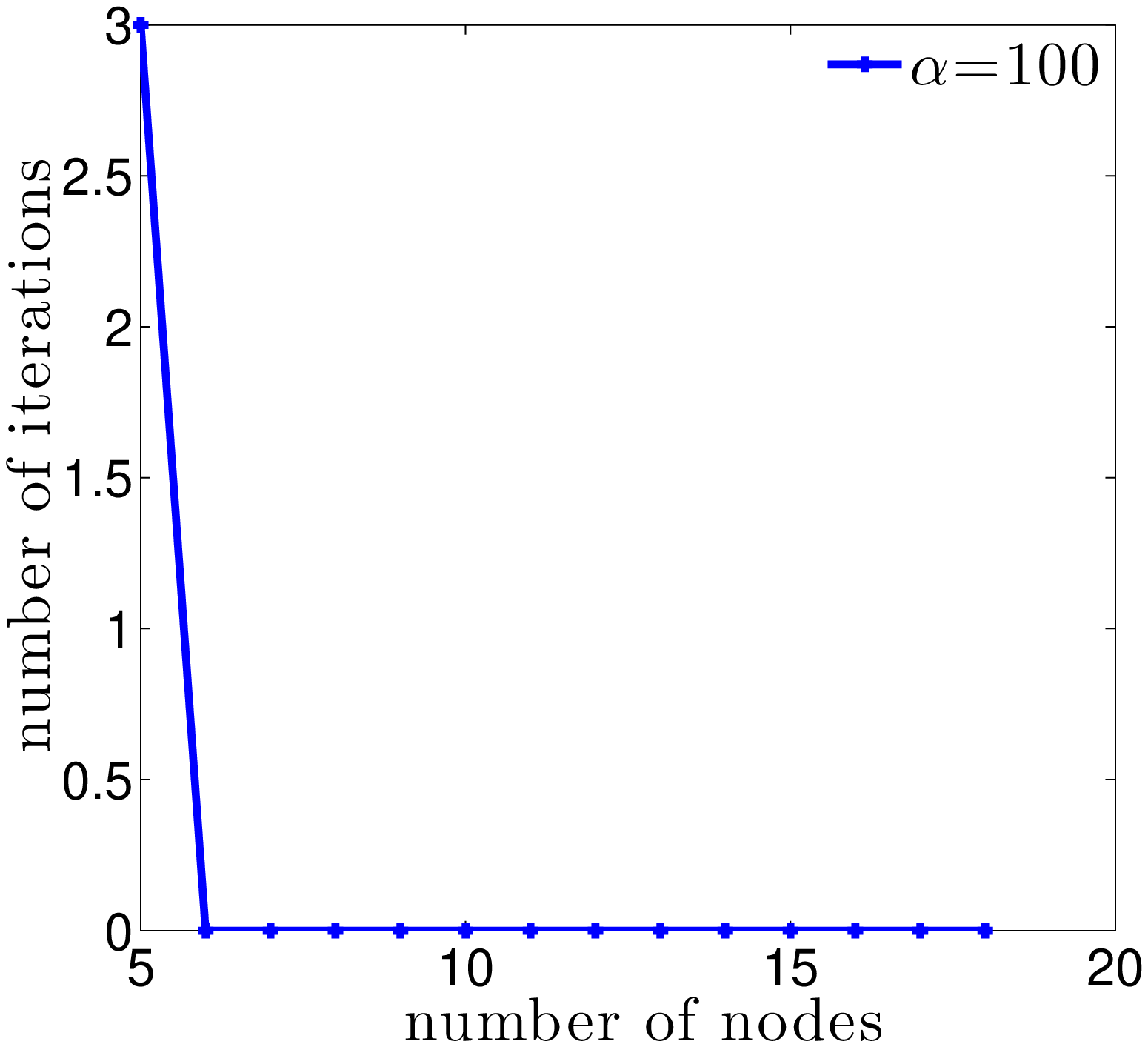}} 
      \subfigure[$\alpha = 500$]{
      \includegraphics[scale=0.45]{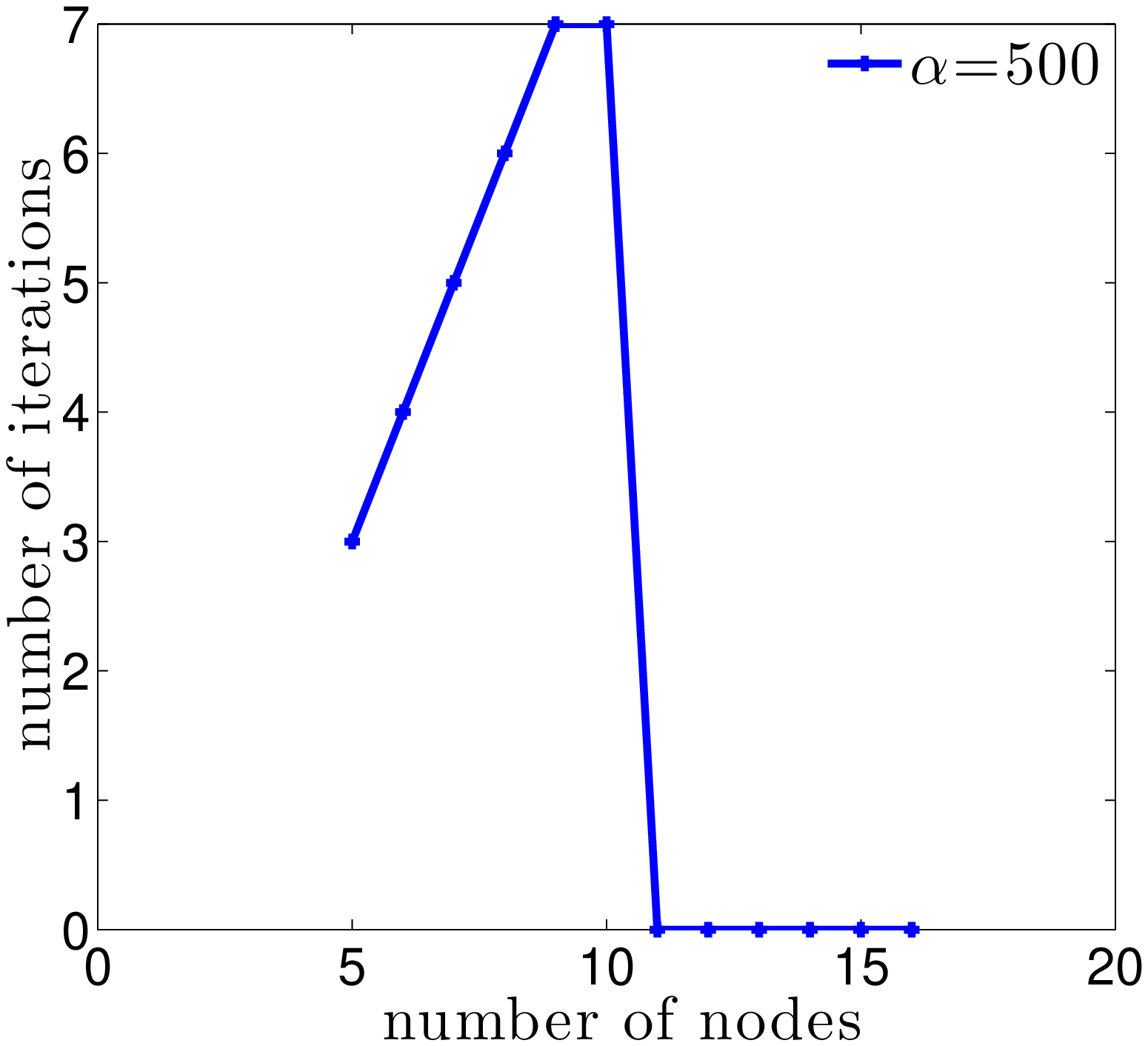}} 
      \subfigure[$\alpha = 1000$]{
      \includegraphics[scale=0.45]{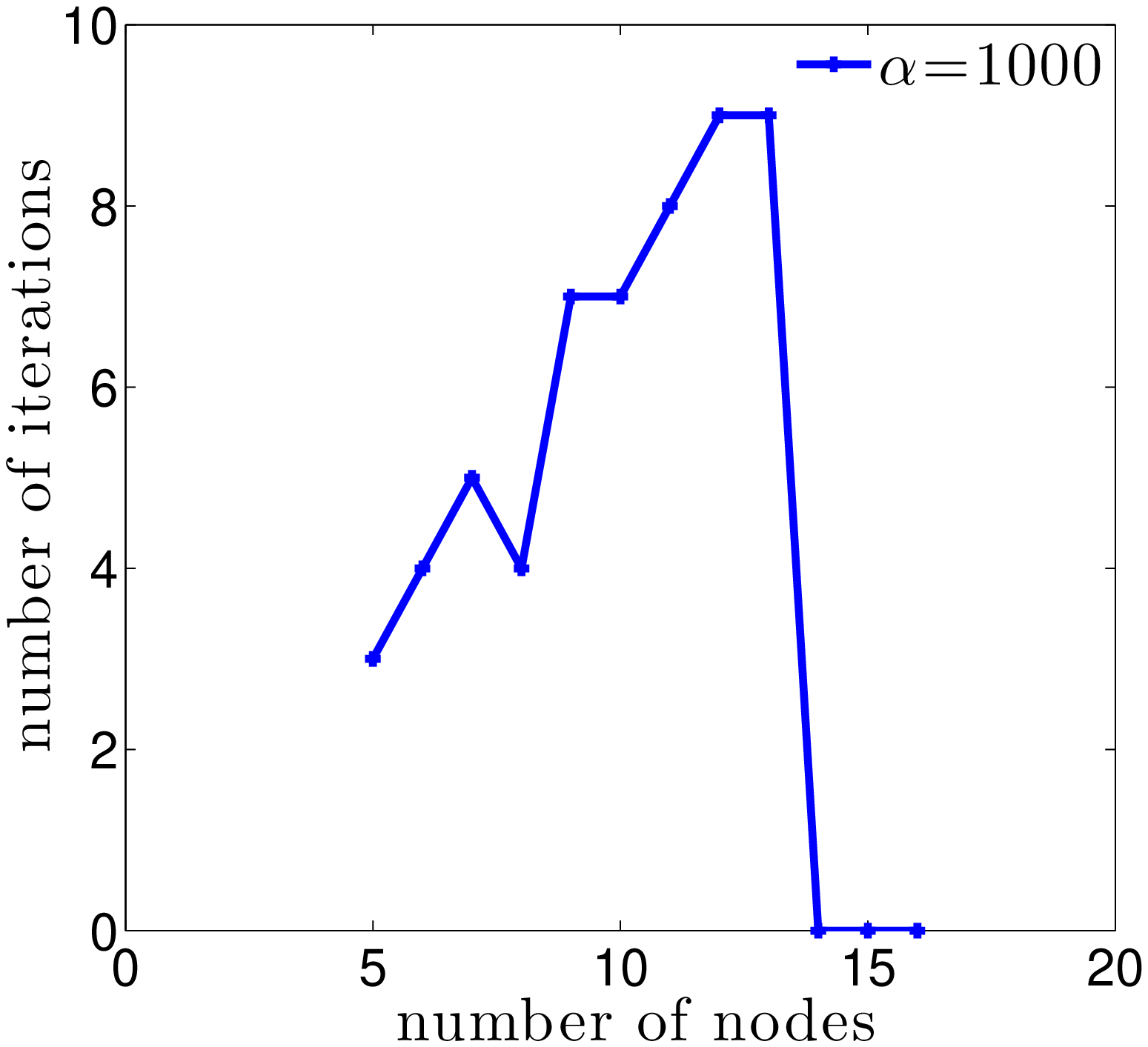}} 
    \caption{One-dimensional problem: These figures present the number of iterations required 
      for the proposed formulation using the active-set strategy for various values of $\alpha$ 
      with respect to the number of nodes. Note that the number of iterations required for the 
      optimization to terminate increases as the value of $\alpha$ increases. After sufficient 
      mesh refinement, there will be no violation of the discrete maximum principle, and there 
      is no need to solve the constrained optimization problem.}
    \label{Fig:Decay_1D_active_set_iterations}
\end{figure}

\begin{figure}[h]
  \centering
      \subfigure[$\alpha = 1$]{
      \includegraphics[scale=0.35]{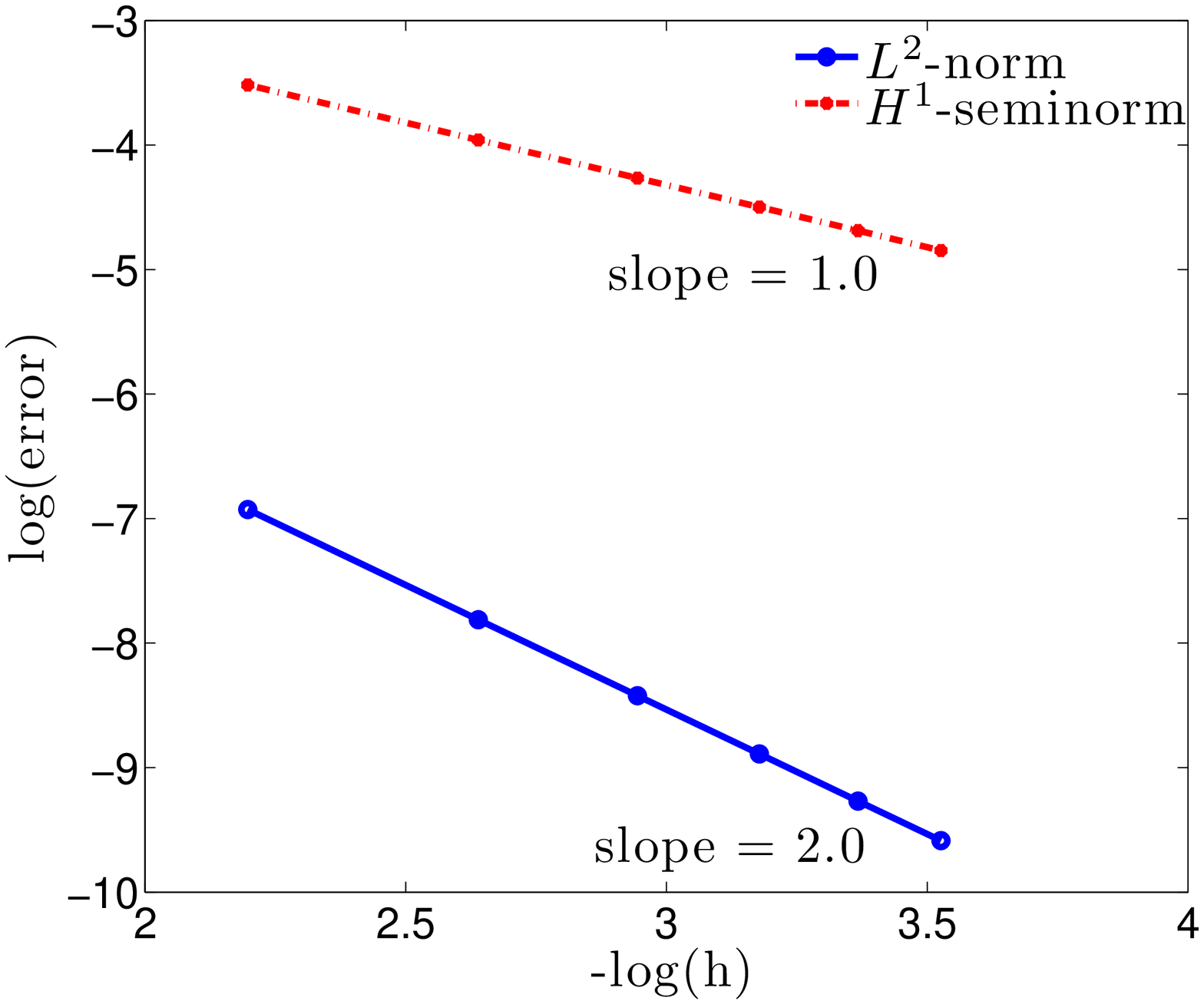}}        
      \subfigure[$\alpha = 100$]{
      \includegraphics[scale=0.35]{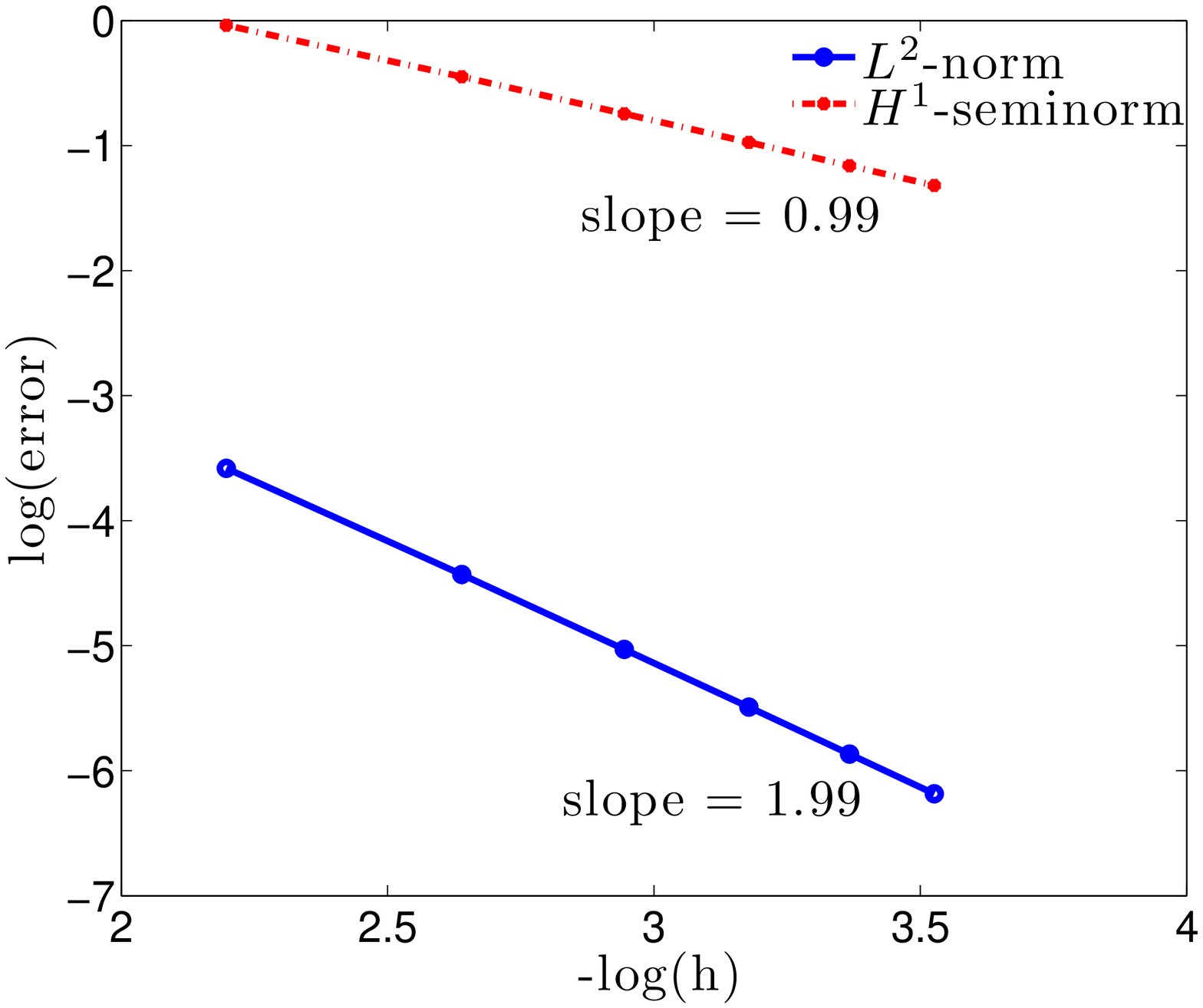}} 
      \subfigure[$\alpha = 500$]{
      \includegraphics[scale=0.35]{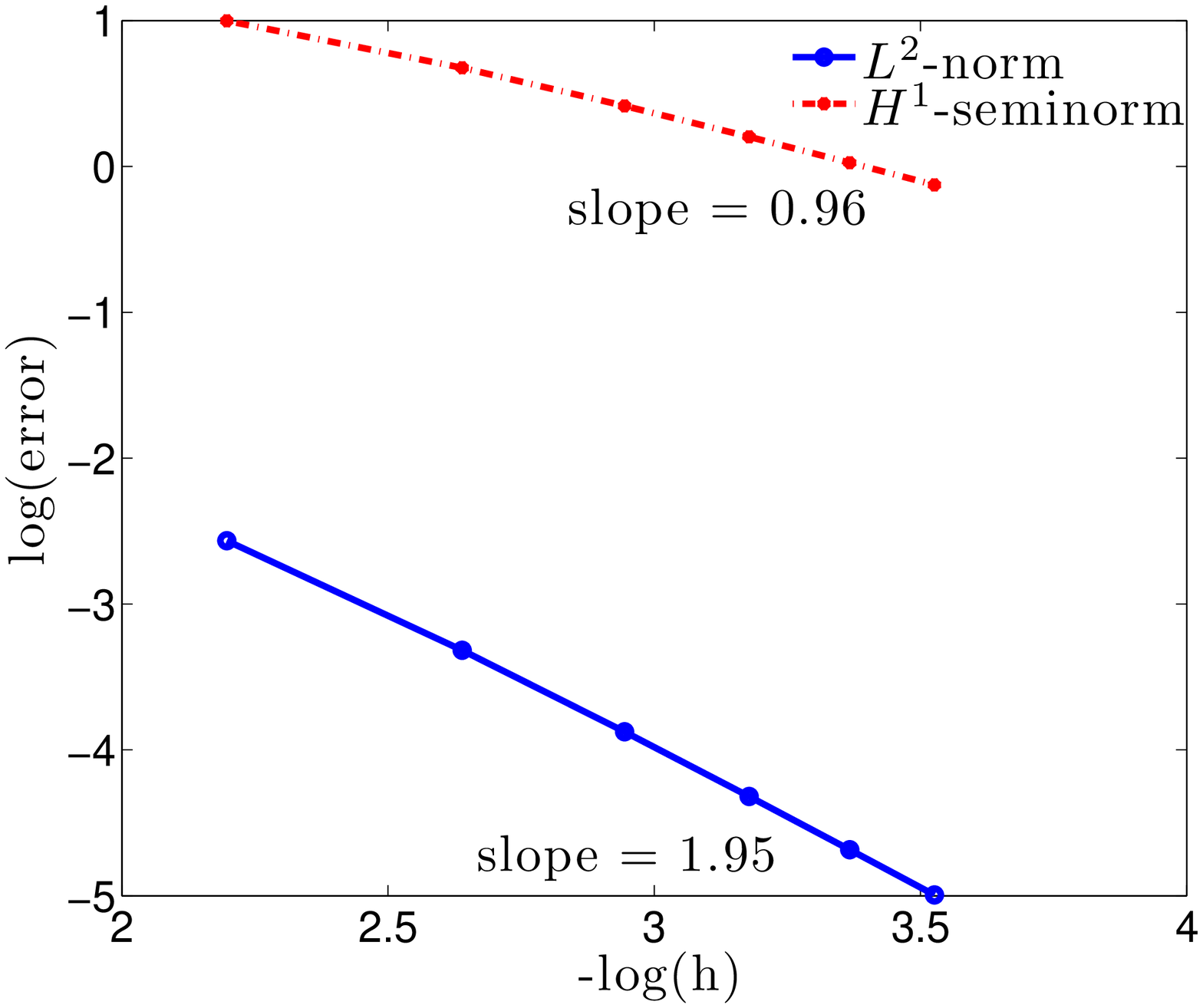}} 
      \subfigure[$\alpha = 1000$]{
      \includegraphics[scale=0.35]{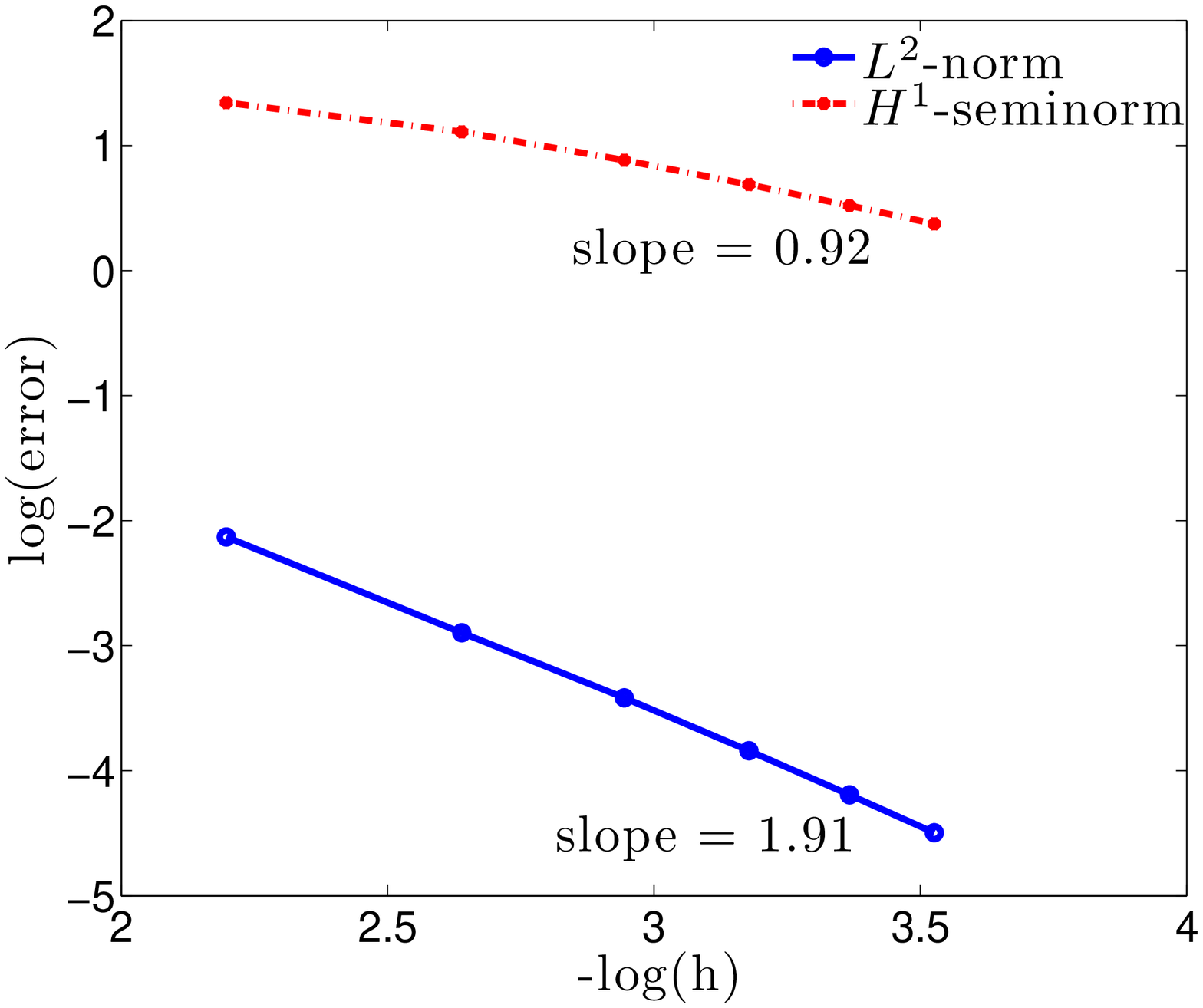}} 
   \caption{One-dimensional problem: This figure presents numerical 
  convergence of the proposed formulation with mesh refinement for 
  various values of decay coefficient. From the figure it is evident 
  that the rates of convergence with respect to mesh refinement in 
  $L^{2}$-norm and $H^{1}$-seminorm are about the same as for the original 
  linear finite element method.} \label{Fig:Decay_error_convergence_1D}
  \end{figure}

In Figure \ref{Fig:Decay_1D_active_set_iterations}, we plot the number of iterations taken 
by the active-set strategy with respect to number of (finite element) nodes. For one-dimensional 
problems, the violation of the discrete maximum principle decreases with mesh refinement, and 
eventually there will no violation of the discrete maximum principle. This can be seen in Figure 
\ref{Fig:Decay_1D_active_set_iterations} as the number of active-set strategy iterations is zero 
for (sufficiently) finer computational meshes for various values of $\alpha$.  In Figure 
\ref{Fig:Decay_error_convergence_1D}, we have shown the convergence of the proposed 
formulation with respect to mesh refinement for various values of $\alpha$, and the 
proposed formulation performed well. 

We now derive sufficient conditions for uniform meshes for one-dimensional problems under 
the classical Galerkin formulation to satisfy the maximum principle. We shall use the following 
results from Matrix Analysis \cite{Quarteroni_Numerical_Mathematics,Varga_Matrix_Iterative_Analysis}: 
Given $\boldsymbol{A} \boldsymbol{x} = \boldsymbol{b}$ with $\boldsymbol{b} \succeq \boldsymbol{0}$, 
sufficient conditions to ensure that $\boldsymbol{x} \succeq \boldsymbol{0}$ are
\begin{enumerate}[(a)]
  \item positive diagonal entries: $A_{ii} > 0$, 
  \item non-positive off-diagonal entries: $A_{ij} \leq 0 \; \forall i \neq j$, and 
  \item strict diagonal dominance by rows: $|A_{ii}| > \sum_{j \neq i} |A_{ij}| \; \forall i$.
\end{enumerate}
\begin{remark}
Note that the aforementioned sufficient conditions to ensure $\boldsymbol{x} \succeq \boldsymbol{0}$ 
are quite stringent, and weaker (sufficient) conditions can be devised. For example, weaker sufficient 
conditions that ensures $\boldsymbol{x} \succeq \boldsymbol{0}$ are: $\boldsymbol{A}$ is invertible, 
and all the entries in $\boldsymbol{A}^{-1}$ are non-negative. 

Another sufficient condition that can be used requires that the matrix $\boldsymbol{A}$ to 
be an M-matrix, which is widely used in the numerical studies on flux and slope limiters 
\cite{Kuzmin_Turek_JCP_2002_v175_p525} and iterative linear solvers \cite{Saad}. An M-matrix 
is a non-singular matrix whose off-diagonal elements are non-positive and all entries of the 
inverse $\boldsymbol{A}^{-1}$ are non-negative. Note that there are many equivalent definitions 
for an M-matrix \cite{Saad,Varga_Matrix_Iterative_Analysis}, and the definition we just outlined 
is quite suitable for our discussion. 

Note that an M-matrix, by definition, has all the entries in its inverse to be non-negative. 
It can be shown that a matrix with positive diagonal entries, non-positive off-diagonal entries, 
and strict diagonal dominance by rows is an M-matrix \cite{Quarteroni_Numerical_Mathematics}. 
We have employed the sufficient conditions outlined just above this remark as they are easy to 
verify, and also suffice our purpose. 
\end{remark}

\begin{figure}[h]
  \centering
  \includegraphics[scale=0.5]{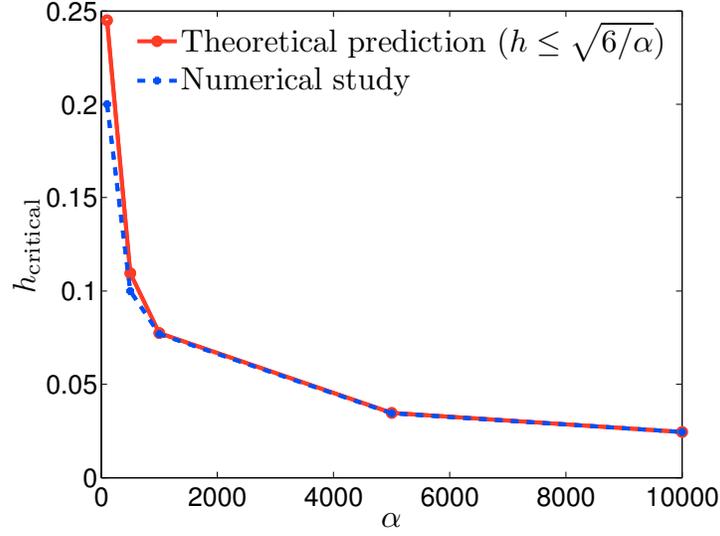}
  \caption{One-dimensional problem: In this figure we compare the sufficient condition 
  derived for uniform one-dimensional problems to satisfy maximum principles with 
  numerical results, and the theoretical prediction is found out to be excellent.} 
  \label{Fig:Decay_1D_h_critical}
\end{figure}

We shall apply the above mathematical result to equation \eqref{Eqn:Helmholtz_discrete}, 
which arises from the finite element discretization of diffusion with decay using the classical 
Galerkin formulation. The computational domain is discretized using equal-sized two-node 
linear finite elements, and let $h$ denotes the size of an element. Since the forcing function 
is assumed to be a source (that is, $f(\mathrm{x}) \geq 0$), and the prescribed Dirichlet 
boundary conditions are non-negative; and we have $\boldsymbol{f} \succeq \boldsymbol{0}$. 
To get sufficient conditions for non-negative nodal concentration, we need to assess the entries 
of the ``stiffness matrix'' $\boldsymbol{K}$.  The entries of the stiffness matrix for an intermediate 
node (say $i$-th node) after the finite element discretization using two-node linear element under 
the classical Galerkin formulation take the following form:
\begin{align}
  \frac{\alpha h}{6}
  \left[\begin{array}{lll} 
      1 & 2 & 1 
    \end{array} \right] 
  \left\{\begin{array}{c}
      c_{i-1} \\
      c_{i} \\
      c_{i+1}
    \end{array}\right\}
  + \frac{D}{h}
  \left[\begin{array}{lll} 
      -1 & 2 & -1 
    \end{array} \right] 
  \left\{\begin{array}{c}
      c_{i-1} \\
      c_{i} \\
      c_{i+1}
    \end{array}\right\}
\end{align}
where $D$ denotes the diffusivity coefficient. Since $\alpha \geq 0$, $D > 0$ and 
$h > 0$; the conditions on positive diagonal entries and strict diagonal dominance 
are satisfied automatically. The condition on non-positive non-diagonal entries 
yields the following equation:
\begin{align}
  \label{Eqn:Decay_1D_mesh_constraint}
  h \leq \sqrt{\frac{6 D}{\alpha}}
\end{align}
In Figure \ref{Fig:Decay_1D_h_critical}, we compare the above theoretical prediction with 
numerical simulations for various values of $\alpha$, and the prediction is excellent.

\begin{remark}
  For the case of pure diffusion (that is, $\alpha = 0$), equation \eqref{Eqn:Decay_1D_mesh_constraint} 
  implies that any uniform mesh using two-node linear finite elements satisfies the maximum principle. 
  However, for diffusion with decay, there is a constraint on the mesh size $h$, which is proportional to 
  $D/\alpha$. That is, for a fixed $D$, the element size has to decrease with an increase in the decay 
  coefficient to meet the maximum principle. This result highlights one of the main differences between 
  diffusion with decay and pure diffusion under the classical Galerkin formulation.
\end{remark}

\subsection{Two-dimensional problem with isotropic medium}
\begin{figure}[h]
  \centering
  \includegraphics[scale=0.55]{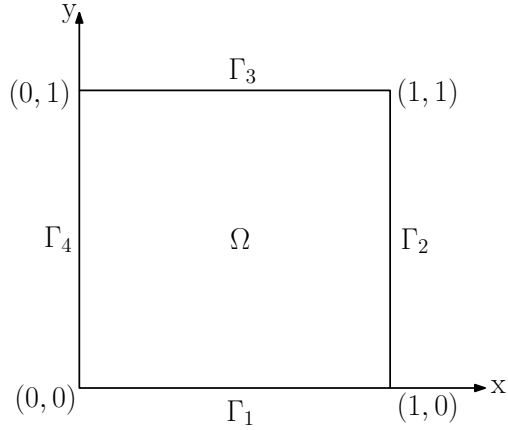}
  \caption{Two-dimensional problem with isotropic medium: The forcing function is taken to be 
    zero. The analytical solution is given by $c(\mathrm{x},\mathrm{y})=(1/\exp (\sqrt{\alpha}))
    (\exp (\sqrt{\alpha} \mathrm{x}) + \exp (\sqrt{\alpha} \mathrm{y}))$. The Dirichlet boundary 
    conditions are $c(\mathrm{x},0)=(1/\exp (\sqrt{\alpha}))(\exp (\sqrt{\alpha} \mathrm{x})+1)$ 
    on $\Gamma_{1}$, $c(1,\mathrm{y}) = 1 + \exp (\sqrt{\alpha}(\mathrm{y}-1))$ on $\Gamma_{2}$, 
    $c(\mathrm{x},1) = \exp(\sqrt{\alpha} (\mathrm{x}-1))+1$ on $\Gamma_{3}$, and $c(0,\mathrm{y}) 
    = (1/\exp (\sqrt{\alpha}))(1+\exp (\sqrt{\alpha} \mathrm{y}))$ on $\Gamma_{4}$.} 
  \label{Fig:Helmholtz_2D_analytical}
\end{figure}

\begin{figure}[h]
  \centering
  \subfigure[Three-node triangular mesh with 512 elements]{
    \includegraphics[scale=0.335]{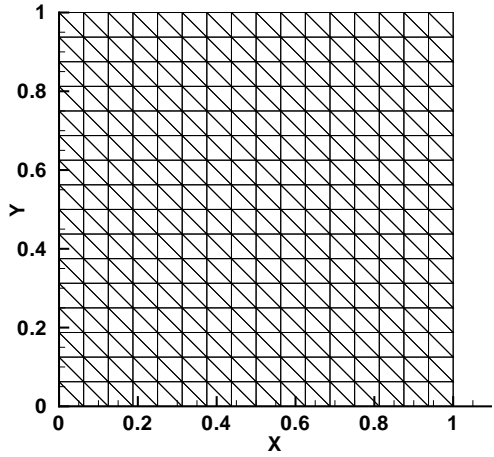}}        
  \subfigure[Classical Galerkin formulation]{
    \includegraphics[scale=0.335]{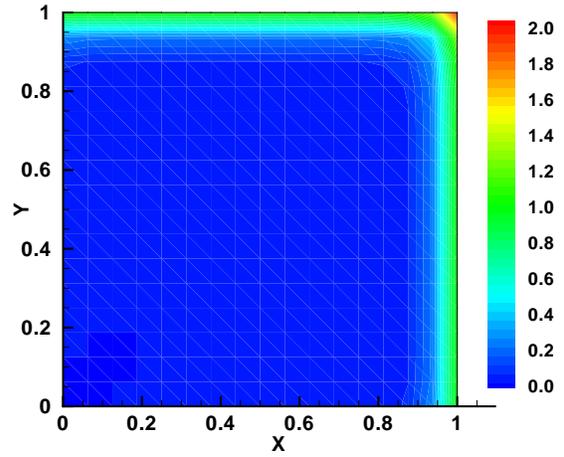}} 
  \caption{Two-dimensional problem with isotropic medium: This figure shows the 
    contours of concentration for $\alpha = 500$ on a \emph{fine computational mesh} 
    under the Galerkin formulation, and there is no violation of the maximum principle}
  \label{Fig:Decay_2D_fine_mesh}
\end{figure}

\begin{figure}      
  \subfigure[Unstructured mesh]{
    \includegraphics[scale=0.335]{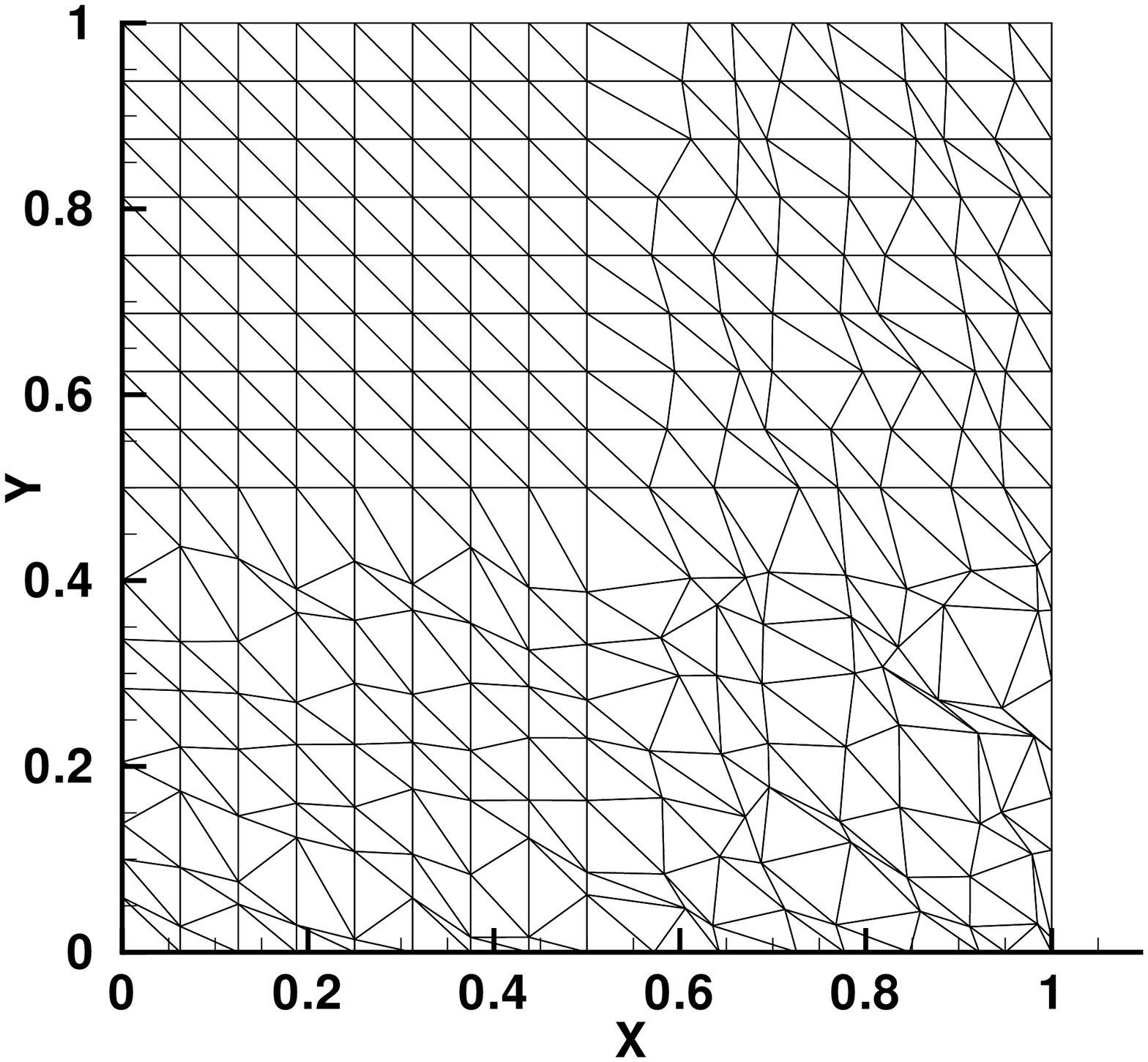}}        
  \subfigure[Classical Galerkin formulation]{
    \includegraphics[scale=0.335]{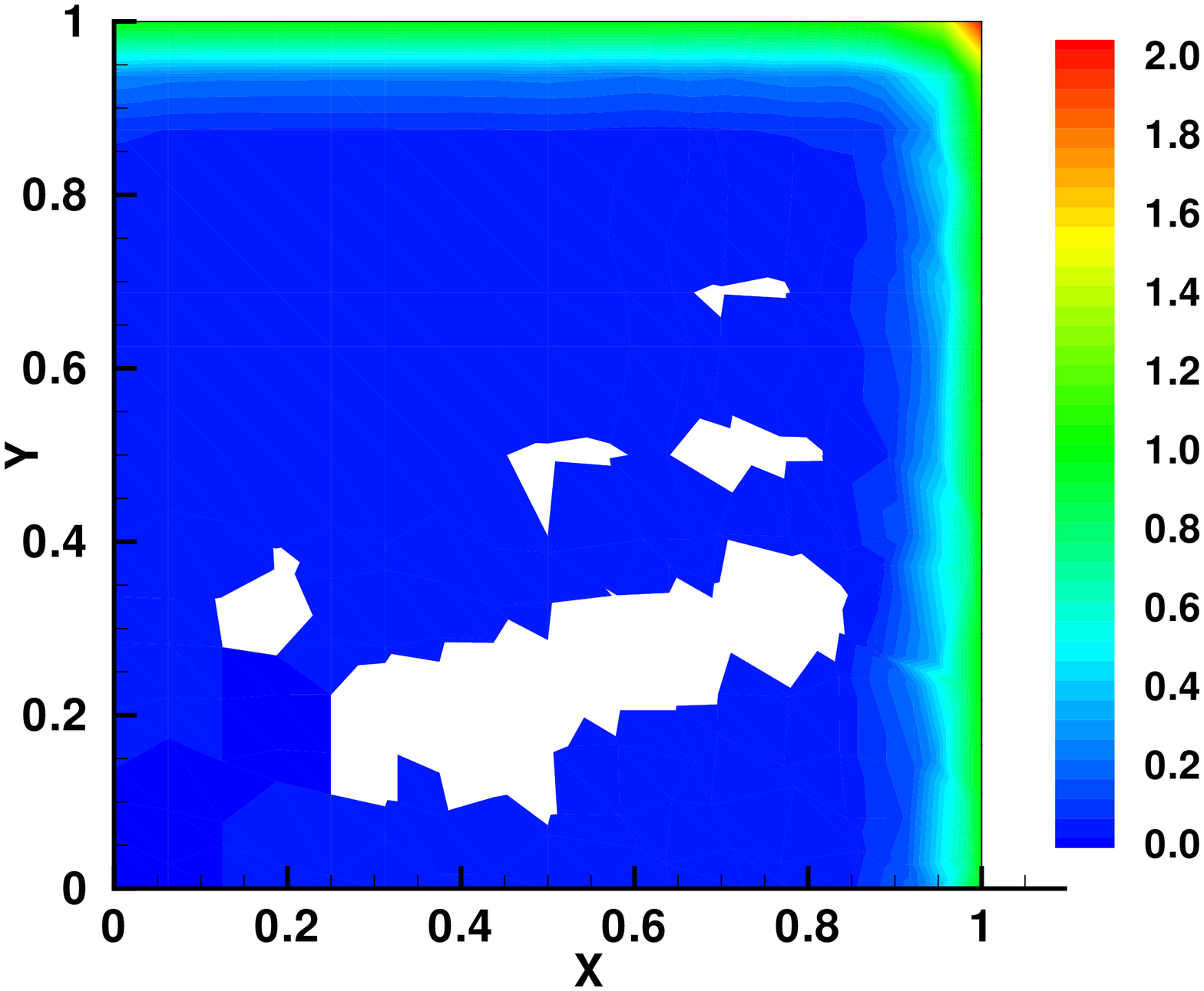}} 
  \subfigure[Analytical solution]{
    \includegraphics[scale=0.335]{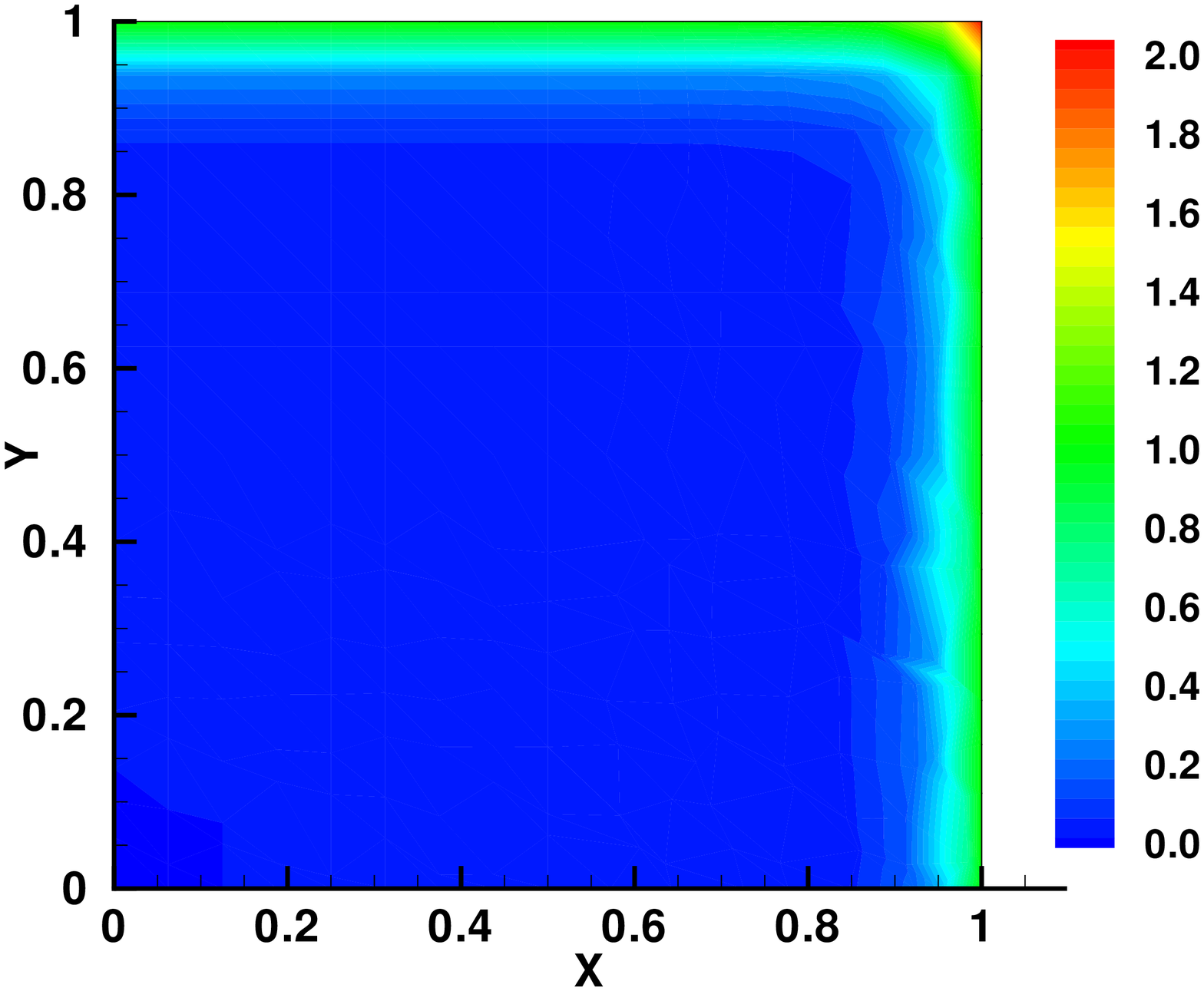}} 
  \subfigure[Proposed formulation]{
    \includegraphics[scale=0.335]{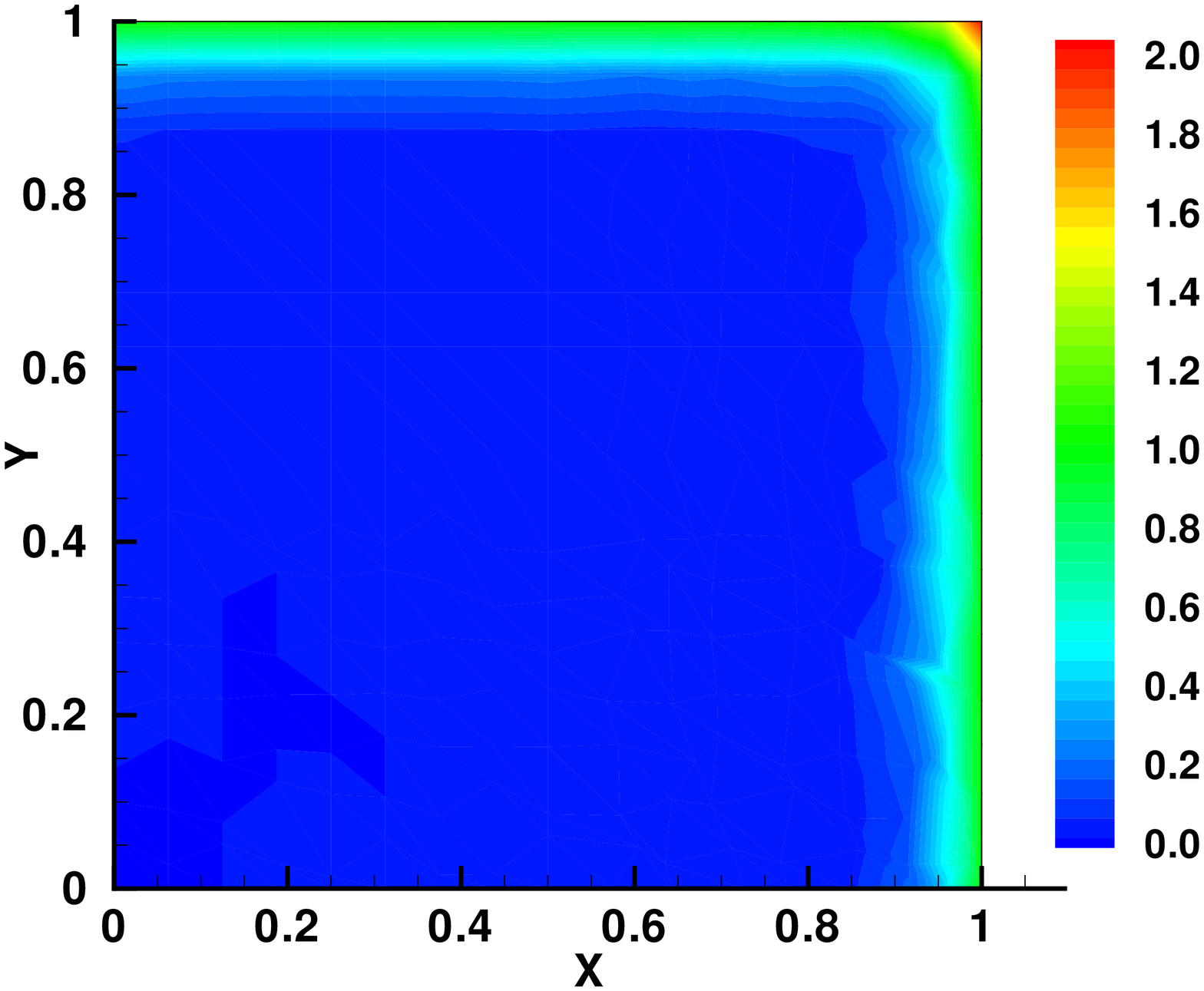}} 
  \caption{Two-dimensional problem with isotropic medium: The problem is solved on a fine 
    unstructured mesh using the classical Galerkin formulation and the proposed formulation. 
    Analytical solution is also shown in the Figure. The decay coefficient is taken to be 
    $\alpha = 500$. Regions that have negative concentrations are indicated in white color. 
    The proposed formulation produced physically meaningful non-negative concentrations, 
    and matched well with the analytical solution. Under the classical Galerkin formulation, 
    approximately $14.2\%$ of the total number of nodes have negative concentrations. 
    Under the classical Galerkin formulation, the  minimum value of concentration is  
    $-0.0466$. Note that the negative concentration occurred mostly in the perturbed 
    mesh region.} \label{Fig:Decay_2D_perturbed_mesh}
\end{figure}

Consider the following two-dimensional problem with homogeneous forcing function:
\begin{align}
  \label{Eqn:Helmholtz_2D_equation}
  \alpha c(\mathrm{x},\mathrm{y})- \mathrm{div}[\mathbb{D} \; \mathrm{grad}[c(\mathrm{x},
  \mathrm{y})]] = 0 \quad \mathrm{in} \; \Omega := (0, 1) \times (0, 1) 
\end{align}
with $\mathbb{D}$ is assumed to be the identity tensor (i.e, isotropic medium) and $\alpha 
\geq 0$. The geometry and boundary conditions for this two-dimensional problem are shown in 
Figure \ref{Fig:Helmholtz_2D_analytical}. The analytical solution is given by
\begin{align}
  \label{Eqn:Helmholtz_2D_analytical}
  c(\mathrm{x},\mathrm{y})= \dfrac{\exp(\sqrt{\alpha} \mathrm{x}) + 
    \exp(\sqrt{\alpha} \mathrm{y})}{\exp (\sqrt{\alpha})}
\end{align}
In this numerical study, we have taken $\alpha = 500$. In Figure \ref{Fig:Decay_2D_coarse_mesh}, 
we show the numerical results obtained using the classical Galerkin formulation and the proposed 
formulation on a coarse computational mesh. The classical Galerkin formulation violates the maximum 
principle, and the proposed formulation produces physically meaningful non-negative concentration 
even on the chosen coarse computational mesh. 
Since the diffusivity tensor is isotropic and the mesh is based on right-angled isosceles 
triangles, the violation of maximum principle under the Galerkin formulation (if it occurs) 
is due to the decay term. Moreover, the violation vanishes with sufficient mesh refinement. 
This fact is illustrated in Figure \ref{Fig:Decay_2D_fine_mesh} wherein we have employed 
a finer computational mesh, and the Galerkin formulation satisfies the maximum principle. 
However, it should be noted that the classical Galerkin formulation violates the maximum 
principle on fine unstructured computational meshes, which is illustrated in Figure 
\ref{Fig:Decay_2D_perturbed_mesh}. 

\begin{figure}[!h]
  \centering
  \subfigure{
    \includegraphics[scale=0.45]{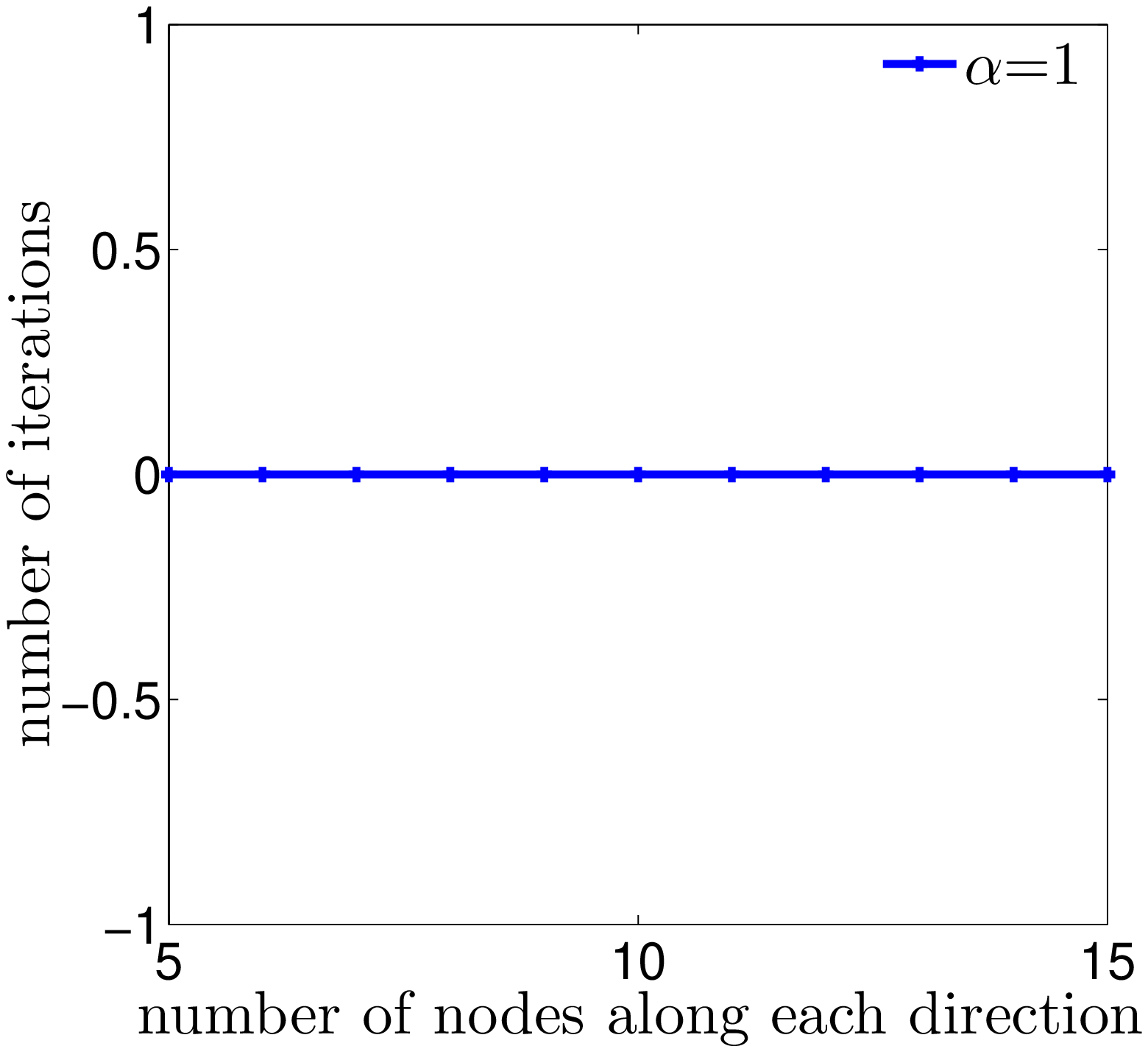}}        
  \subfigure{
    \includegraphics[scale=0.45]{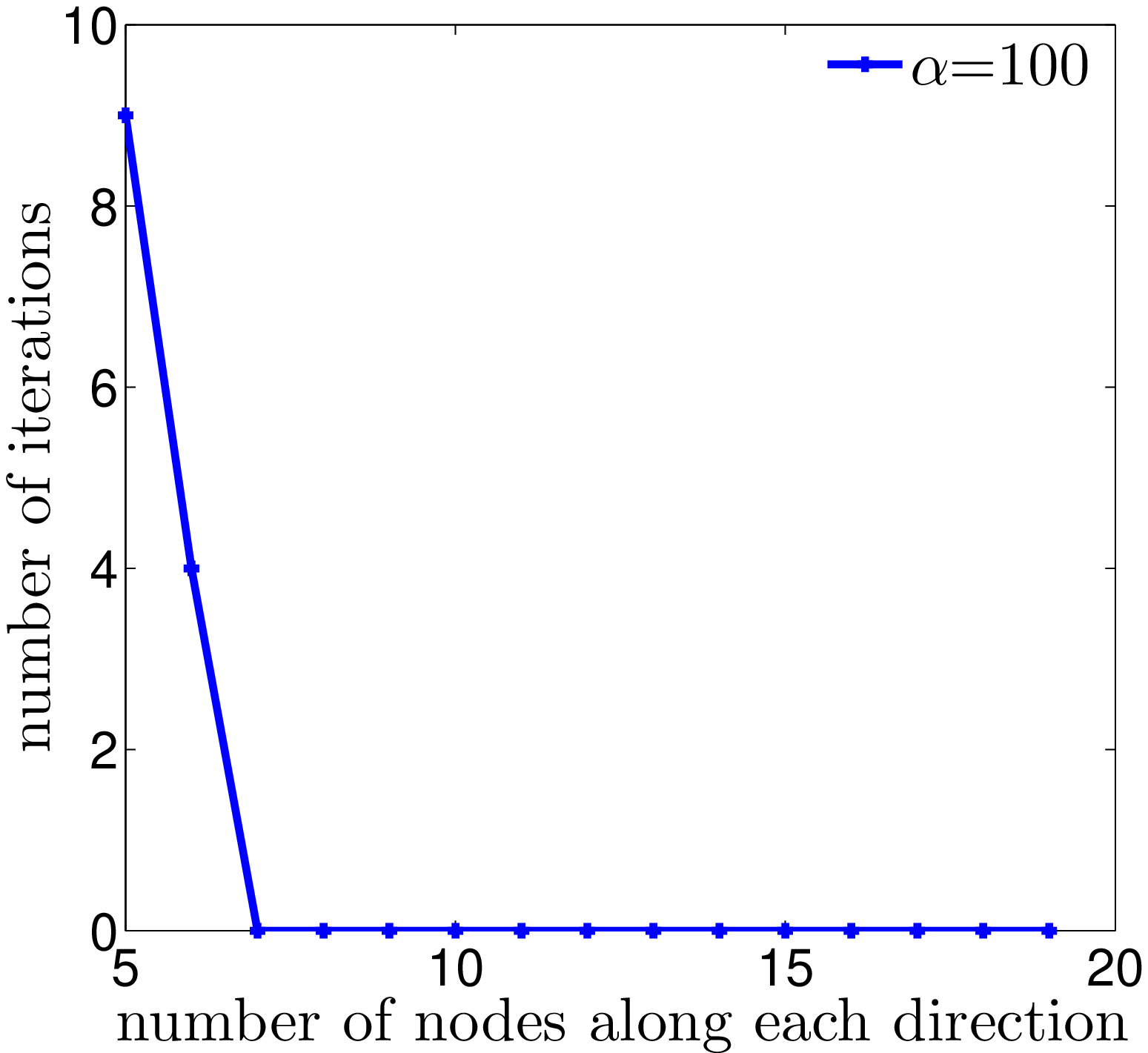}} 
  \subfigure{
    \includegraphics[scale=0.45]{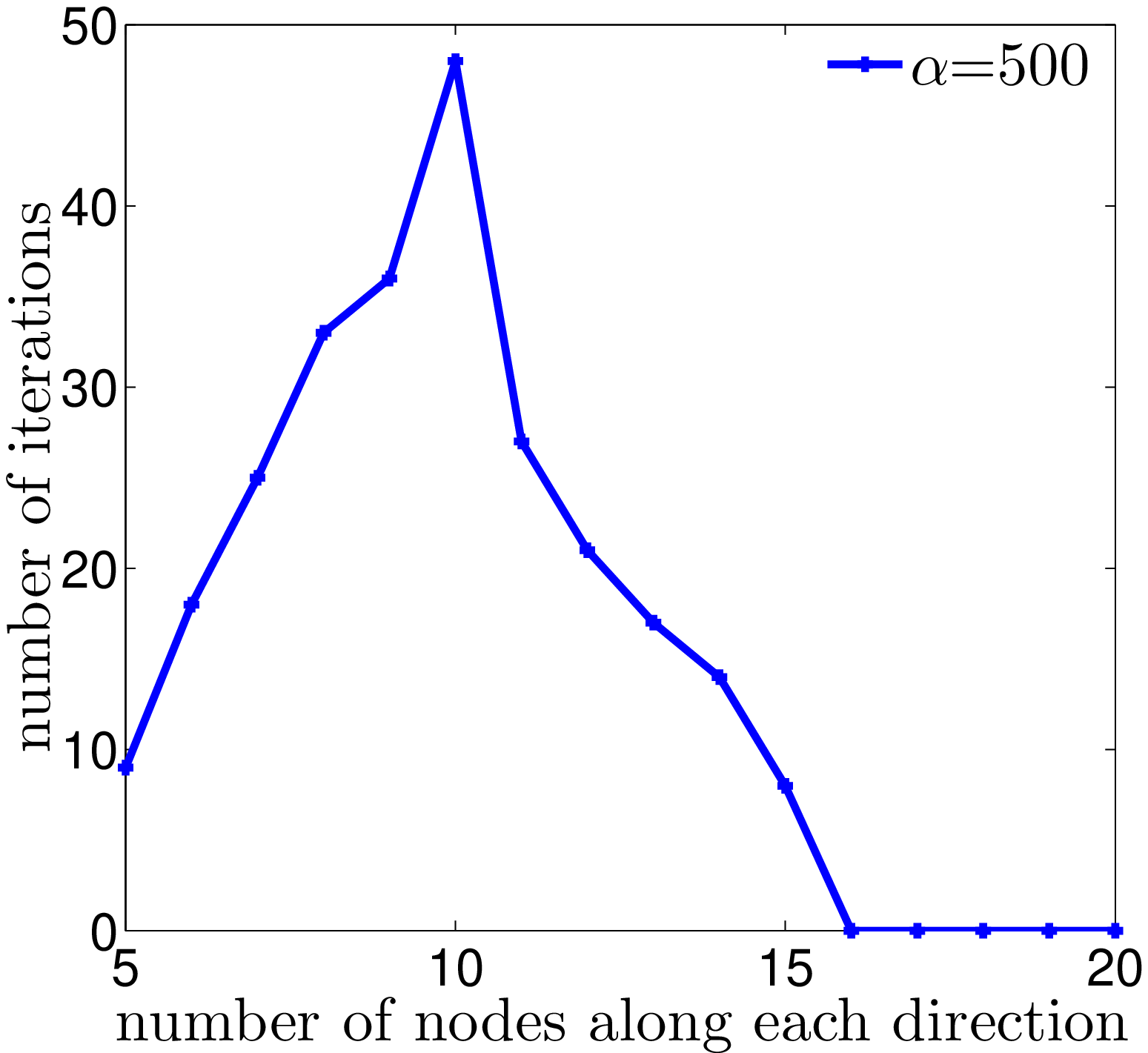}} 
  \subfigure{
  \includegraphics[scale=0.45]{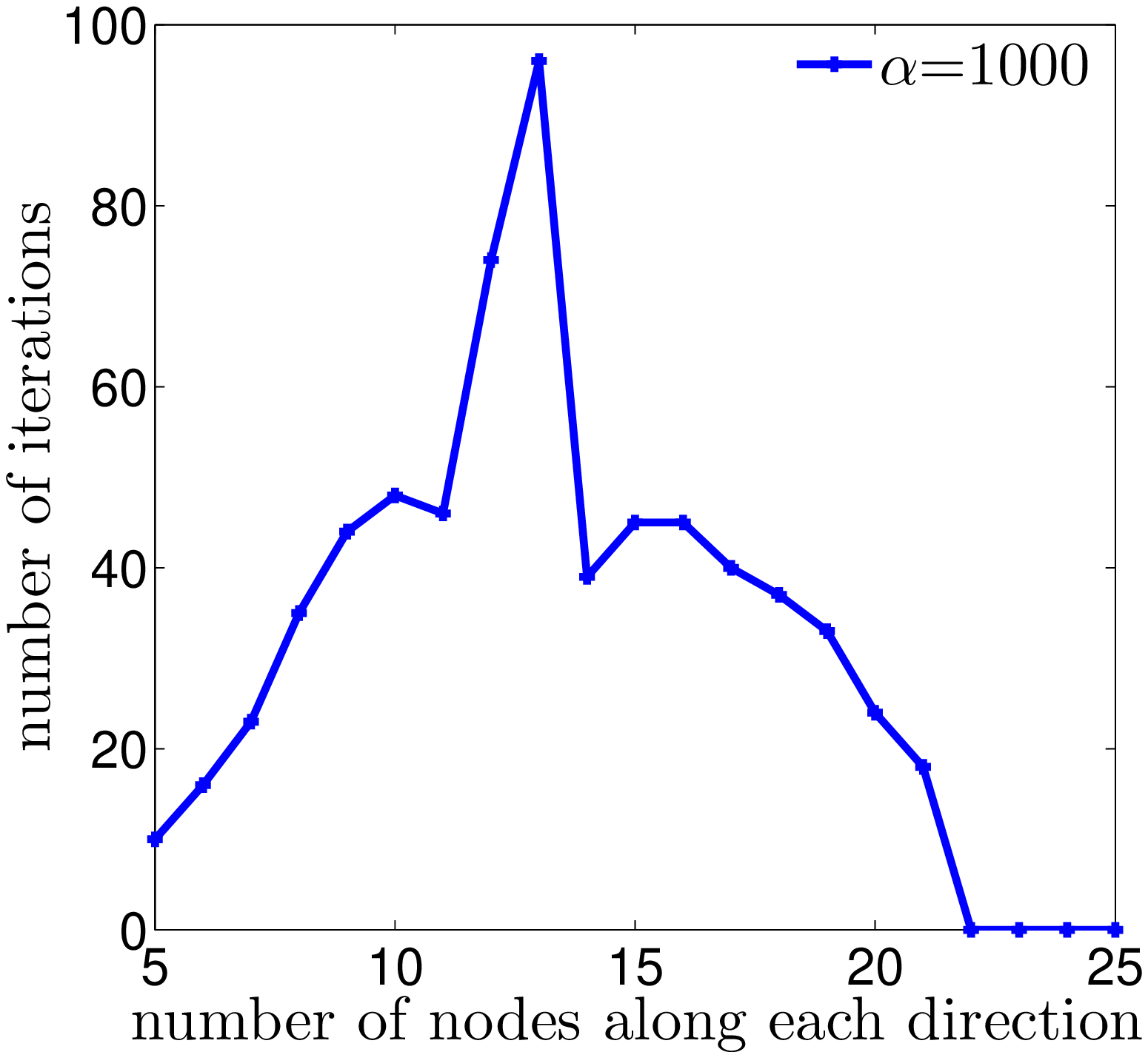}} 
  \caption{Two-dimensional problem with isotropic medium: These figures present the 
    number of iterations required for the proposed formulation using active-set 
    strategy at various values of $\alpha$ with respect to the number of nodes along 
    each side of the computational domain (which is same in both x and y directions). 
    Note that the number of iterations required for the active-set strategy to terminate 
    increases as $\alpha$ increases. Again for this case, there is no violation of the 
    discrete maximum principle after sufficient mesh refinement, and there is no need 
    to solve the constrained optimization problem.} \label{Fig:Decay_2D_active_set_iterations}
\end{figure}

\begin{figure}[!h]
  \centering
      \subfigure[$\alpha = 1$]{
      \includegraphics[scale=0.335]{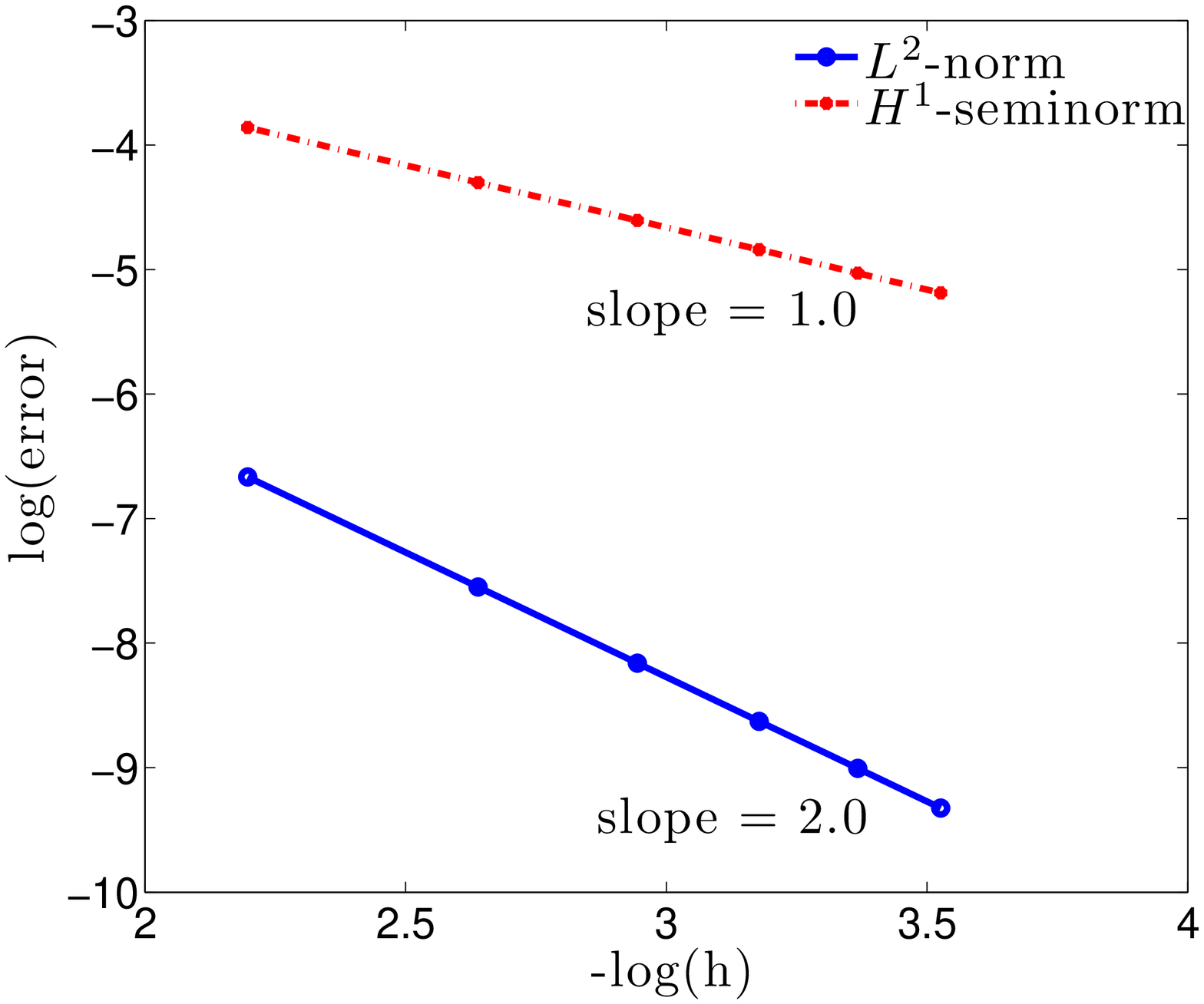}}        
      \subfigure[$\alpha = 100$]{
      \includegraphics[scale=0.335]{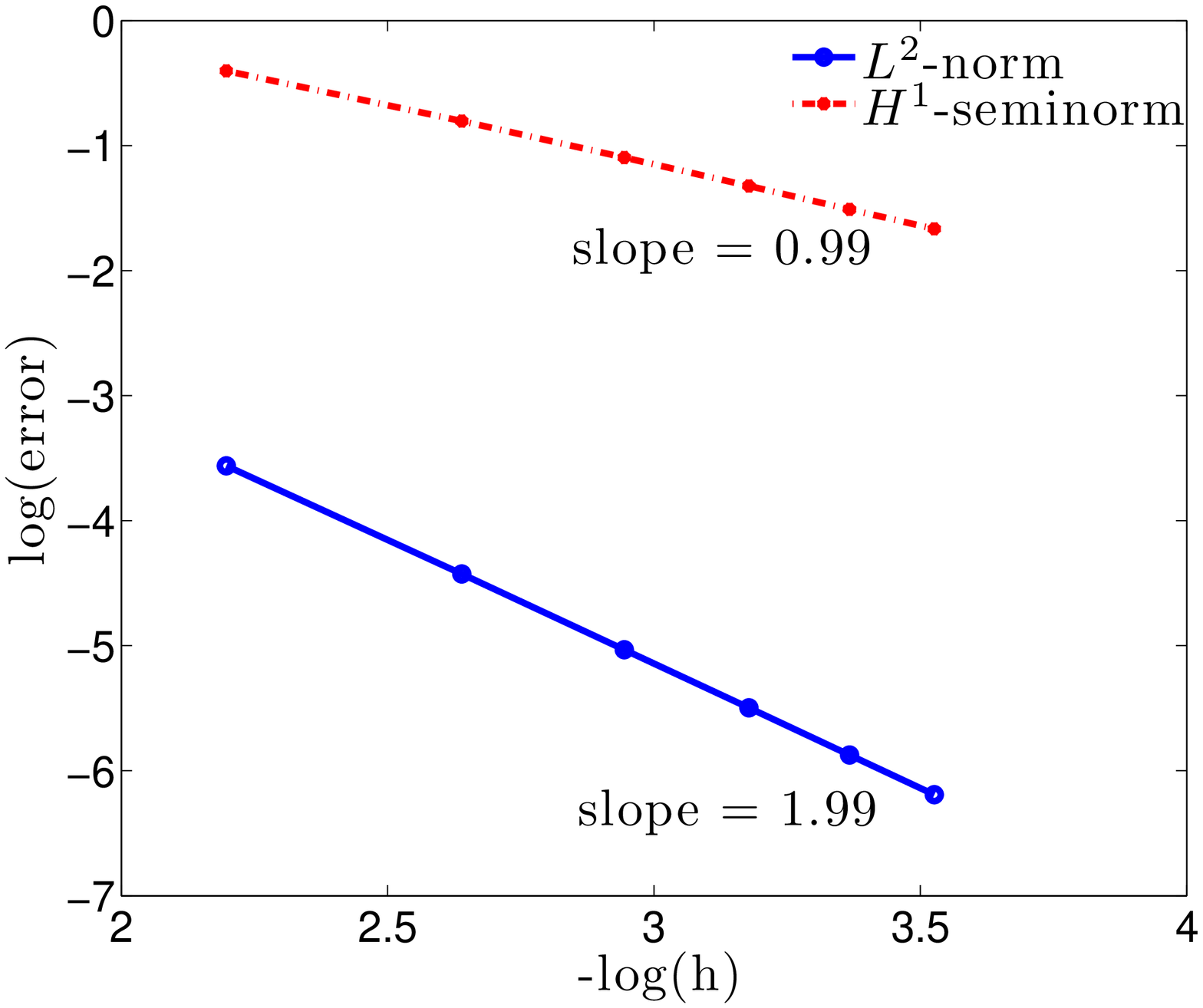}} 
      \subfigure[$\alpha = 500$]{
      \includegraphics[scale=0.335]{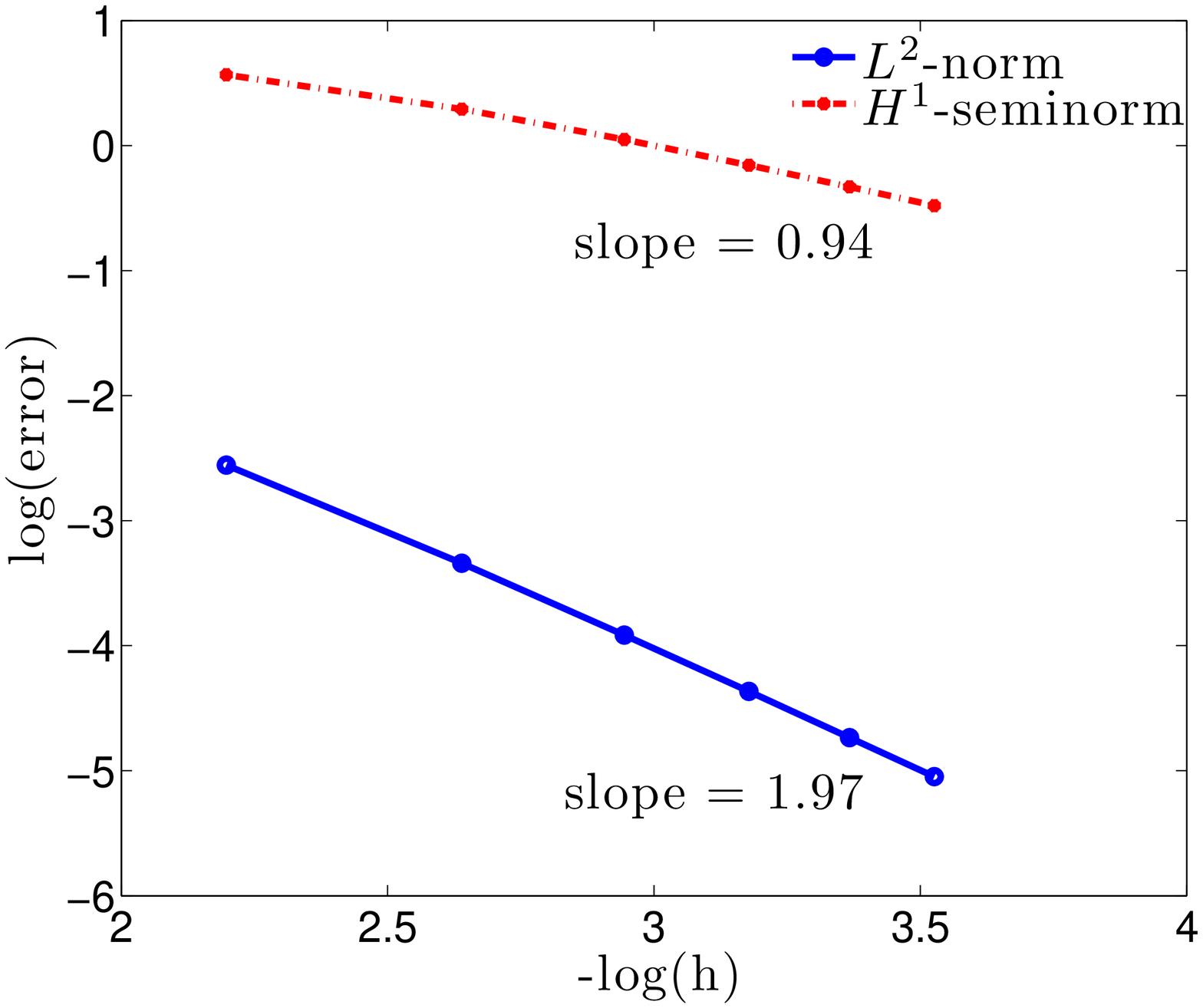}} 
      \subfigure[$\alpha = 1000$]{
      \includegraphics[scale=0.335]{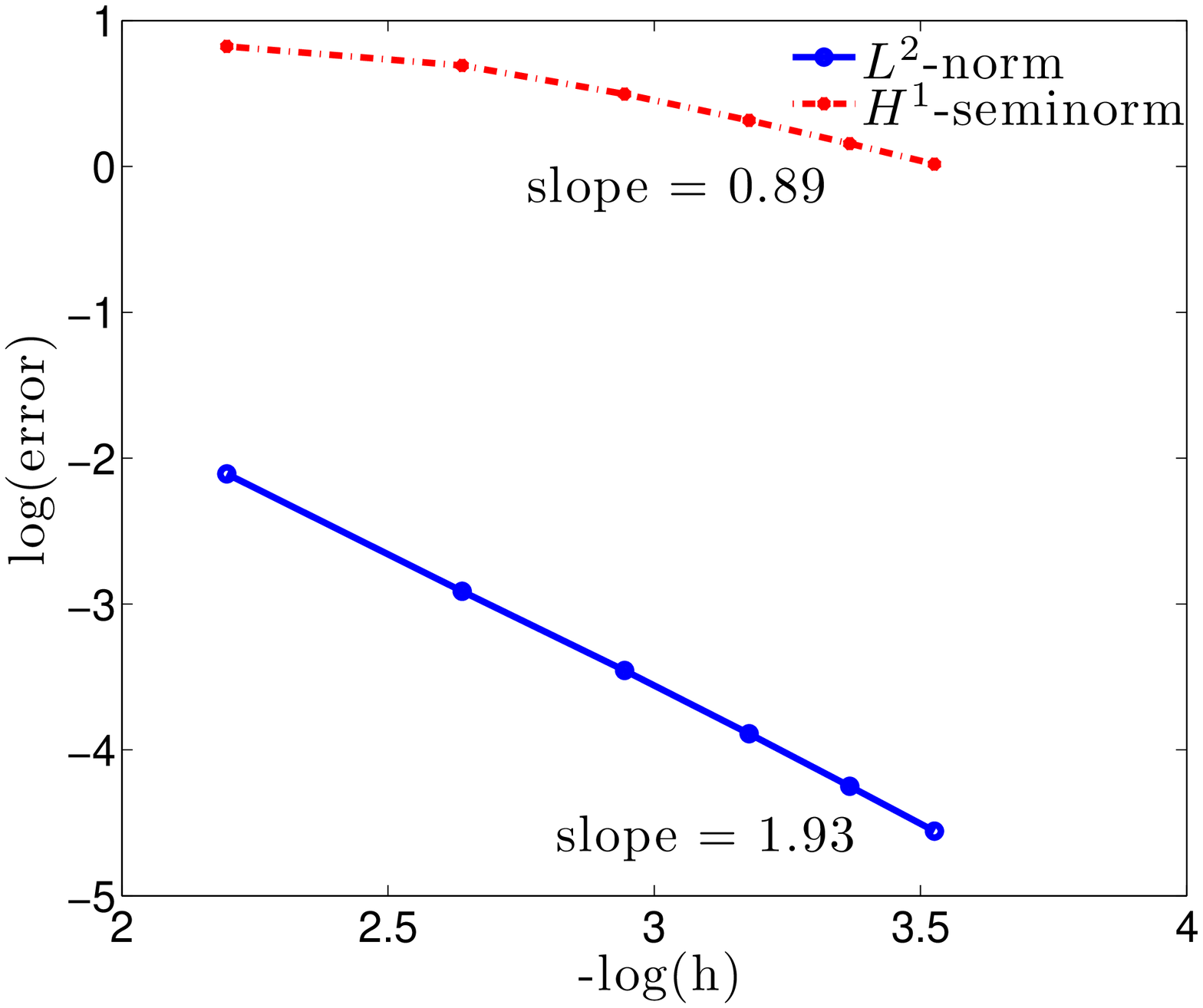}} 
\caption{Two-dimensional problem with isotropic medium: This figure presents numerical 
  convergence of the proposed formulation with mesh refinement for various values of decay 
  coefficient. From the figure it is evident that the rates of convergence with respect to mesh 
  refinement in $L^{2}$-norm and $H^{1}$-seminorm are about the same as for the original 
  linear finite element method.} \label{Fig:Decay_error_convergence_2D}
\end{figure}

In Figure \ref{Fig:Decay_2D_active_set_iterations}, we plot the number of iterations taken by 
the active-set strategy with respect to number of (finite element) nodes. As discussed earlier, 
since the medium is isotropic, the violation of the discrete maximum principle again decreases 
with mesh refinement, and eventually there will no violation of the discrete maximum principle. 
This can be seen in Figure \ref{Fig:Decay_2D_active_set_iterations} as the number of 
active-set strategy iterations is zero for (sufficiently) finer computational meshes for various 
values of $\alpha$.
In Figure \ref{Fig:Decay_error_convergence_2D}, we perform numerical convergence studies of 
the proposed formulation for various values of $\alpha$, and the proposed formulation performed 
well. 

\subsection{Two-dimensional problems with anisotropic medium}
Consider anisotropic diffusion in a bi-unit square plate $\Omega = (0, 1) \times (0, 1)$. 
The anisotropic diffusivity tensor is taken as follows:
\begin{align}
  \label{Eqn:Decay_diffusivity_tensor_homogeneous}
  \mathbb{D} = 
  \left(\begin{array}{cc}
      \cos(\theta) & \sin(\theta) \\ 
      -\sin(\theta) & \cos(\theta)
    \end{array}\right)
  \left(\begin{array}{cc}
      k_1 & 0 \\ 
      0 & k_2
    \end{array}\right)
  \left(\begin{array}{cc}
      \cos(\theta) & -\sin(\theta) \\ 
      \sin(\theta) & \cos(\theta)
    \end{array}\right)
\end{align}
with $\theta = \pi/6$, $k_1 = 10000$ and $k_2 = 1$. The forcing function is taken to 
be zero (that is, $f(\mathrm{x}, \mathrm{y}) = 0)$, and the decay coefficient is taken 
as $\alpha = 1$. We prescribed Dirichlet boundary conditions on the whole boundary, and 
the prescribed concentrations are as follows: the left, right and top sides of the 
computational domain have a prescribed concentration of zero (that is, $c^{\mathrm{p}}
(\mathrm{x} = 0,\mathrm{y}) = c^{\mathrm{p}}(\mathrm{x} = 1, \mathrm{y}) = c^{\mathrm{p}}
(\mathrm{x}, \mathrm{y} = 1)=0$), and the bottom side of the computational domain has 
a prescribed concentration of $c^{\mathrm{p}}(\mathrm{x},\mathrm{y} = 0) = \mathrm{sin}
(\pi \mathrm{x})$. The prescribed data in this problem meet all the conditions in Theorem 
\ref{Theorem:Decay_maximum_principle}, and from the maximum principle we can infer the 
following:
\begin{align}
  0 \leq c(\mathbf{x}) \leq 1 \quad \forall \mathbf{x} \in \bar{\Omega}
\end{align}

\begin{figure}
  \centering
  \subfigure{
    \includegraphics[scale=0.3]{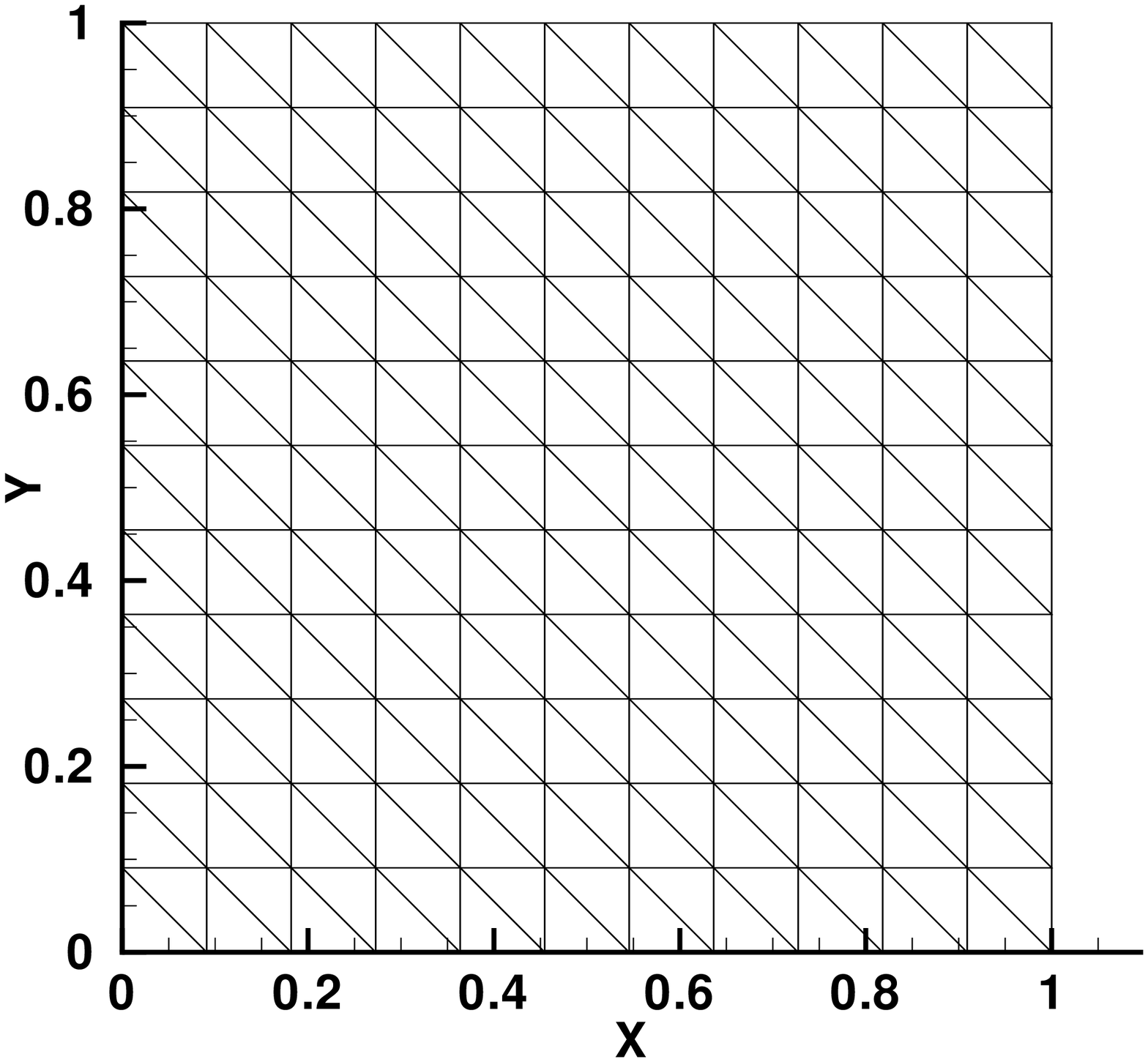}} 
  \subfigure{
    \includegraphics[scale=0.3]{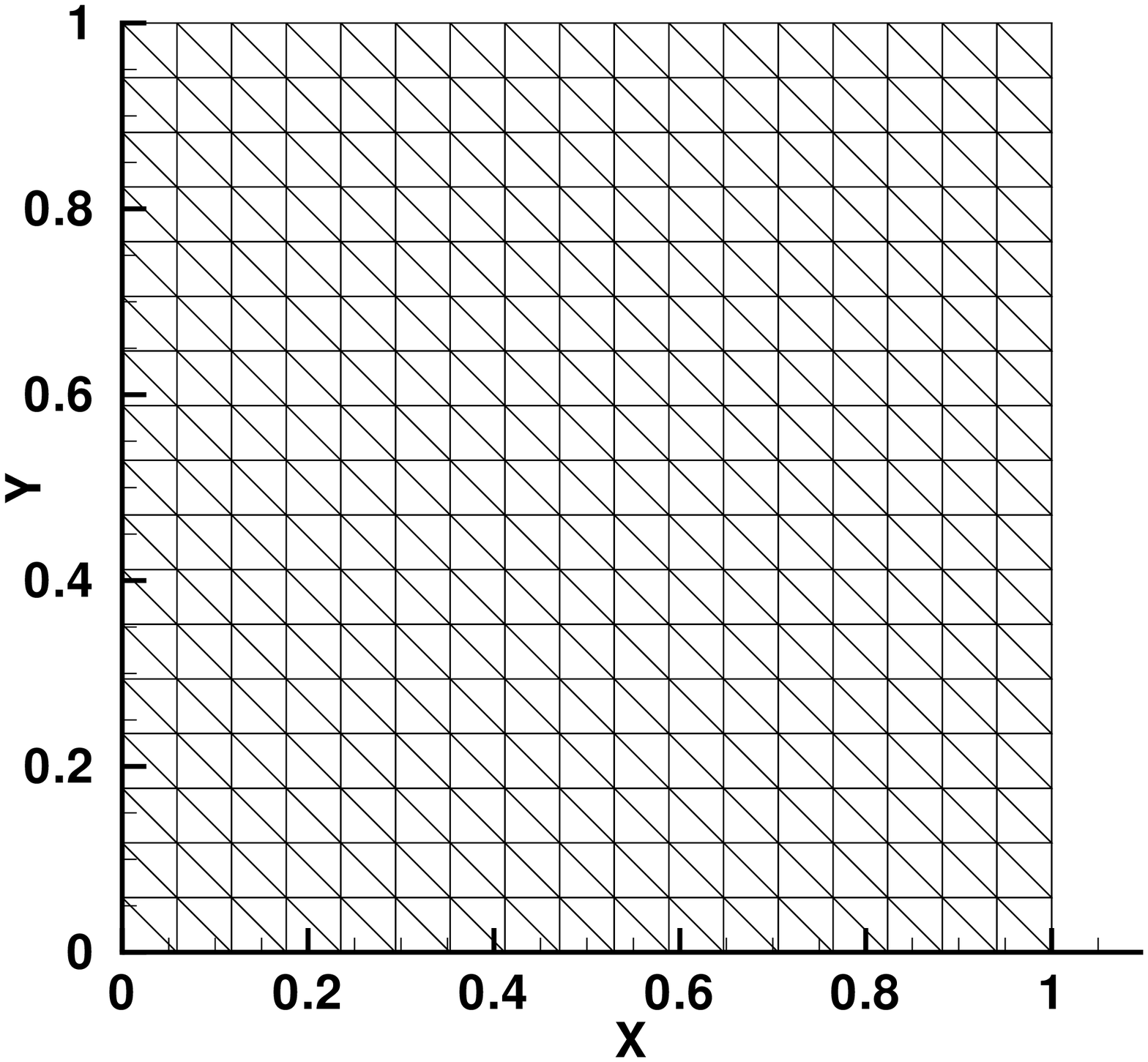}}        
  \subfigure{
    \includegraphics[scale=0.3]{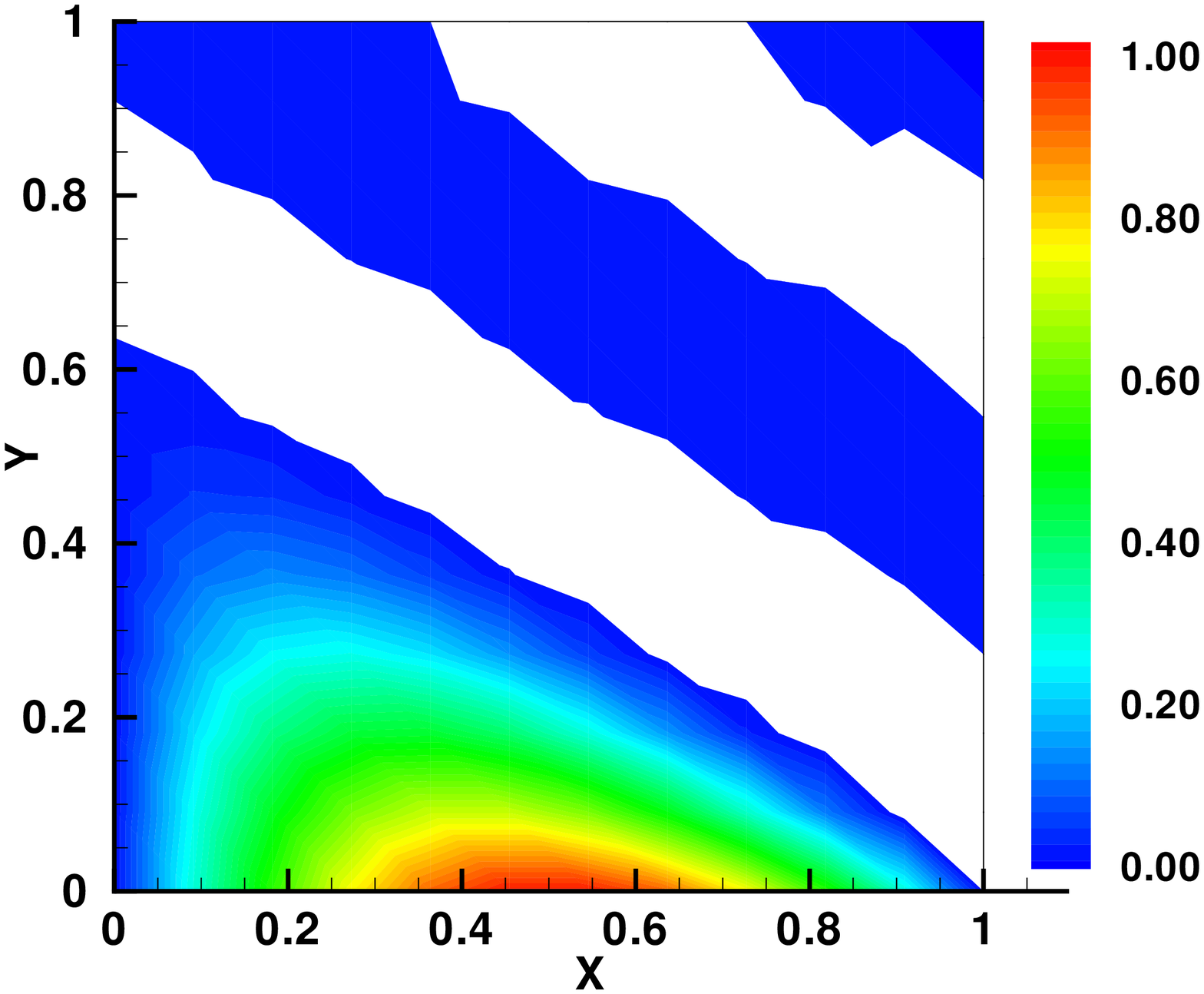}} 
  \subfigure{
    \includegraphics[scale=0.3]{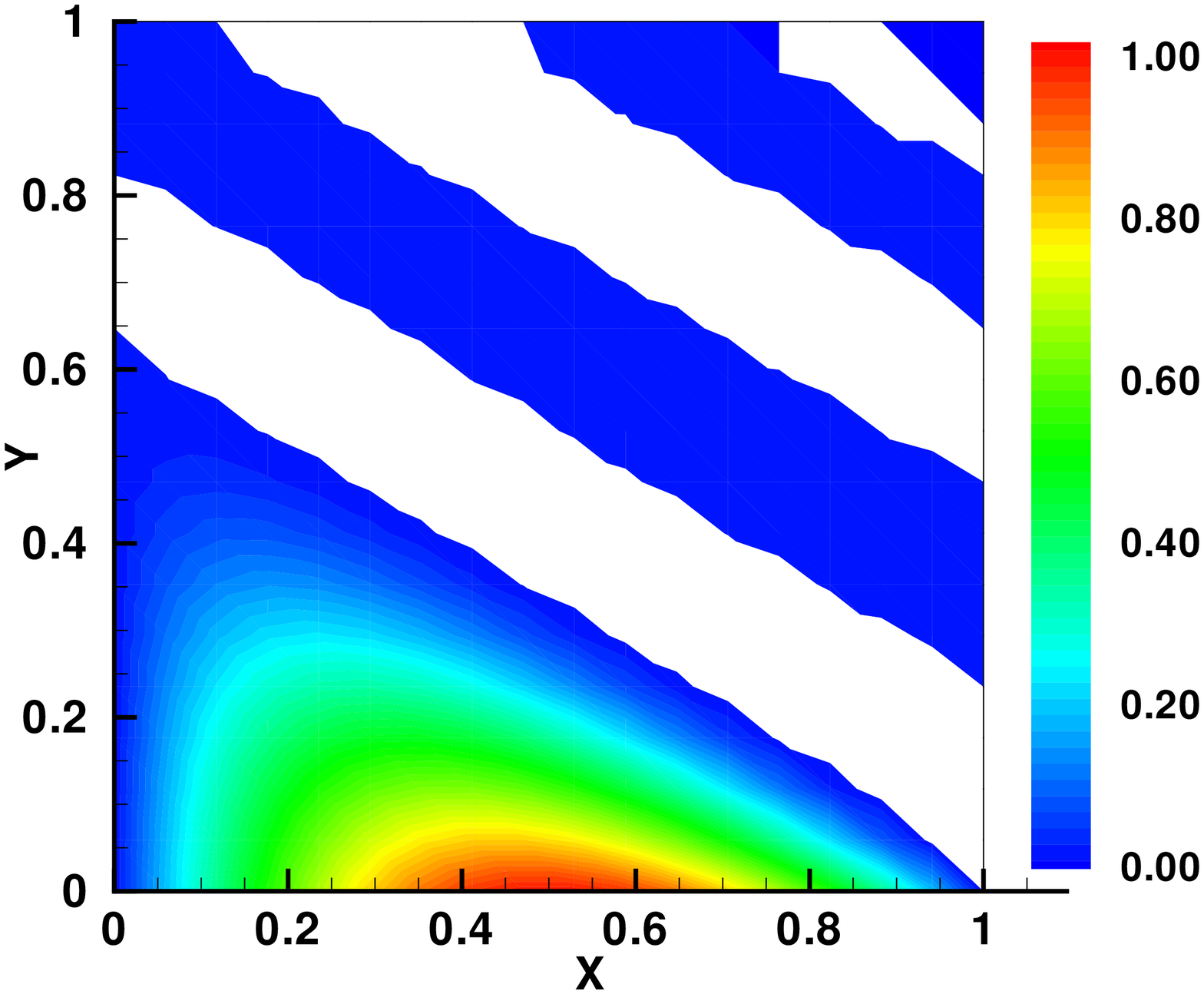}}        
  \subfigure{
    \includegraphics[scale=0.3]{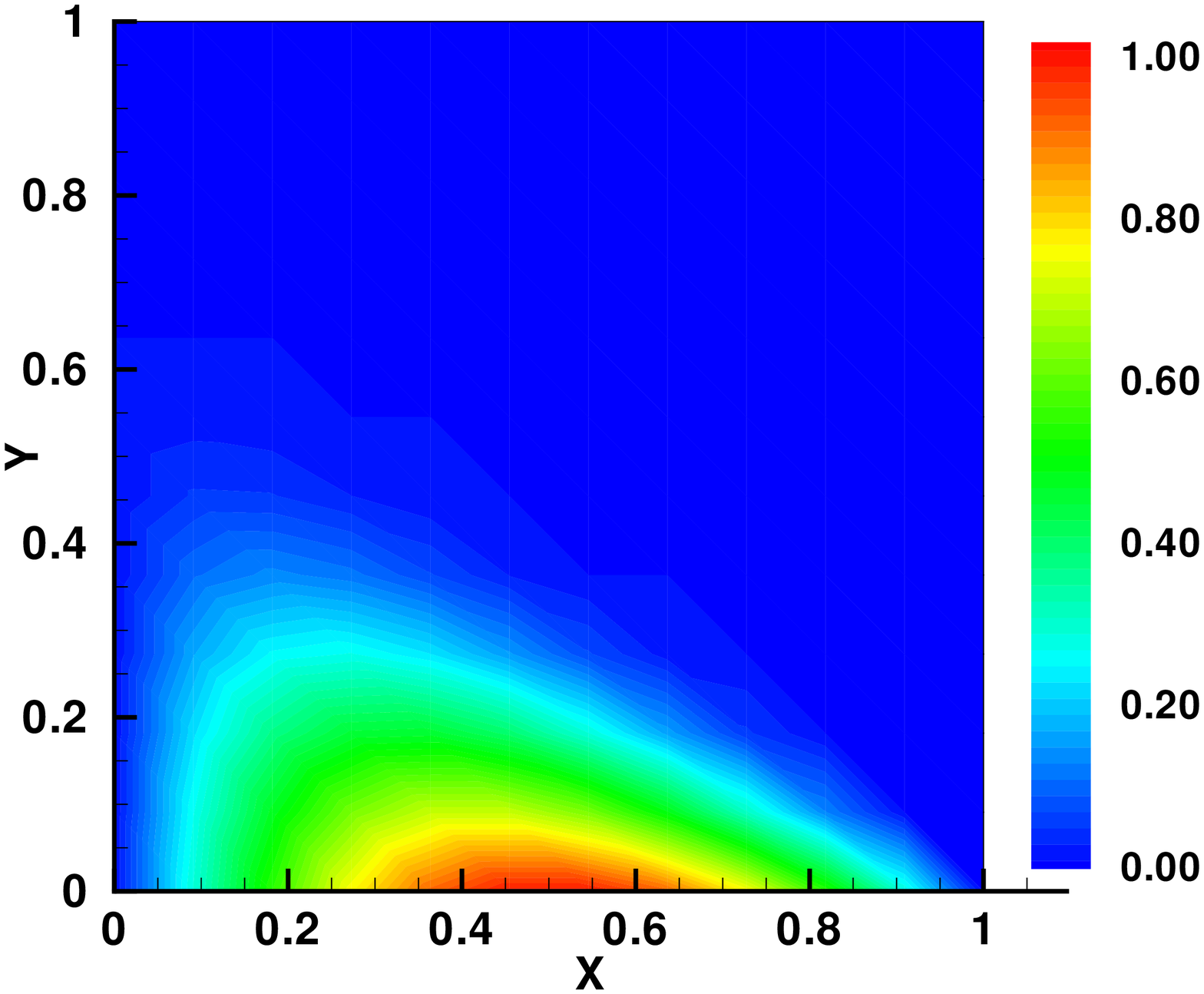}} 
  \subfigure{
    \includegraphics[scale=0.3]{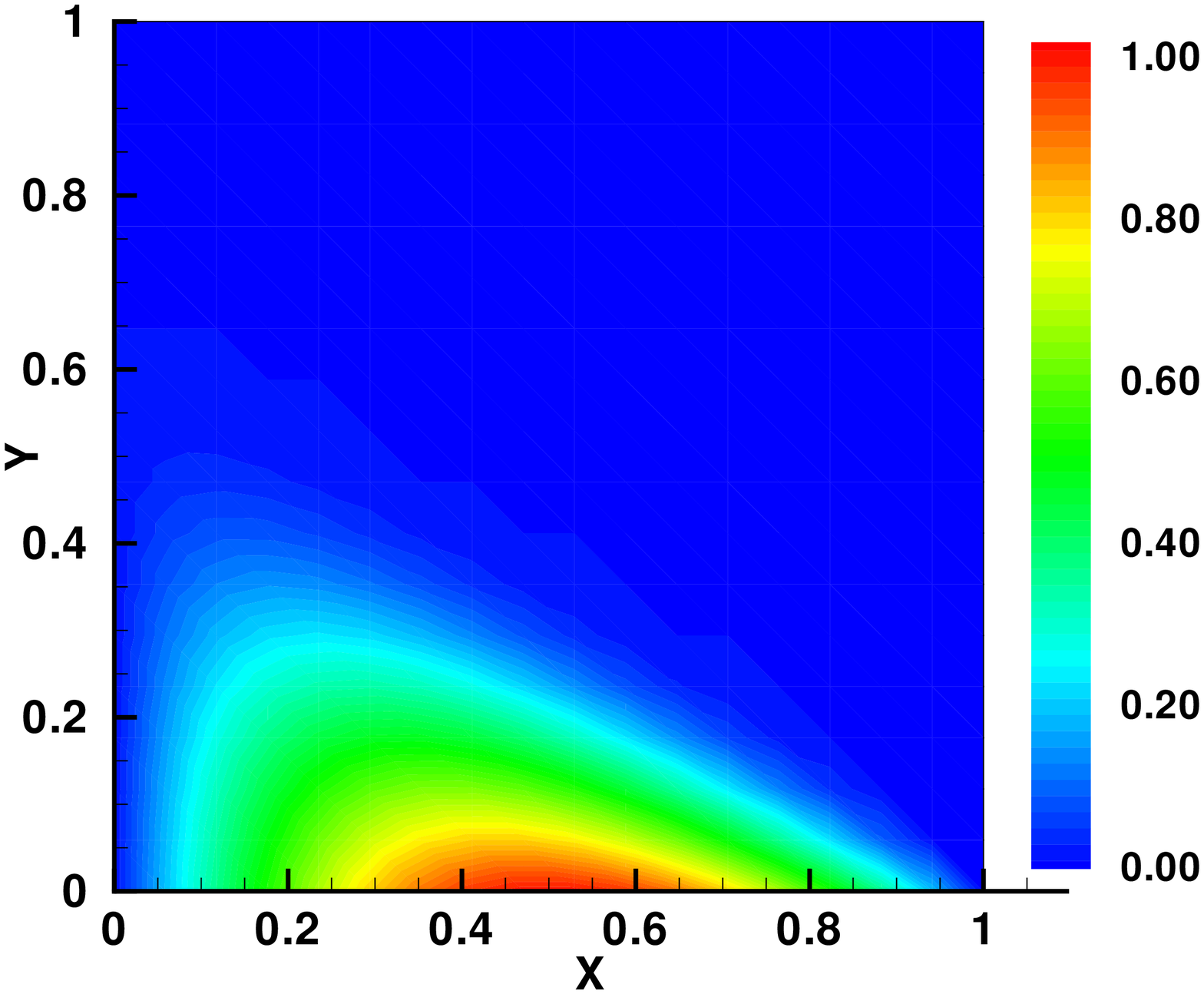}}  
  \caption{Two-dimensional problem with anisotropic medium: The problem is solved using the 
    Galerkin formulation (middle) and the proposed formulation (bottom). The left and right 
    figures are, respectively, using $12 \times 12$ and $18 \times 18$ \emph{three-node 
      triangular meshes}. Regions that have negative concentrations are indicated in white 
    color. Under the Galerkin formulation, $27.78 \%$ (for $12 \times 12$ mesh) and $30.86\%$ 
    of the total number of nodes have negative nodal concentration. The minimum concentrations 
    are $-0.035$ (for $12 \times 12$ mesh) and $-0.022$ (for $18 \times 18$ mesh).}
  \label{Fig:Decay_2D_anisotropic_T3}
\end{figure}

\begin{figure}
  \centering
      \subfigure{
      \includegraphics[scale=0.3]{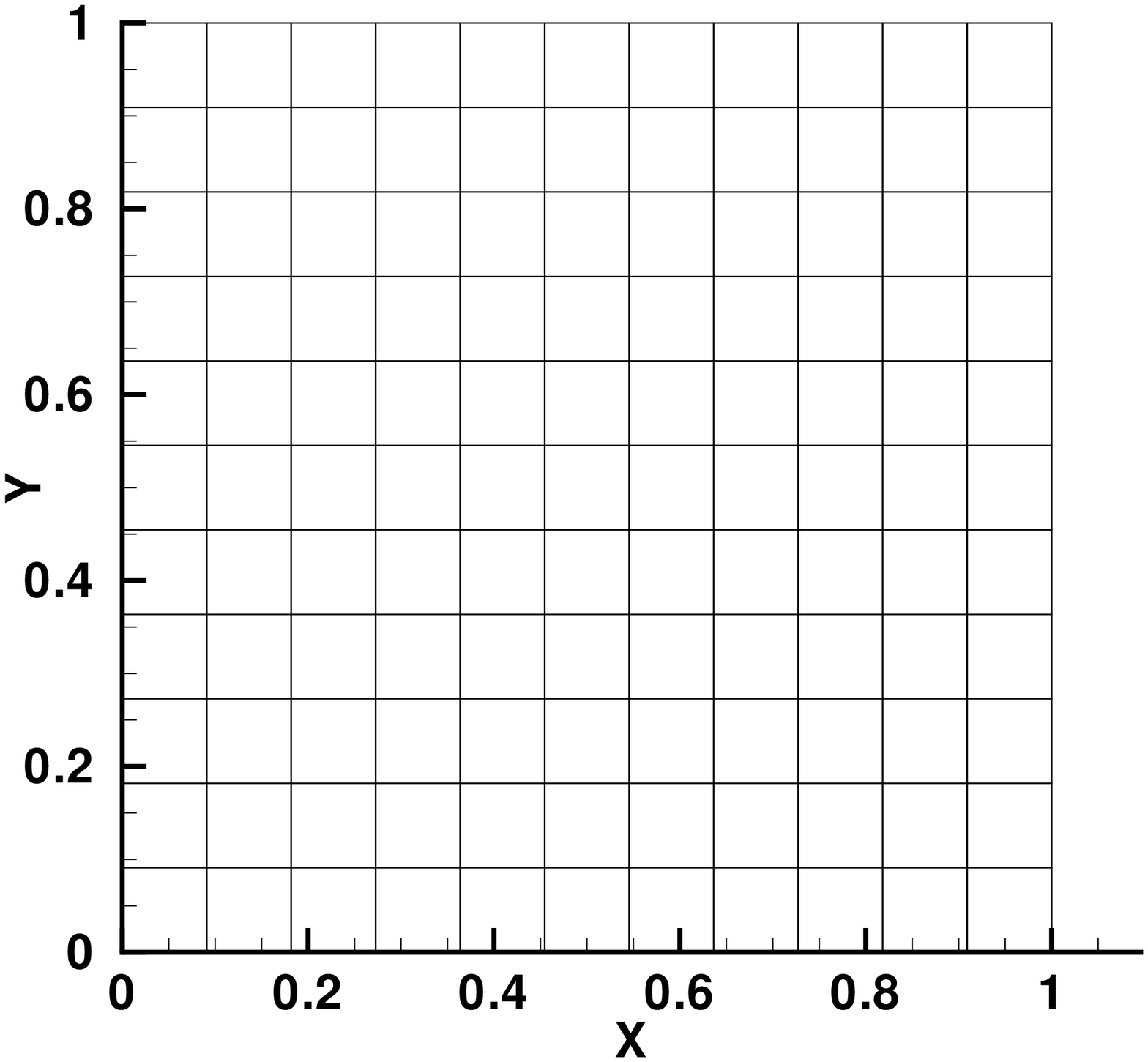}} 
      \subfigure{
      \includegraphics[scale=0.3]{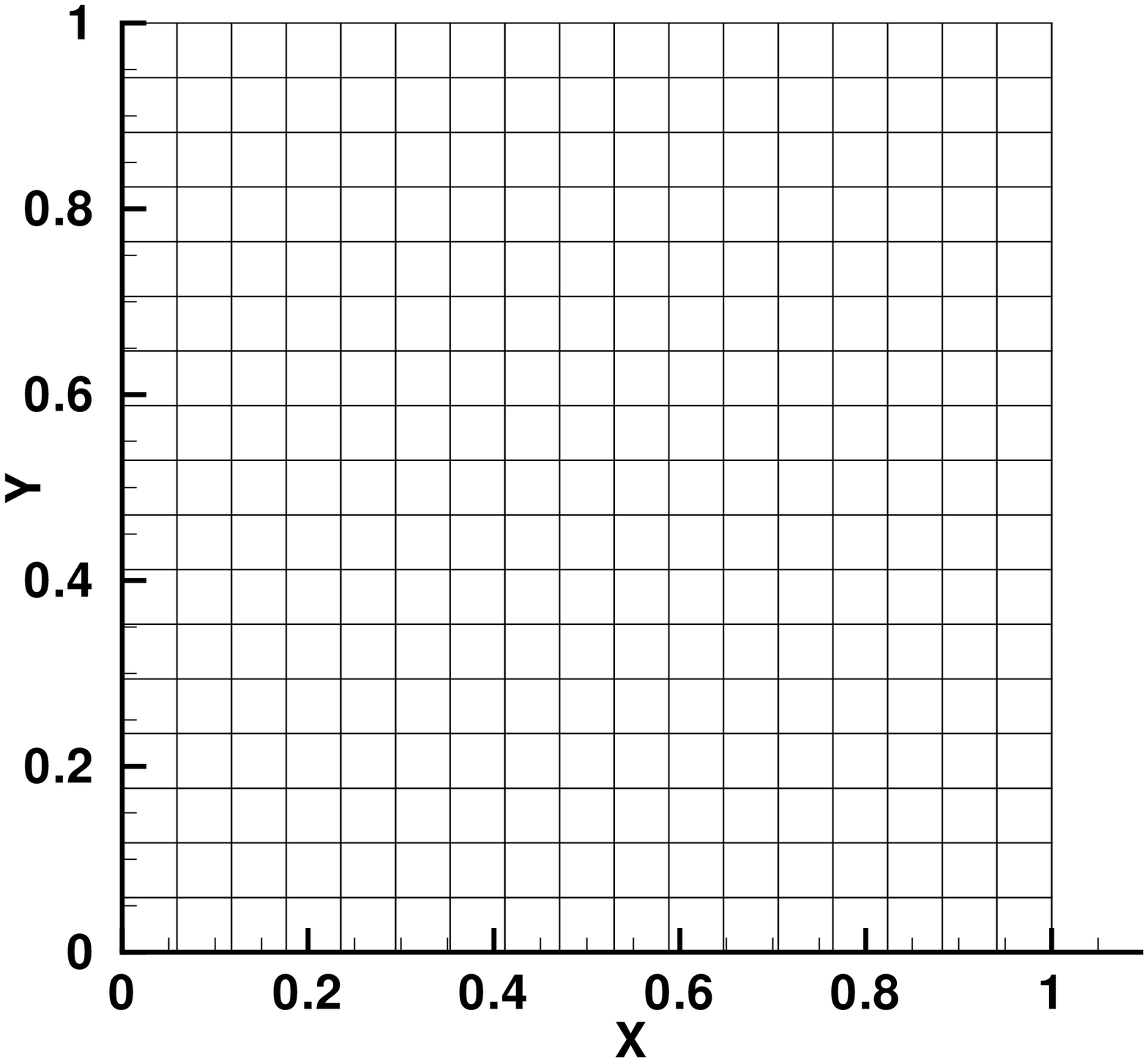}}        
       \subfigure{
      \includegraphics[scale=0.3]{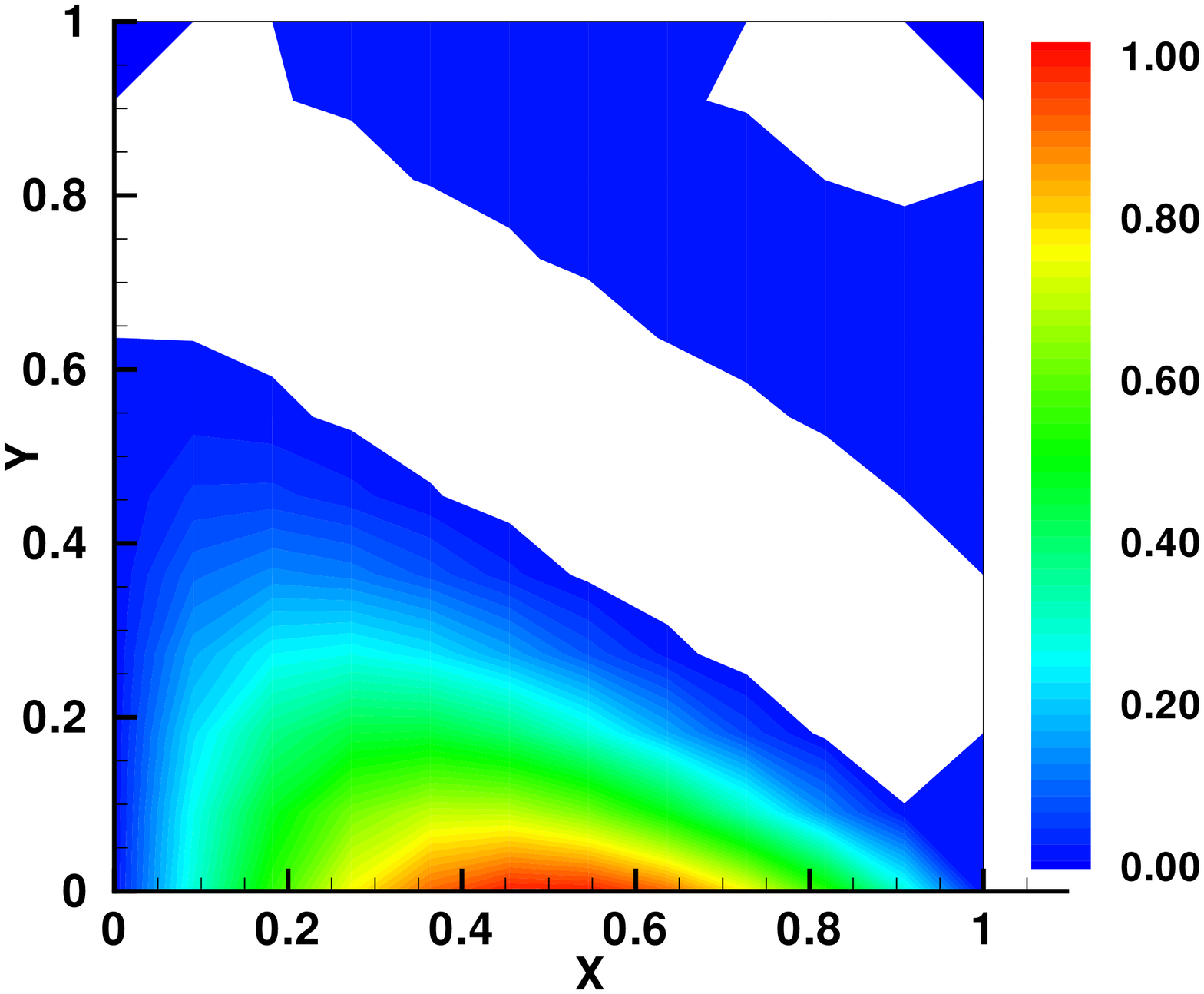}} 
      \subfigure{
      \includegraphics[scale=0.3]{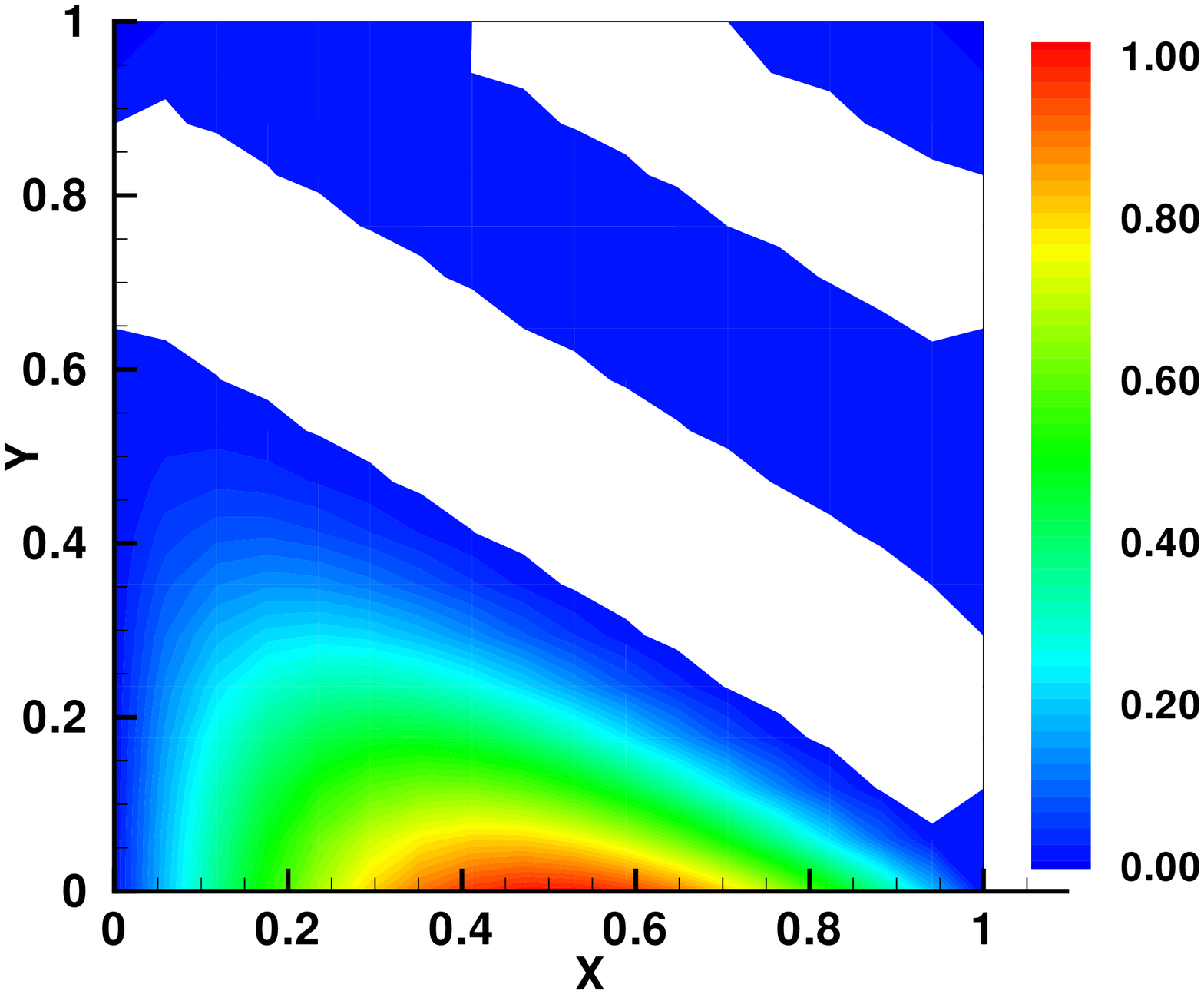}}        
      \subfigure{
      \includegraphics[scale=0.3]{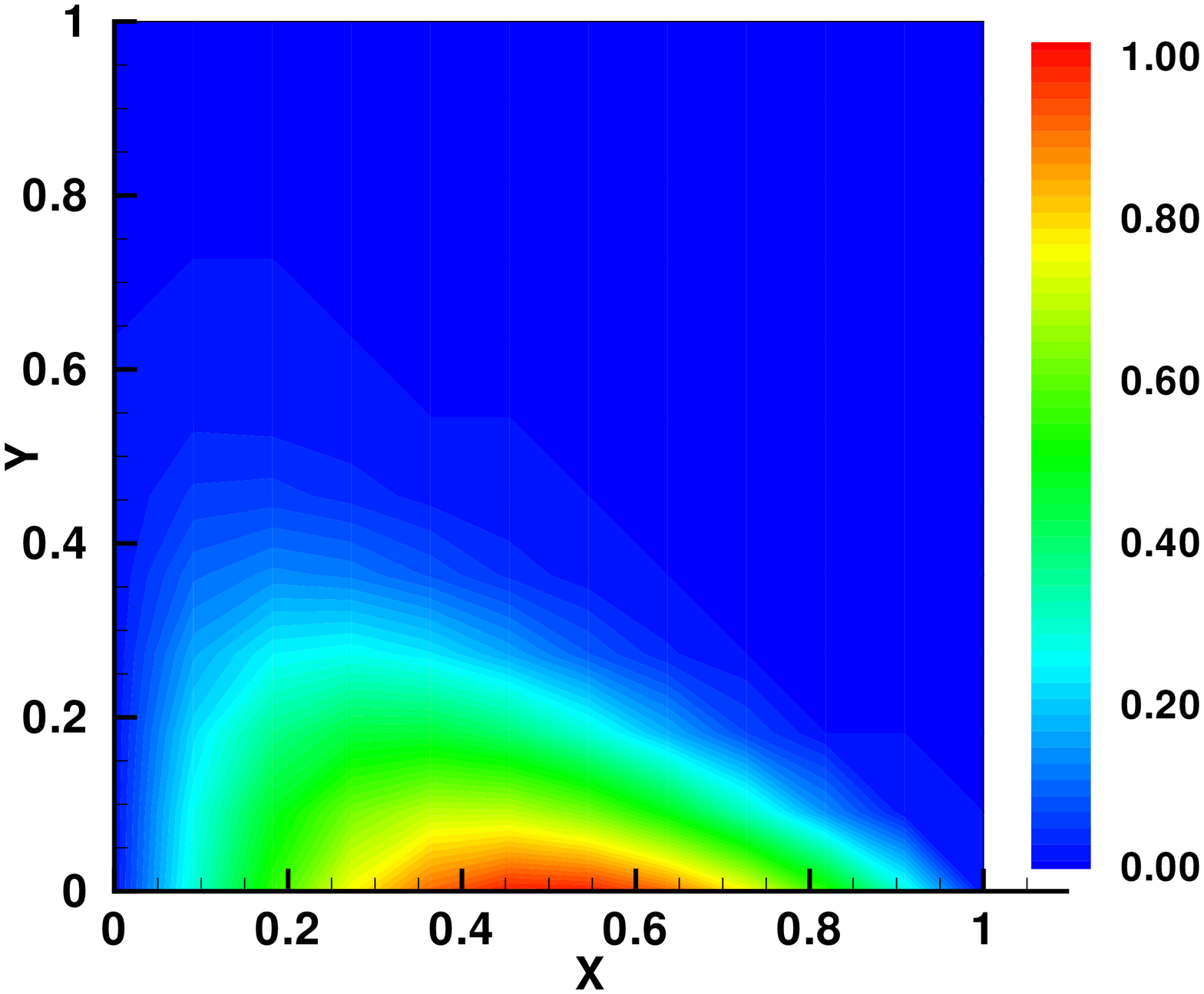}} 
      \subfigure{
      \includegraphics[scale=0.3]{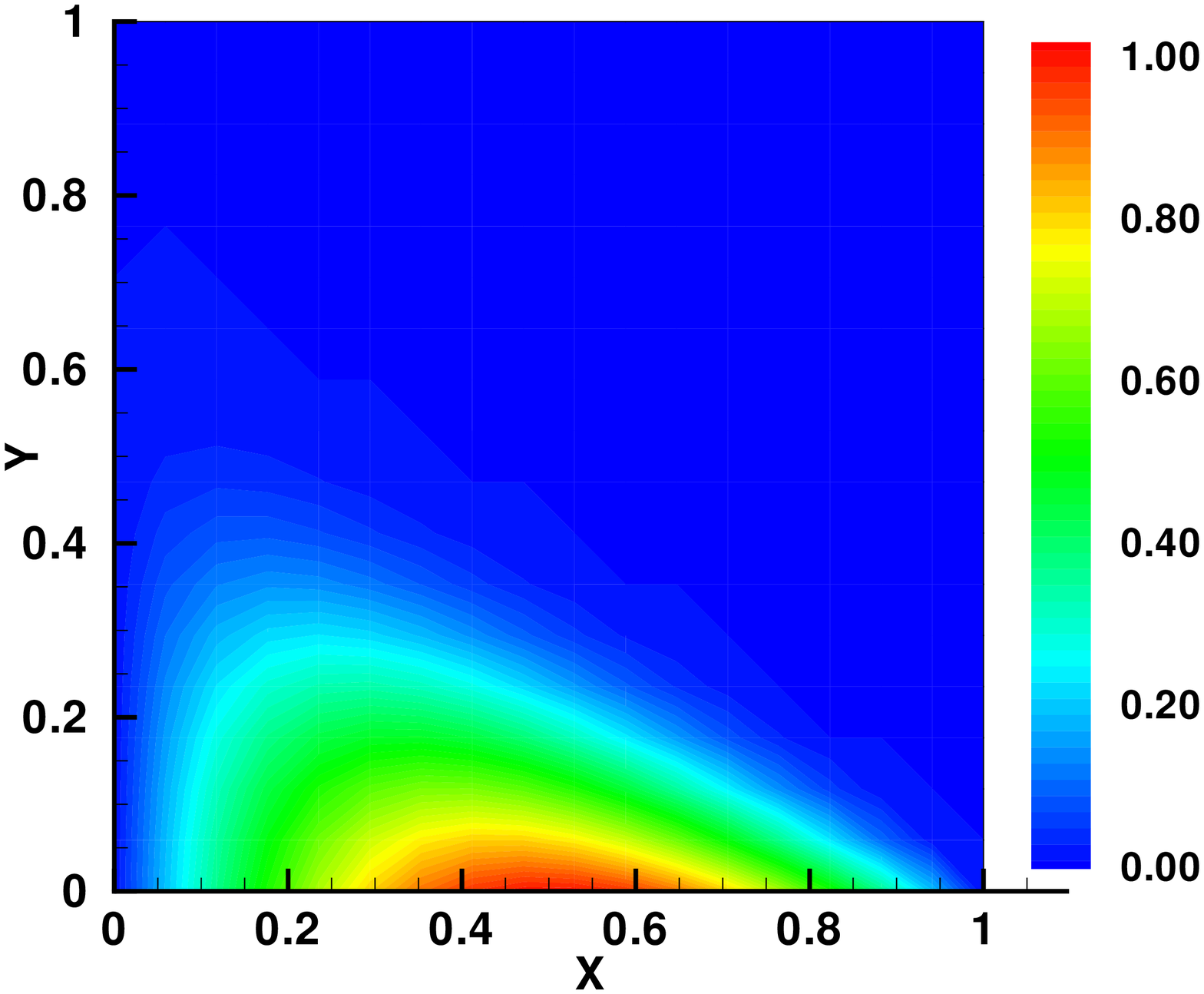}}  
    \caption{Two-dimensional problem with anisotropic medium: The problem is solved using 
      the Galerkin formulation (middle) and the proposed formulation (bottom). The left and right 
      figures are, respectively, using $12 \times 12$ and $18 \times 18$ \emph{four-node 
      quadrilateral meshes}. Regions that have negative concentrations are indicated in white 
      color. Under the Galerkin formulation, $29.17 \%$ (for $12 \times 12$ mesh) and $31.17\%$ 
      of the total number of nodes have negative nodal concentration. The minimum concentrations 
      are $-0.020$ (for $12 \times 12$ mesh) and $-0.017$ (for $18 \times 18$ mesh).}
    \label{Fig:Decay_2D_anisotropic_Q4}
\end{figure}

\begin{figure}[!h]
  \centering
  \includegraphics[scale=0.45]{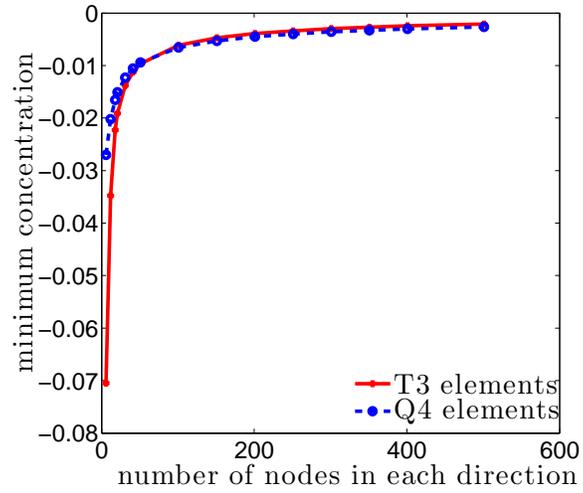}
  \caption{Two-dimensional problem with anisotropic medium: Variation of the minimum 
    concentration with respect to mesh refinement for three-node triangular (T3) and 
    four-node quadrilateral (Q4) meshes.} \label{Fig:Decay_2D_anisotropic_min_conc_mesh_refinement}
\end{figure}

\begin{figure}[htbp]
  \centering
  \subfigure{
    \includegraphics[scale=0.45]{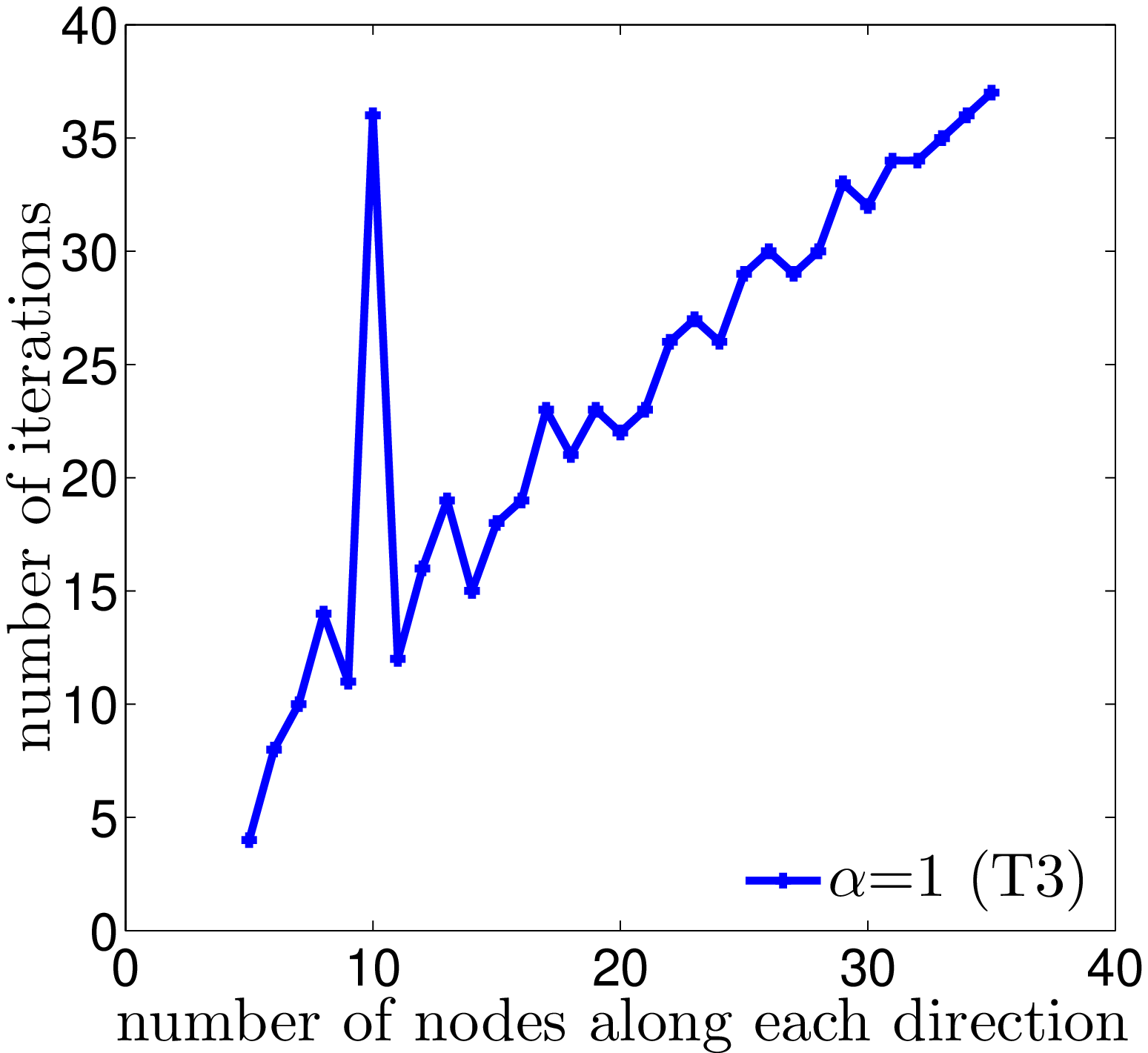}}        
  \subfigure{
    \includegraphics[scale=0.45]{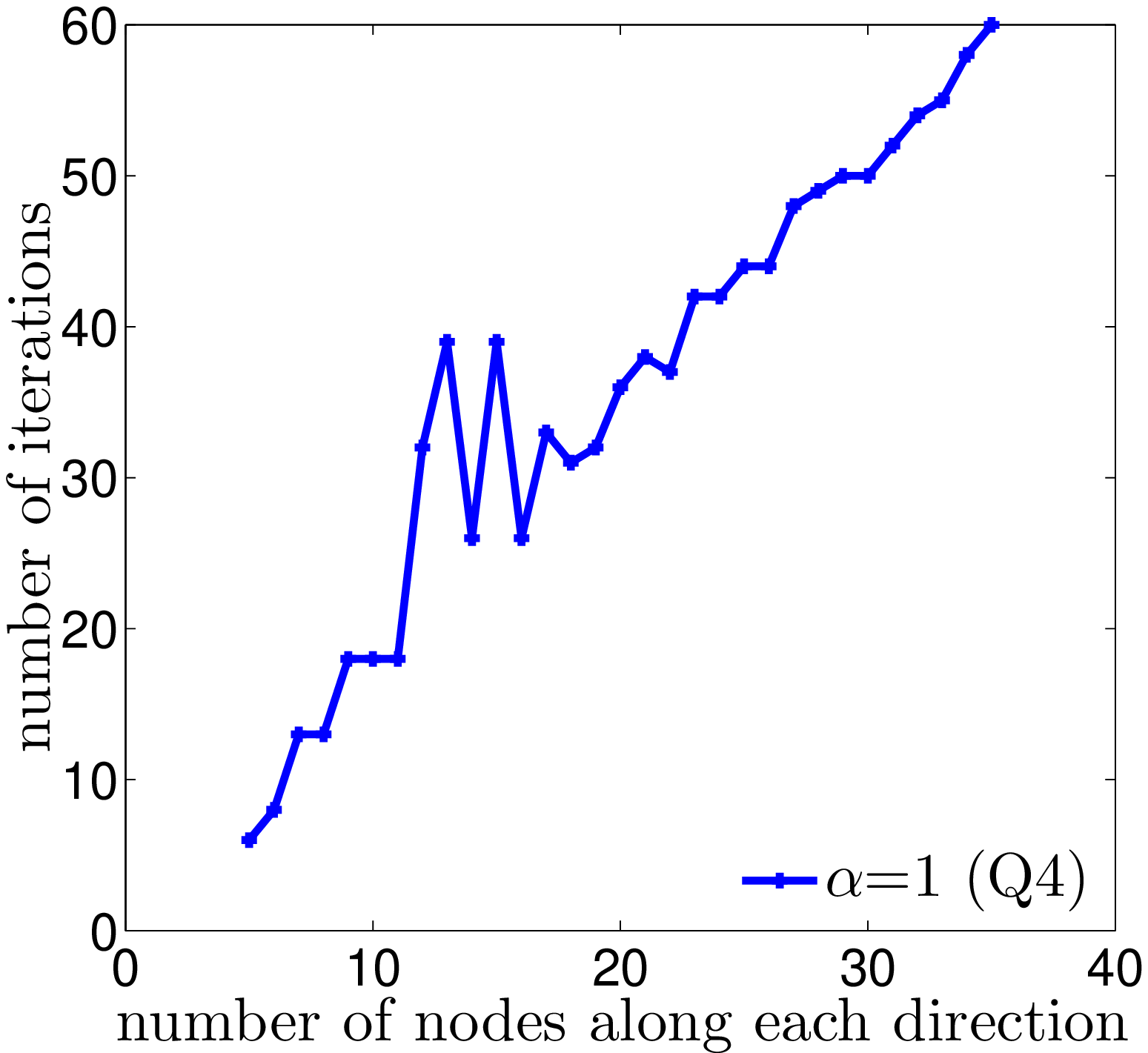}} 
  \subfigure{
    \includegraphics[scale=0.45]{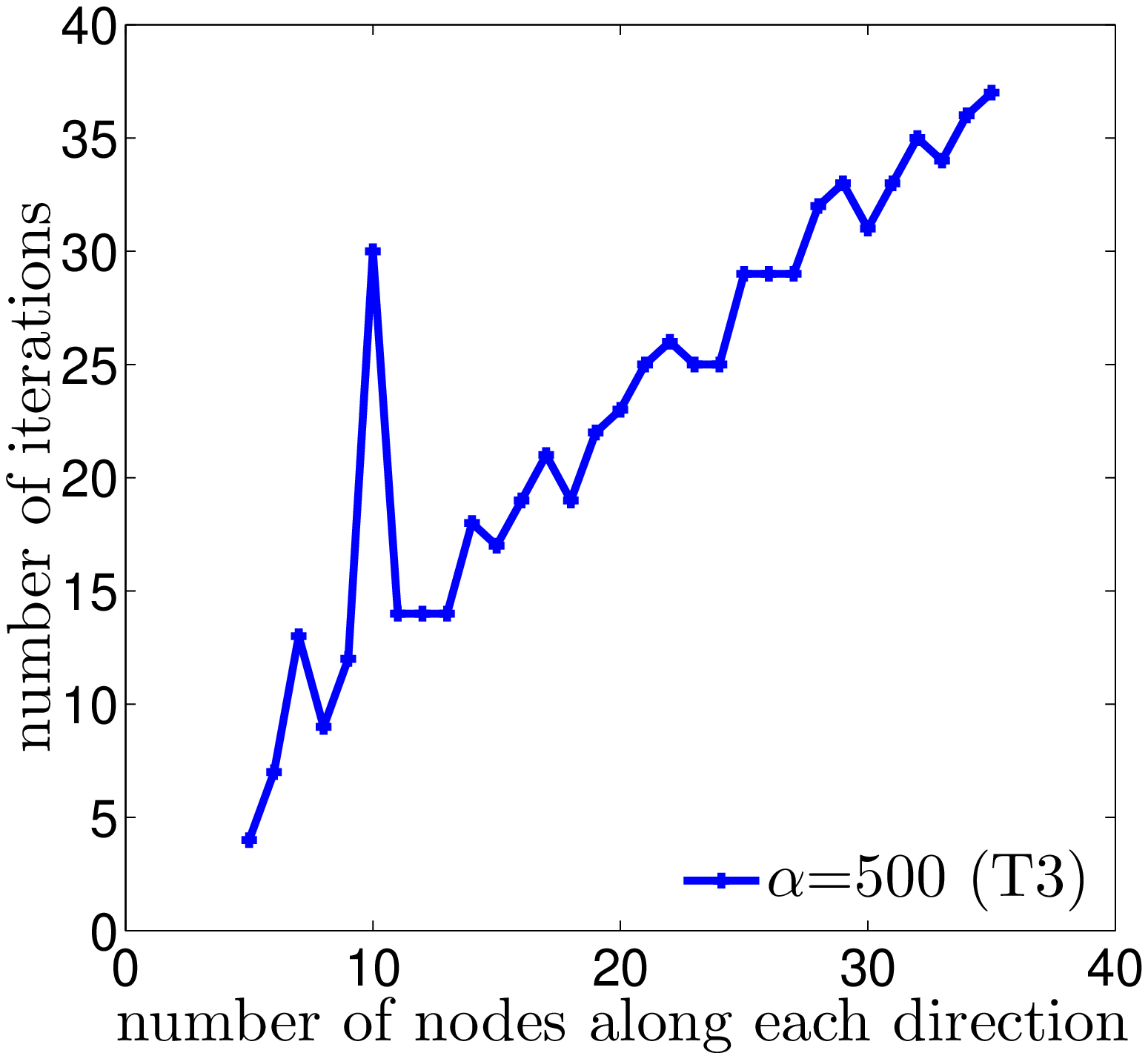}}        
  \subfigure{
    \includegraphics[scale=0.45]{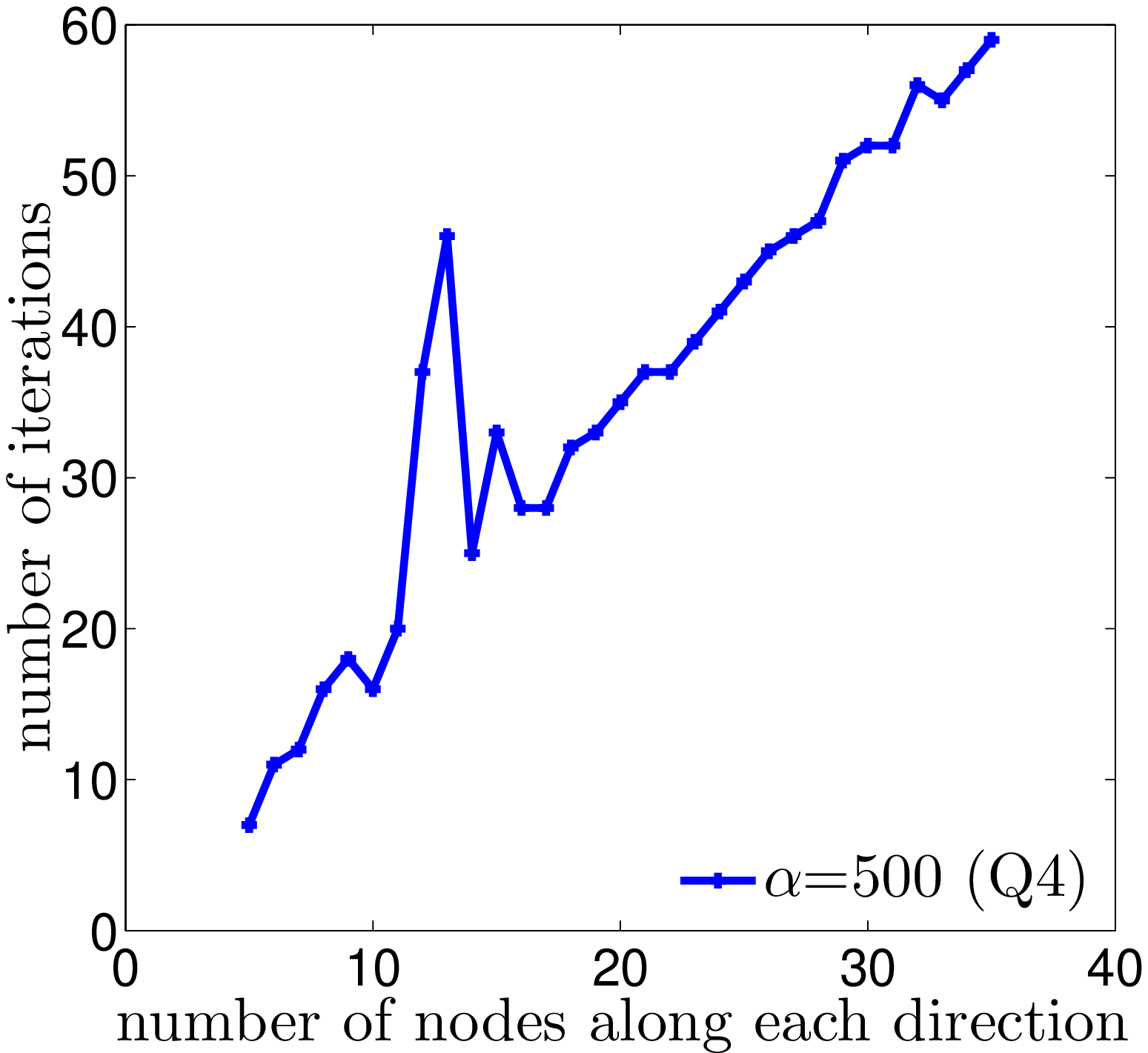}}  
  \caption{Two-dimensional problem with anisotropic medium: This figure shows the number of iterations 
    taken by the active-set strategy under the proposed formulation for two different values of decay 
    coefficient: $\alpha = 1$ (top) and $\alpha = 500$ (bottom). The number of iterations are shown for 
    both three-node triangular (left) and four-node quadrilateral (right) meshes. Equal number of nodes 
    are employed along both x and y directions. Because of the anisotropy, the violation of the maximum 
    principle does not vanish with mesh refinement even for smaller values of decay coefficient (in 
    this case, $\alpha = 1$).} \label{Fig:Decay_anisotropy_iterations_vs_XSeed_alpha_1}
\end{figure}

The problem is solved using two different computational meshes, and the numerical results 
are shown in Figure \ref{Fig:Decay_2D_anisotropic_T3} (for three-node triangular mesh) 
and Figure \ref{Fig:Decay_2D_anisotropic_Q4} (for four-node quadrilateral mesh). 
The amount of the violation of the maximum principle spatially for various mesh 
refinements is illustrated in Tables \ref{Table:Decay_square_plate_T3} and 
\ref{Table:Decay_square_plate_Q4}. 
In Figure \ref{Fig:Decay_2D_anisotropic_min_conc_mesh_refinement}, we have 
shown the variation of minimum concentration with respect to mesh refinement. From 
these figures, it is evident that the negative concentration (which is the indication 
of the violation of the maximum principle) reached constant values for both 
three-node triangular and four-node quadrilateral meshes, and the violation 
existed irrespective of the mesh refinement.
This is the main difference between the violation due to the decay term and the violation 
due to anisotropy. The violation of the maximum principle due to the decay term decreases 
with respect to mesh refinement, and eventually vanishes with mesh refinement. This fact is 
further illustrated in Figure \ref{Fig:Decay_anisotropy_iterations_vs_XSeed_alpha_1}, which 
shows the number of iterations taken by the active-set strategy for various values of $\alpha$. 

\begin{table}[t]
  \centering
  \caption{Two-dimensional problem with anisotropic medium: Violation of maximum 
    principle with respect to mesh refinement using three-node triangular elements. 
    \label{Table:Decay_square_plate_T3}}
  \begin{tabular}{|c|c|c|} \hline
    \textbf{mesh} & \textbf{\# of negative nodes} & \textbf{\% of nodes violated} \\ \hline
    6 $\times$ 6   & 6 & 19.44                \\ \hline
    12 $\times$ 12  & 34 & 27.78          \\ \hline
    18 $\times$ 18  & 90 & 30.86          \\ \hline
    21 $\times$ 21  & 127 & 31.74        \\ \hline
    31 $\times$ 31  & 301 & 33.71        \\ \hline
    41 $\times$ 41  & 546 & 34.56        \\ \hline
    51 $\times$ 51  & 854 & 34.58        \\ \hline
    101 $\times$ 101 & 3500 & 35.53  \\ \hline
  \end{tabular}
\end{table}

\begin{table}[t]
  \centering
  \caption{Two-dimensional problem with anisotropic medium: Violation of maximum 
    principle with respect to mesh refinement using four-node quadrilateral elements. 
    \label{Table:Decay_square_plate_Q4}}
  \begin{tabular}{|c|c|c|} \hline
    \textbf{mesh} & \textbf{\# of negative nodes} & \textbf{\% of nodes violated} \\ \hline
    $6 \times 6$ & 6 & 25.00                  \\ \hline
    $12 \times 12$  & 36 & 29.17          \\ \hline
    $18 \times 18$ & 92 & 31.17           \\ \hline
    $21 \times 21$ & 127 & 31.52         \\ \hline
    31 $\times$ 31 & 291 & 32.98         \\ \hline
    41 $\times$ 41 & 530 & 34.68         \\ \hline
    51 $\times$ 51 & 853 & 35.22         \\ \hline
    101 $\times$ 101 & 3462 & 35.44  \\ \hline
  \end{tabular}
\end{table}

\subsection{Two-dimensional problem with a square hole}
\begin{figure}[!h]
  \centering
  \includegraphics[scale=0.335]{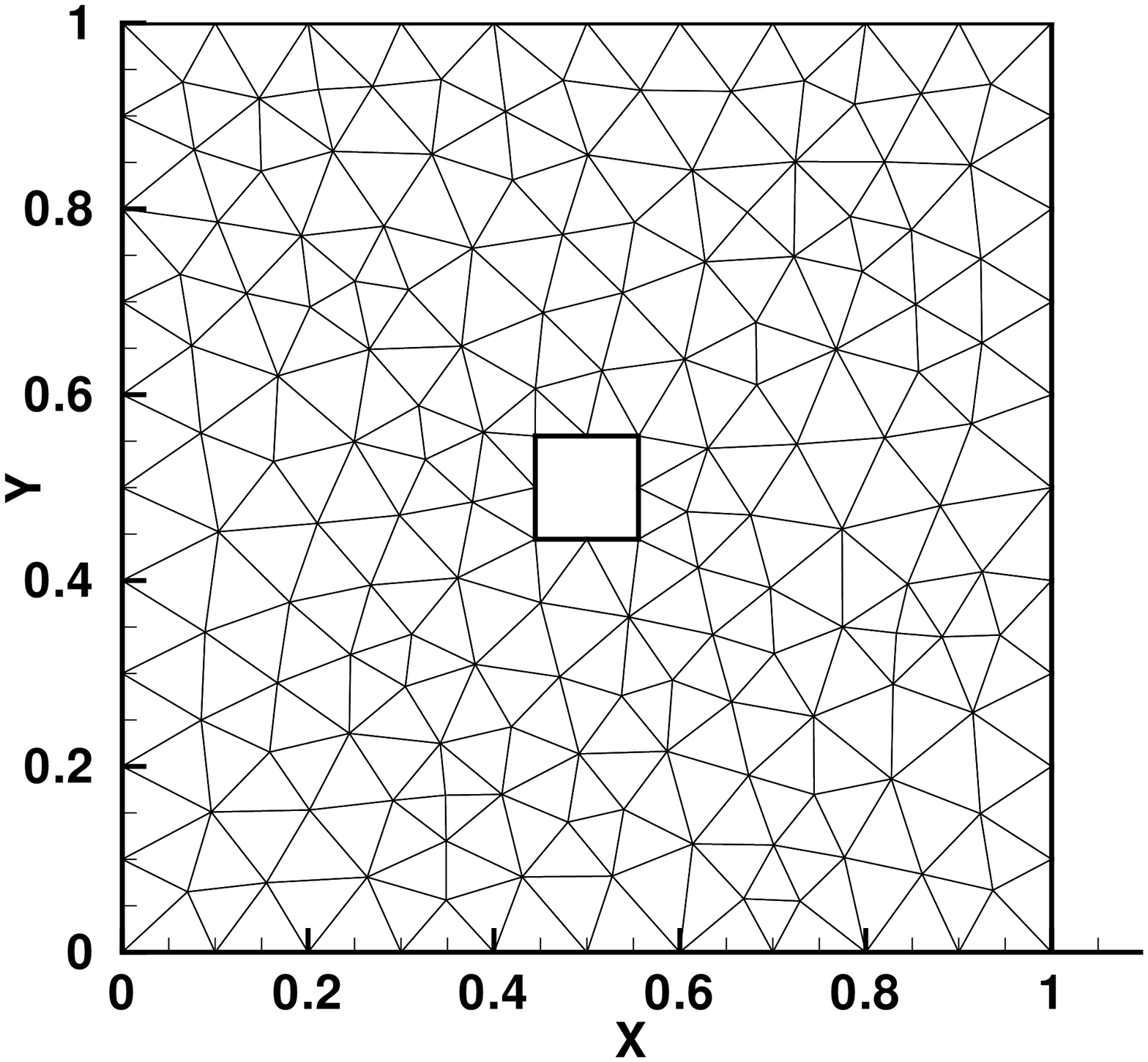}       
  \caption{Two-dimensional problem with a square hole: Computational mesh using 
    three-node triangular finite elements.} \label{Fig:Decay_plate_with_hole_mesh}
\end{figure}

\begin{figure}[!h]
  \centering
  \subfigure{
    \includegraphics[scale=0.335]{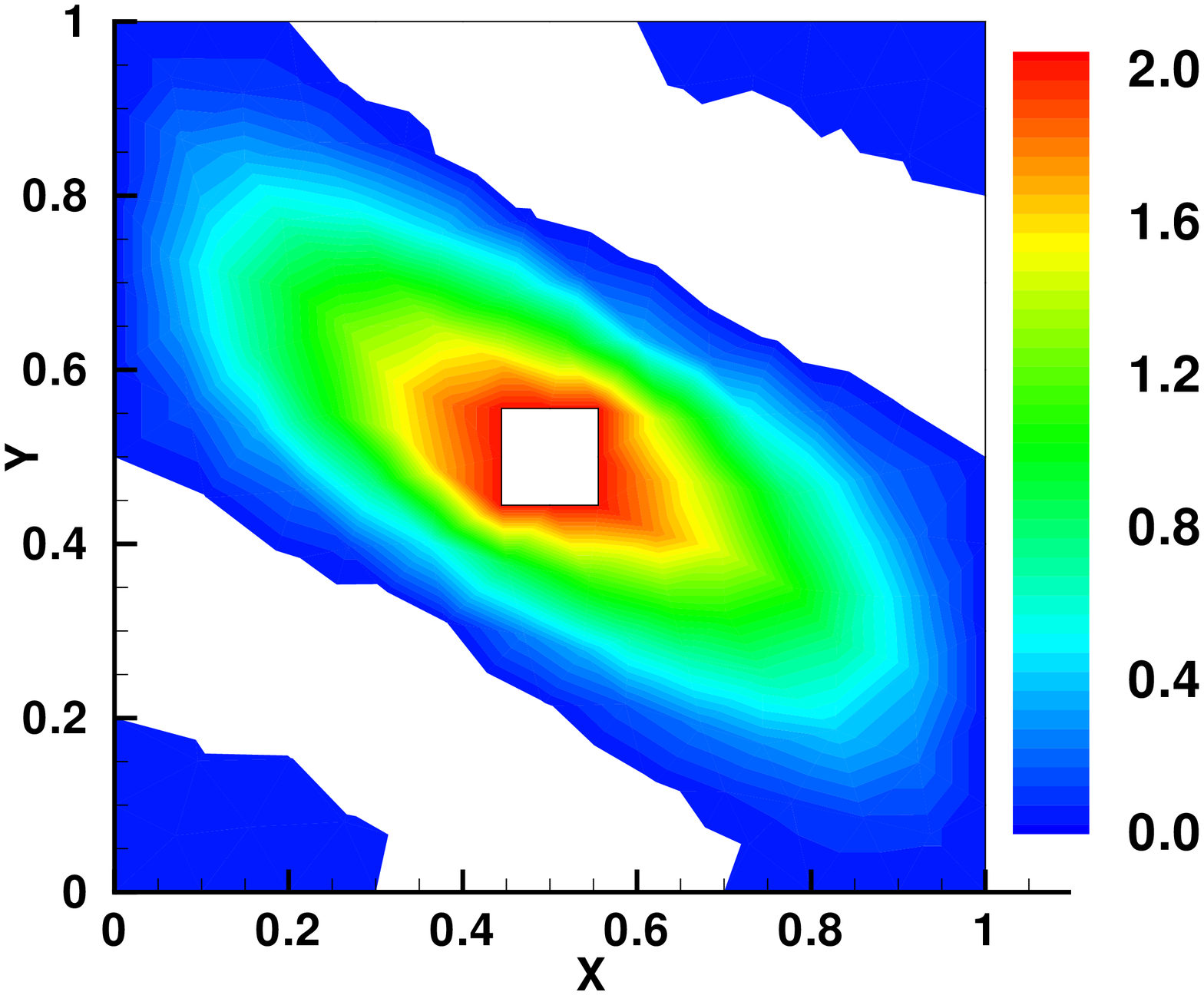}}
  \subfigure{
    \includegraphics[scale=0.335]{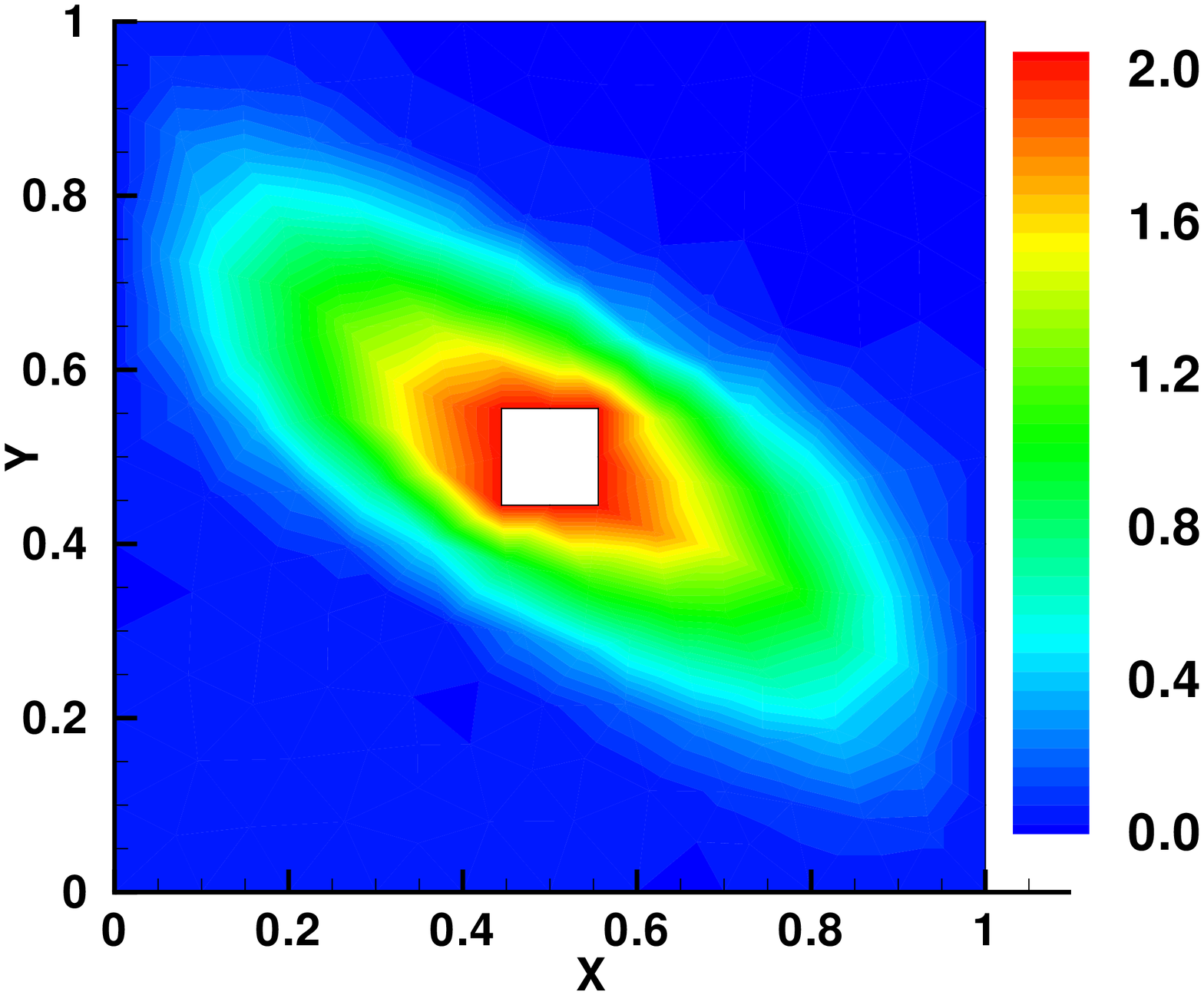}}
  \caption{Two-dimensional problem with a square hole: Contours of the concentration 
    obtained using the Galerkin formulation (left) and the proposed formulation 
    (right) are shown in this figure. Regions that have negative concentrations 
    are indicated in white color. The proposed formulation produced physically 
    meaningful non-negative values for the concentration. Under the Galerkin 
    formulation, approximately $26.92\%$ of the total number of nodes have 
    negative nodal concentrations. The minimum value of the concentration 
    (which occurred inside the domain) is $-0.0916$.} 
  \label{Fig:Decay_plate_with_hole_NR}
\end{figure}

The computational domain is a bi-unit square plate $\Omega := (0, 1) \times (0, 1)$ with a square 
hole of dimension $[4/9, 5/9] \times [4/9, 5/9]$. On the outer boundary we prescribe $c^{\mathrm{p}}
(\mathrm{x},\mathrm{y}) = 0$, and on the inner boundary we prescribe $c^{\mathrm{p}}(\mathrm{x},
\mathrm{y}) = 2$. The forcing function is taken to be zero (that is, $f(\mathrm{x},\mathrm{y}) = 0$). The 
diffusivity tensor is same as the one employed in the previous subsection (see equation 
\eqref{Eqn:Decay_diffusivity_tensor_homogeneous}). The computational mesh employed in this 
numerical simulation is shown in Figure \ref{Fig:Decay_plate_with_hole_mesh}. Numerical results 
obtained using the Galerkin formulation and the proposed formulation are shown in Figure 
\ref{Fig:Decay_plate_with_hole_NR}, and the proposed formulation performed well.

\subsection{Two-dimensional problem with heterogeneous anisotropic medium}
This test problem is similar to the one proposed in Reference \cite{LePotier_CRM_2005_v341_p787}, 
which addressed pure anisotropic diffusion equation. This test problem is considered as 
a good benchmark problem for testing numerical formulations for violation/satisfaction 
of discrete maximum principle. In this test problem, the diffusivity tensor is anisotropic 
and heterogeneous (that is, it varies spatially), and is given by
\begin{align}
  \label{Eqn:Decay_diffusivity_tensor_heterogeneous}
  \mathbb{D}(\mathrm{x},\mathrm{y}) = 
  \left(\begin{array}{cc}
      \mathrm{y}^{2} + \epsilon \mathrm{x}^{2} & -(1-\epsilon) \mathrm{x} \mathrm{y}\\ 
      -(1-\epsilon) \mathrm{x} \mathrm{y} & \mathrm{x}^{2} + \epsilon \mathrm{y}^2 
    \end{array}\right)
\end{align}
where $\epsilon = 10^{-4}$. The domain is a bi-unit square plate: $\Omega = (0,1) \times (0,1)$. 
Homogeneous Dirichlet boundary conditions are prescribed on the entire boundary. The forcing 
function is taken to be $f(\mathrm{x},\mathrm{y}) = 1$ if $(\mathrm{x}, \mathrm{y}) \in 
[3/8,5/8] \times [3/8,5/8]$, and zero otherwise. Since the forcing function is non-negative, 
and homogeneous Dirichlet boundary conditions are prescribed on the whole boundary, from the 
maximum principle we have that the concentration is non-negative in the whole domain (that is, 
$c(\mathbf{x}) \geq 0$ in $\bar{\Omega}$). 

\begin{figure}[htbp]
  \centering
      \subfigure{
      \includegraphics[scale=0.335]{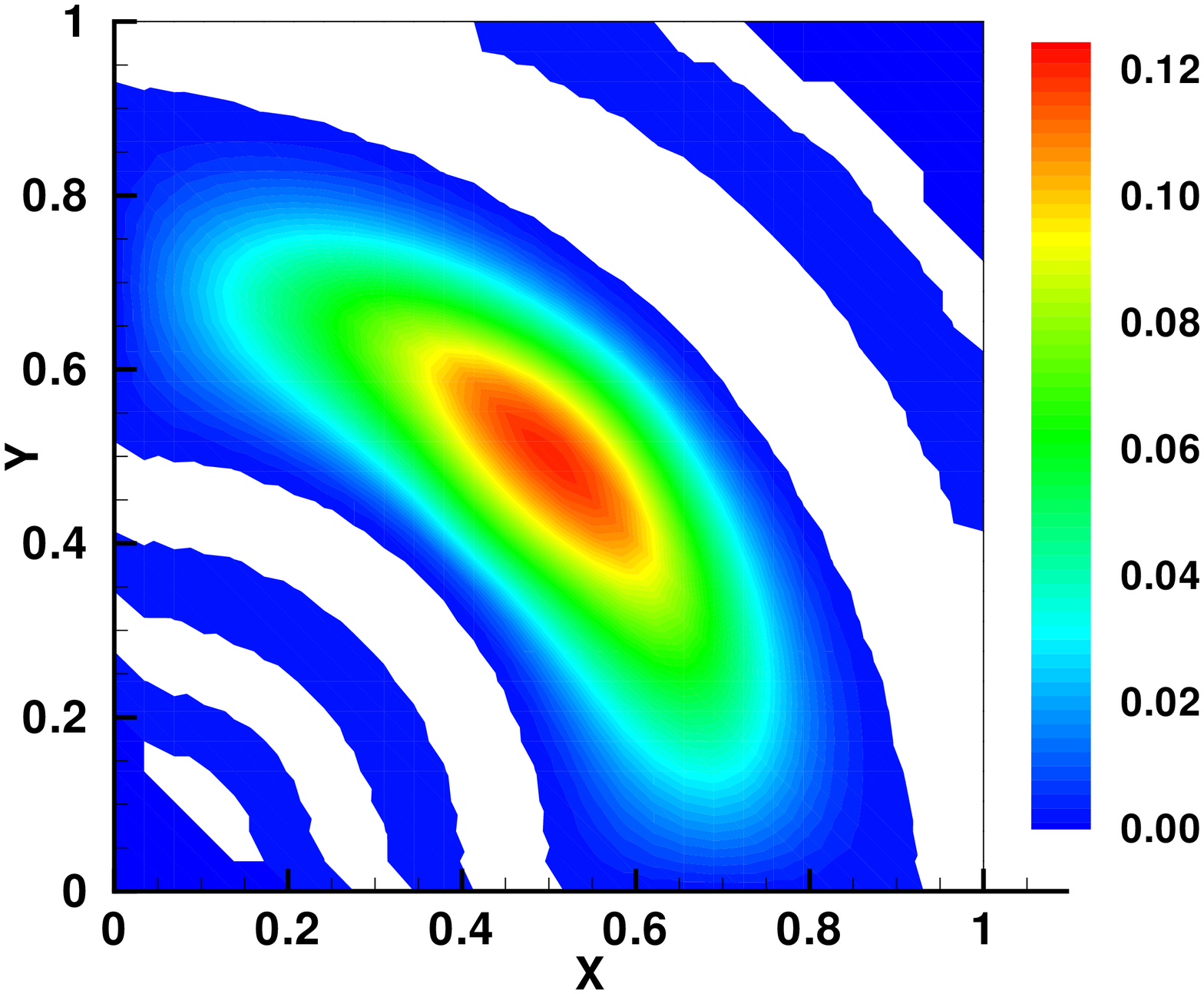}}        
    \subfigure{
      \includegraphics[scale=0.335]{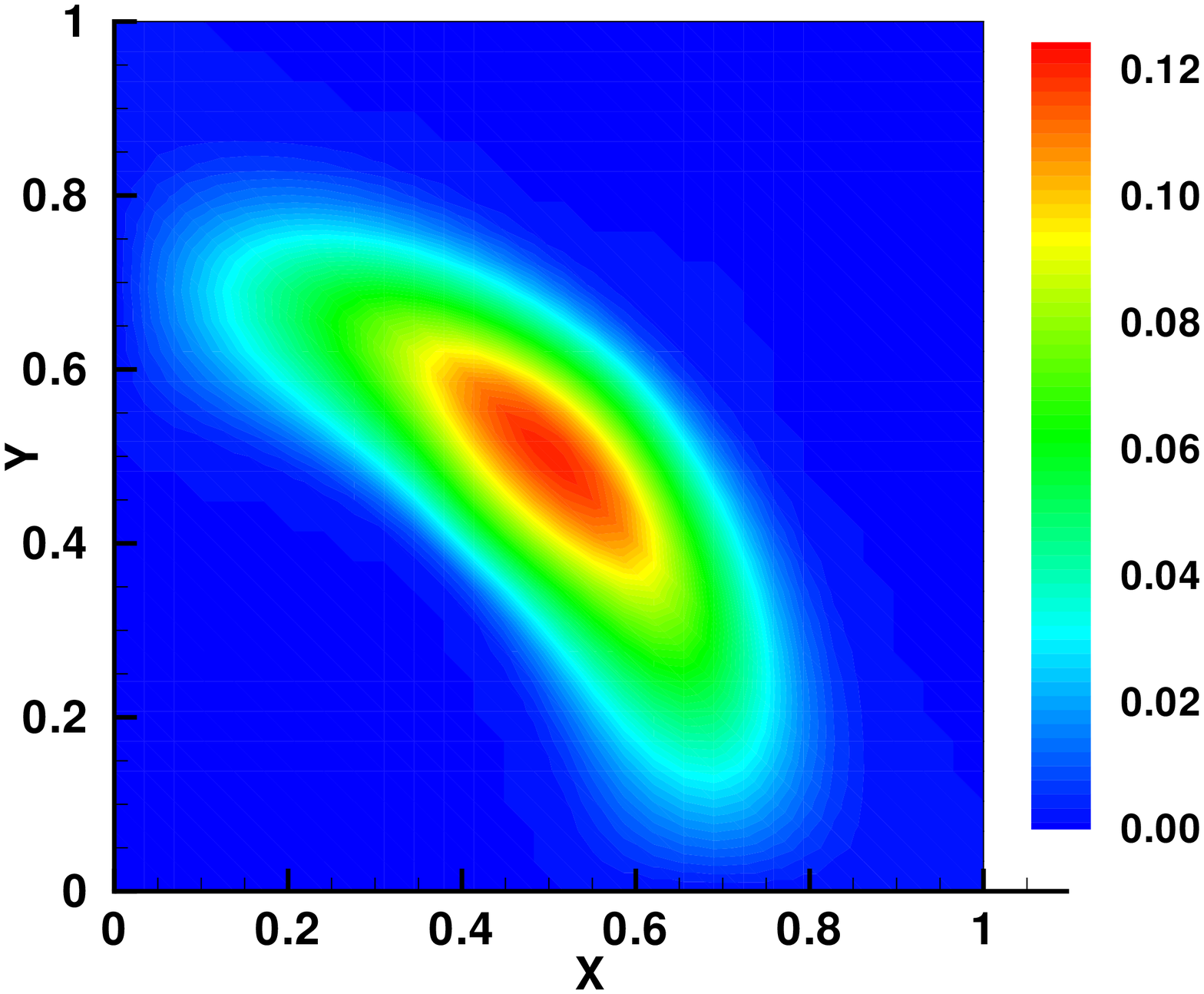}} 
    \caption{Heterogeneous anisotropic medium: This figure shows the concentration obtained 
      using the Galerkin formulation (left) and the proposed formulation (right) for a decay 
      coefficient of $\alpha = 1$. The regions that have negative concentrations are indicated 
      in white color. The proposed formulation produced physically meaningful non-negative 
      values for the concentration. Under the Galerkin formulation, approximately $31.4\%$ 
      of the total number of nodes have negative nodal concentrations. The minimum value of 
      the concentration is $-0.0012$. In the case of anisotropic medium, the violation 
      of the maximum principle will occur even for smaller values of decay coefficient. 
      Moreover, the violation, in general, will not vanish with the mesh refinement.} 
    \label{Fig:Decay_2D_T3_heterogeneous}
\end{figure}

\begin{figure}[htbp]
  \centering
  \includegraphics[scale=0.335]{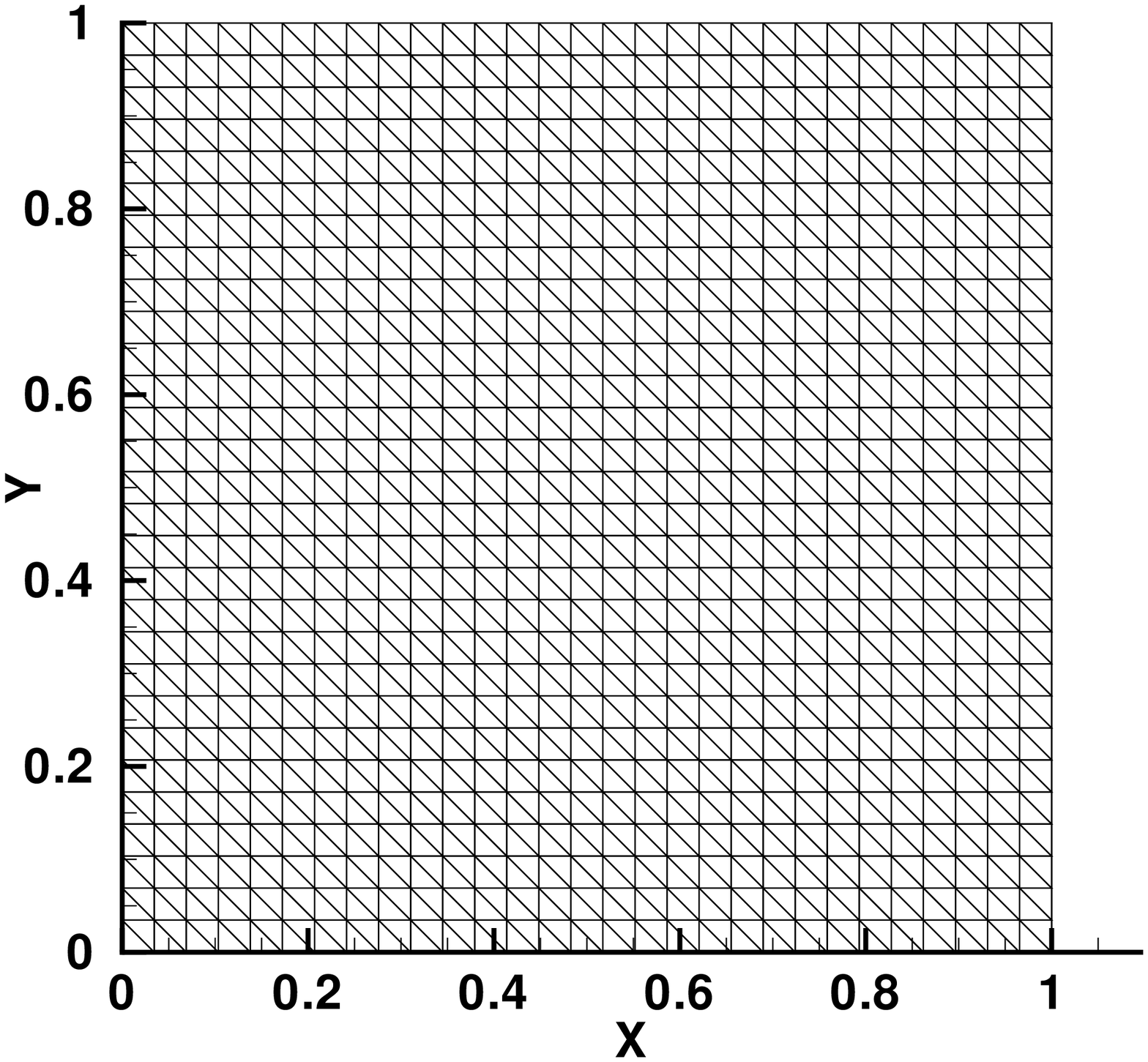}        
  \caption{Heterogeneous anisotropic medium: Computational mesh using three-node 
    triangular finite elements.} \label{Fig:Decay_2D_T3_heterogeneous_mesh}
\end{figure}

\begin{figure}[htbp]
  \centering
  \subfigure{
    \includegraphics[scale=0.45]{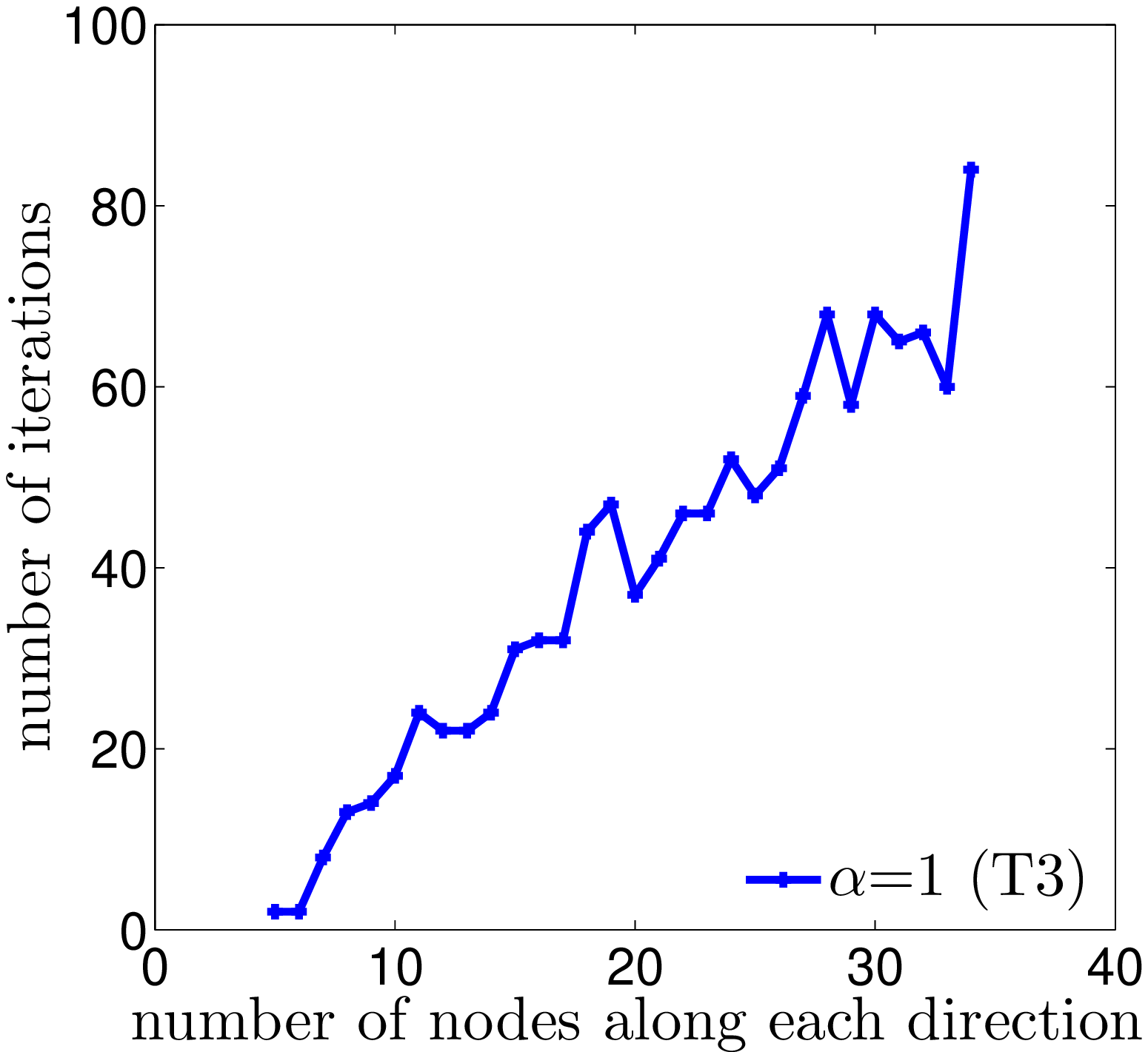}}        
  \subfigure{
    \includegraphics[scale=0.45]{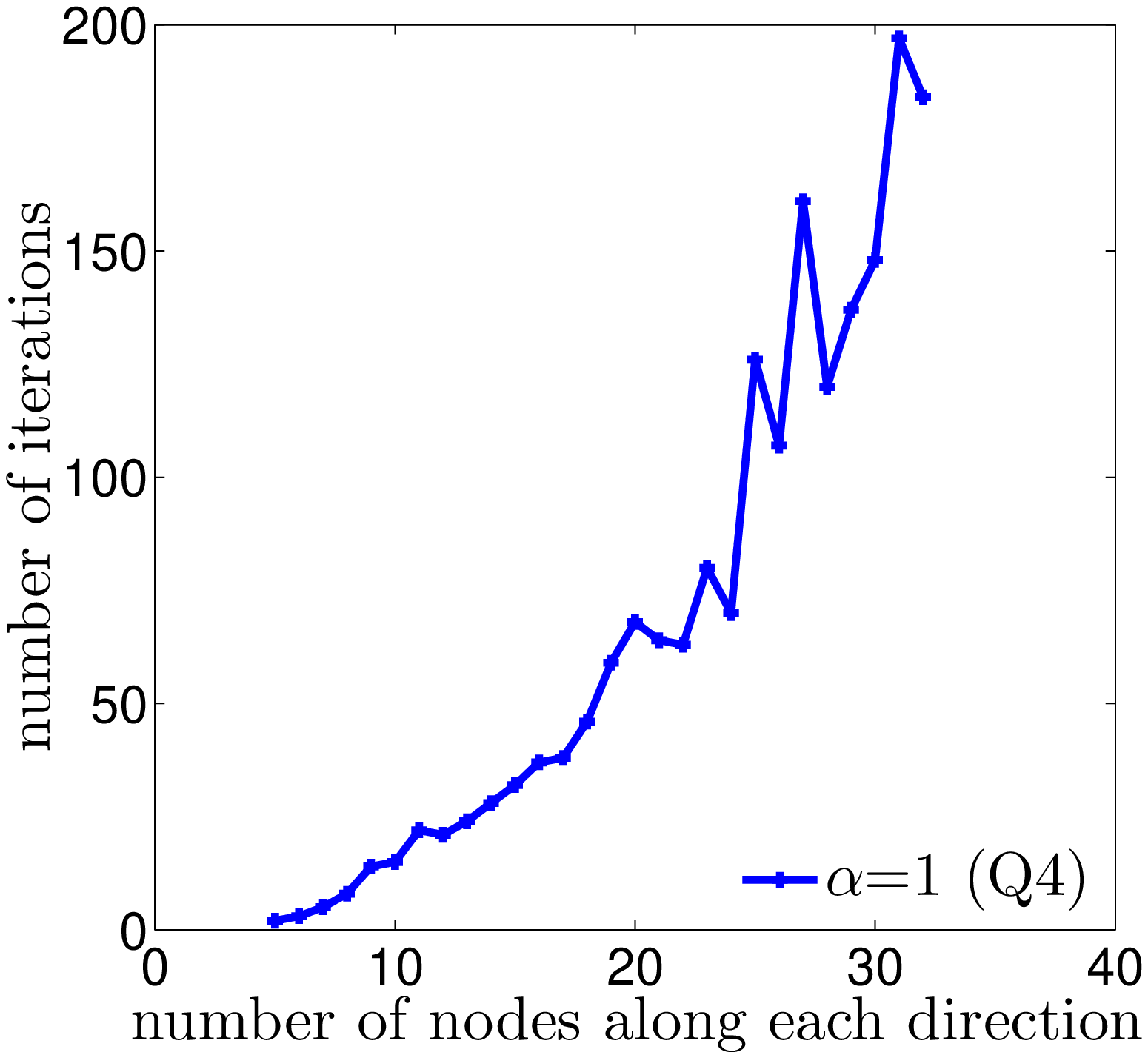}} 
  \subfigure{
    \includegraphics[scale=0.45]{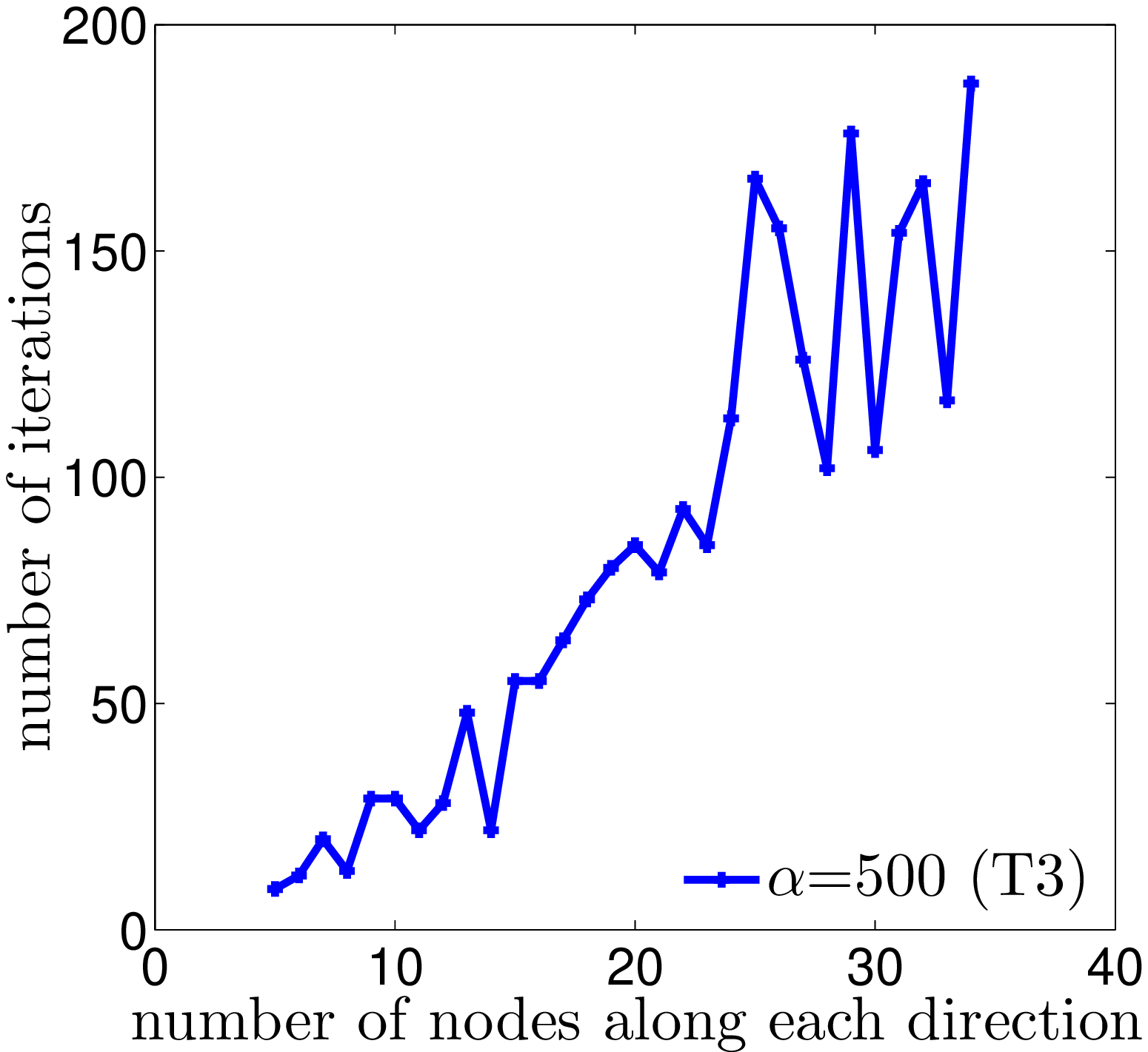}}        
  \subfigure{
    \includegraphics[scale=0.45]{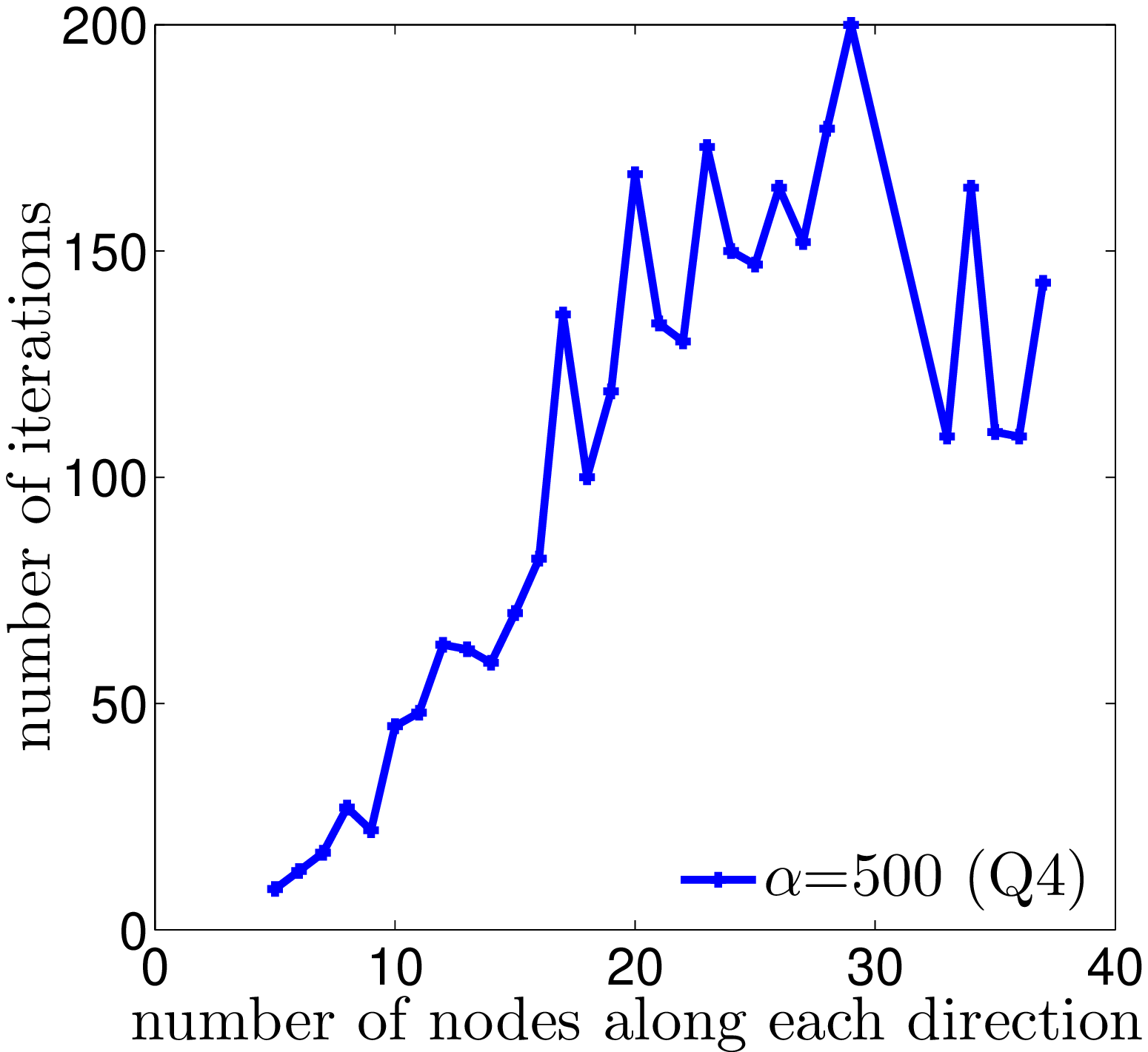}} 
  \caption{Heterogeneous anisotropic medium: This figure shows the number of iterations 
    taken by the active-set strategy under the proposed formulation for two different 
    values of decay coefficient: $\alpha = 1$ (top) and $\alpha = 500$ (bottom). The 
    number of iterations are shown for both three-node triangular (left) and four-node 
    quadrilateral (right) meshes. Equal number of nodes are employed along both x and y 
    directions. Because of the anisotropy and heterogeneity, the violation of the maximum 
    principle does not vanish with mesh refinement even for smaller values of decay 
    coefficient (in this case, $\alpha = 1$).} 
  \label{Fig:Decay_heterogeneity_iterations_vs_XSeed}
\end{figure}

The numerical results that are obtained using the Galerkin formulation and the proposed formulation 
are shown in Figure \ref{Fig:Decay_2D_T3_heterogeneous}. The computational mesh that is employed in 
the numerical simulations is shown in Figure \ref{Fig:Decay_2D_T3_heterogeneous_mesh}. Figure 
\ref{Fig:Decay_heterogeneity_iterations_vs_XSeed} shows the number of iterations taken by the 
active-set strategy under the proposed formulation for both three-node triangular and four-node 
quadrilateral meshes. In the case of heterogeneous anisotropic medium, violation of discrete 
maximum principle occurs under the Galerkin formulation even for lower values of decay coefficient 
(e.g., even for $\alpha = 1$), and the violation does not vanish even with mesh refinement.

\begin{remark}
  The proposed methodology works even for low-order three-dimensional finite 
  elements like four-node tetrahedron element, eight-node brick element, and 
  six-node wedge element. Herein, a three-dimensional problem is not solved 
  as there are no computational challenges other than standard book-keeping.
\end{remark}

\section{CONCLUSIONS}
\label{Sec:Decay_Conclusions}
In this paper, we have presented a methodology for enforcing the 
non-negative constraint and maximum principles for anisotropic 
diffusion with decay. The proposed method is obtained by adding 
constraints to the variational structure of the classical Galerkin 
formulation, and can handle general computational grids with low-order 
finite elements. The resulting equations form a convex quadratic 
programming problem, and are solved by employing the active-set 
strategy. Numerical experiments have shown that the rates of 
convergence with respect to mesh refinement in $L^{2}$-norm 
and $H^{1}$-seminorm are about the same as for the original 
linear finite element method. Various representative numerical 
examples are presented to illustrate the good performance of 
the proposed formulation. 


\section*{ACKNOWLEDGMENTS}
The research reported herein was supported by Texas Engineering Experiment Station 
(TEES). This support is gratefully acknowledged. The opinions expressed in this paper 
are those of the authors and do not necessarily reflect that of the sponsor. 

\bibliographystyle{plain}
\bibliography{Master_References,Books}

\clearpage 
\newpage

\clearpage
\newpage

\clearpage
\newpage

\clearpage
\newpage

\end{document}